\documentclass{article}

\usepackage[final, nonatbib]{neurips_2023}  

\usepackage[utf8]{inputenc} 
\usepackage[T1]{fontenc}    
\usepackage{hyperref}       
\usepackage{url}            
\usepackage{booktabs}       
\usepackage{amsfonts}       
\usepackage{nicefrac}       
\usepackage{microtype}      
\usepackage{xcolor}         

\usepackage[
    backend=biber,
    style=alphabetic, maxalphanames=5, 
    maxbibnames=100, maxsortnames=100, 
    sorting=nyt,
    natbib=true,
    url=false,doi=false,isbn=false,eprint=false
]{biblatex}
\usepackage{graphicx}
\usepackage{subcaption}
\graphicspath{ {./images/} }
\usepackage{enumitem}
\usepackage[normalem]{ulem}
\usepackage[toc,page]{appendix}

\hypersetup{
    bookmarksnumbered,
    colorlinks=true,
    linkcolor = blue,
    urlcolor  = blue,
    citecolor = teal
}

\newcommand\fnsurl[1]{{\footnotesize\url{#1}}}

\usepackage{myquicksetup}
\numberwithin{definition}{section}
\numberwithin{theorem}{section}
\numberwithin{corollary}{section}
\numberwithin{proposition}{section}
\numberwithin{lemma}{section}
\numberwithin{claim}{section}
\numberwithin{fact}{section}
\numberwithin{remark}{section}
\numberwithin{example}{section}
\numberwithin{equation}{section}
\mathtoolsset{showonlyrefs} 

\newtheorem{assumption}{Assumption}

\newcommand{\figureref}[1]{\autoref{#1}}
\newcommand{\theoremref}[1]{\autoref{#1}}

\newcommand{\propositionref}[1]{\autoref{#1}}

\newcommand{\lemmaref}[1]{\autoref{#1}}
\newcommand{\definitionref}[1]{\autoref{#1}}
\newcommand{\sectionref}[1]{\autoref{#1}}
\newcommand{\appendixref}[1]{\autoref{#1}}
\newcommand{\tableref}[1]{\autoref{#1}}

\usepackage{mylivemacros}

\usepackage{ccnondegen} 
\usepackage{centernot}
\usepackage{tensor}
\usepackage{makecell}
\usepackage{thmtools}  
\usepackage{thm-restate}

\newif\ifextended  
\extendedtrue

\usepackage{subfiles} 

\usepackage{times}

\hypersetup{
    pdftitle={Local Convergence of Gradient Methods for Min-Max Games: Partial Curvature Generically Suffices},
    pdfauthor={Guillaume Wang, Lénaïc Chizat},
    pdfkeywords={%
        Gradient methods,
        min-max optimization,
        spectral analysis,
        last-iterate convergence
    }
}

\title
{Local Convergence of Gradient Methods for Min-Max Games: Partial Curvature Generically Suffices}

\author{%
    Guillaume Wang \\
    Institute of Mathematics\\
    École polytechnique fédérale de Lausanne\\
    Station Z, CH-1015 Lausanne \\
    \texttt{guillaume.wang@epfl.ch} \\
    \And
    Lénaïc Chizat \\
    Institute of Mathematics\\
    École polytechnique fédérale de Lausanne\\
    Station Z, CH-1015 Lausanne \\
    \texttt{lenaic.chizat@epfl.ch}
}

\begin{document}

\maketitle


\begin{abstract}
    We study the convergence to local Nash equilibria of gradient methods for two-player zero-sum differentiable games.
    It is well-known that
    such dynamics converge locally when $S \succ 0$ and may diverge when $S=0$, where $S\succeq 0$ is the symmetric part of the Jacobian at equilibrium that accounts for the ``potential'' component of the game.
    We show that these dynamics also converge as soon as $S$ is nonzero (\emph{partial curvature}) and the eigenvectors of the antisymmetric part $A$ are in general position
    with respect to the kernel of $S$.
    We then study the convergence rates when $S \ll A$ and prove that they typically depend on the \emph{average} of the eigenvalues of $S$, instead of the minimum as an analogy with minimization problems would suggest.
    To illustrate our results, we consider the problem of computing mixed Nash equilibria of continuous games.
    We show that, thanks to partial curvature, conic particle methods -- which optimize over both weights and supports of the mixed strategies -- generically converge faster than fixed-support methods.
    For min-max games, it is thus beneficial to add degrees of freedom ``with curvature'': this can be interpreted as yet another benefit of over-parameterization.
\end{abstract}



\section{Introduction} \label{sec:intro}

Min-max optimization is notoriously subtler than minimization, even in convex-concave settings.
While many of the proof techniques for minimization have a natural equivalent in the min-max world, some common intuitions fail to transfer.
The picture is clear for strongly convex-strongly concave (SC-SC) min-max games: all the classical gradient methods converge exponentially (for small enough step-sizes) with worst-case convergence rates dependent on the strong convexity and strong concavity parameters $\mu_x$, $\mu_y$.
But, most famously perhaps, for bilinear min-max games the last iterate of simultaneous Gradient Descent-Ascent~(GDA)
diverges, while the Proximal Point~(PP) method converges, and alternating GDA and the continuous-time limit of all of these algorithms -- Gradient Flow~(GF) -- exhibit a cycling behavior.

Based on those two extreme cases, it still seems that part of our intuition from minimization, where
the last-iterate convergence rate is indeed determined by the strong convexity parameter, is preserved.
In this paper, we argue that this intuition is in fact overly pessimistic, and that 
gradient methods behave in general more favorably in the min-max setting than in the minimization setting.

\begin{figure}[t]
    \begin{subfigure}[t]{\textwidth}
        \centering
    	\includegraphics[width=0.24\textwidth]{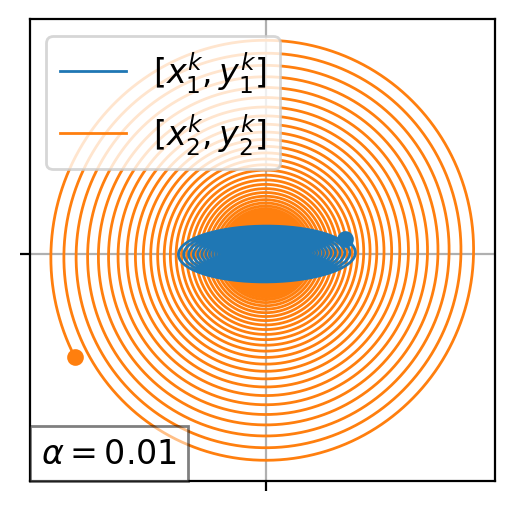}
    	\includegraphics[width=0.24\textwidth]{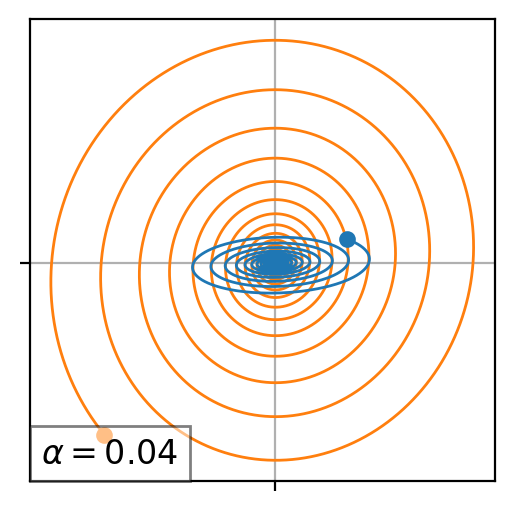}
    	\includegraphics[width=0.24\textwidth]{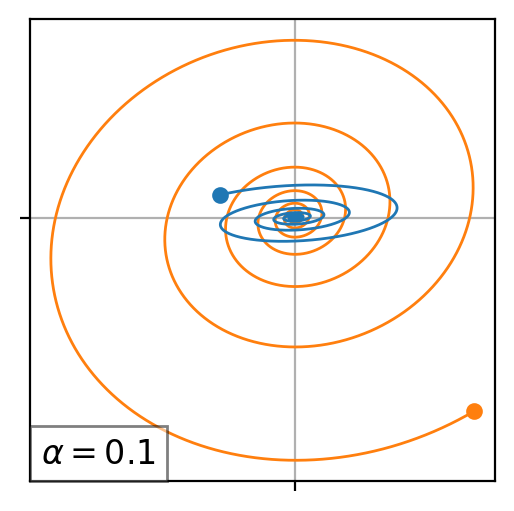}
    	\includegraphics[width=0.24\textwidth]{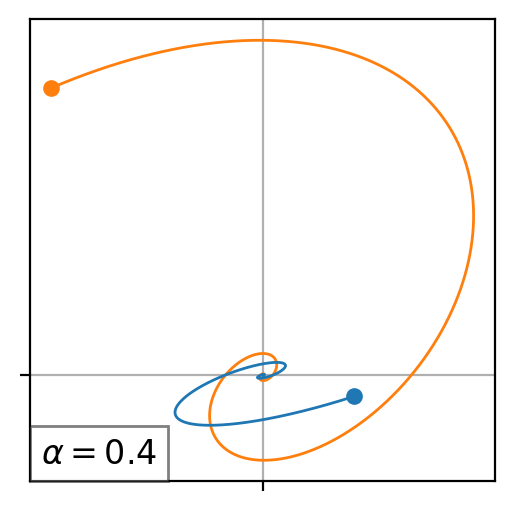}
        \\[-0.5em]
        \caption{Local convergence of the GF iterates for a fixed draw of $P \in \RR^{2 \times 2}$ and various values of $\alpha$. 
        Only the final phase of the dynamics is shown, so here we see the iterates evolve along the ``dominant'' eigenspace of~$M$ (the subspace along which convergence is the slowest).}
        \label{fig:intro:regbilinear:iterates}
    \end{subfigure}
    \\
    \begin{subfigure}[t]{0.48\textwidth}
        \centering
    	\includegraphics[width=\textwidth]{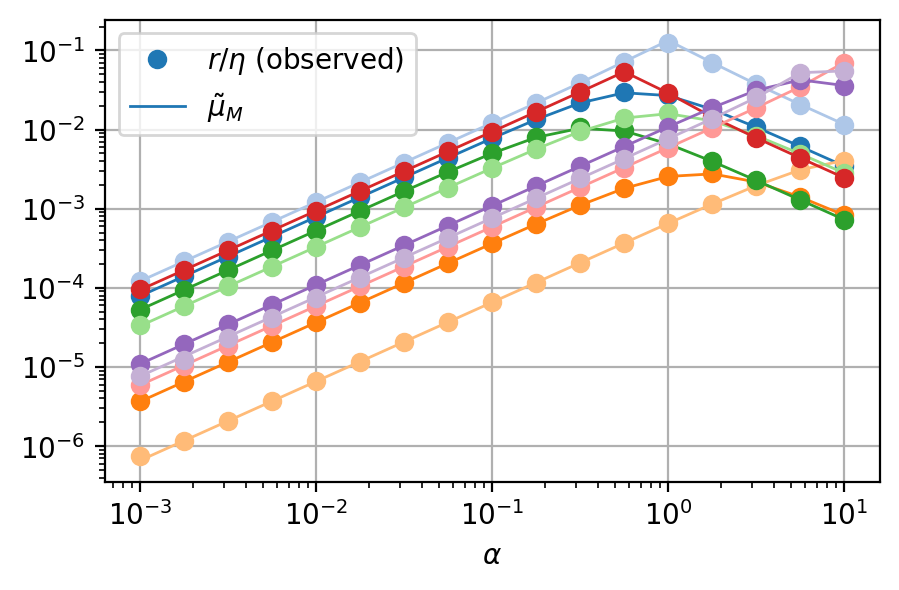}
        \\[-0.5em]
        \caption{Observed and predicted ($\tmu_M$) normalized convergence rate $r/\eta$ of GDA with a small step-size $\eta$, i.e., $\norm{x^k}+\norm{y^k} = \Theta((1-r)^k)$, vs.\ regularization strength $\alpha$ (the higher the faster). Each color represents one draw of $P \in \RR^{2 \times 2}$.}
        \label{fig:intro:regbilinear:rates}
    \end{subfigure}
    \hfill
    \begin{subfigure}[t]{0.48\textwidth}
        \centering
    	\includegraphics[width=\textwidth]{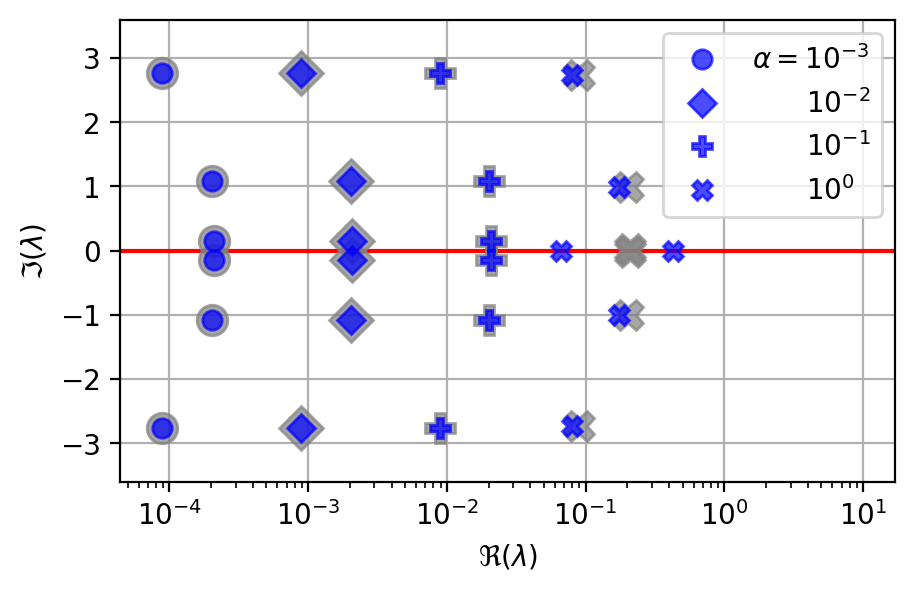}
        \\[-0.5em]
        \caption{Spectrum of the Jacobian $M$
        for a fixed draw of $P \in \RR^{3 \times 3}$ (green) and its approximation by Equation~\eqref{eq:quantitative:ord1} (gray)}
        \label{fig:intro:regbilinear:spectrum}
    \end{subfigure}
    \caption{Convergence of gradient methods on a random bilinear game regularized by $\frac{\alpha}{2} x_1^2$ for small step-sizes. 
    The fact that $\tmu_M$, and so $r/\eta$, scale linearly with $\alpha$ (for small $\alpha$) is explained by \propositionref{prop:quantitative:ord1}.}
    \label{fig:intro:regbilinear}
\end{figure}

Let us first look at a basic example. Consider a bilinear min-max game regularized by a quadratic term \emph{only in one scalar variable}:
\begin{equation}
    \min_{x \in \RR^d} \max_{y \in \RR^d} x^\top P y + \frac{\alpha}{2} x_1^2.
\end{equation}
Here there is no strong convexity in either player (so $\mu_x=\mu_y=0$), yet as observed on \figureref{fig:intro:regbilinear:iterates},~\figureref{fig:intro:regbilinear:rates},
when $P\in \RR^{d\times d}$ is random with independent standard normal entries,
GDA (with a small step-size) typically converges at an exponential rate
that scales linearly with $\alpha$.%
\footnote{
    Julia code for the numerical experiments is available at 
    {\fnsurl{https://github.com/guillaumew16/local_cvgce_minmax}}.
}

As we will see, this phenomenon is a 
consequence of the existing theory for general smooth min-max games 
\begin{equation}\label{eq:intro:smooth-minmax}
    \min_{x \in \RR^n} \max_{y \in \RR^m} f(x, y)
\end{equation}
with a local Nash equilibrium (NE), or local saddle point,
$z^* = (x^*, y^*)$.
It is well-known from the dynamical systems literature that GF converges locally to $z^*$ if
\begin{equation}\label{eq:intro:muM}
    \qquad
    \tilde \mu_M \coloneqq \min_{\lambda \in \spectrum(M)} \Re(\lambda) > 0,
    \qquad
    \text{with}~~
    M = \begin{bmatrix}
        \nabla_{xx}^2 f & \nabla_{xy}^2 f \\
        -\nabla_{xy}^2 f^\top & -\nabla_{yy}^2 f
    \end{bmatrix} (z^*)
    \in \RR^{(n+m)\times (n+m)}
\end{equation}
the Jacobian of the skewed gradient field at $z^*$,
and where $\spectrum(\cdot)$ denotes the spectrum, i.e., the set of eigenvalues.
For our starting basic example, this quantity%
\footnote{In numerical analysis and stability theory, the quantity $-\tilde \mu_{M}$ is known as the \emph{spectral abscissa} of $-M$.} 
can be visualized in \figureref{fig:intro:regbilinear:spectrum}. 
Moreover, $\tmu_M$ also
characterizes the convergence behavior of gradient methods 
\ifextended%
    (GDA, PP, etc.) 
\fi
in the leading order in the step-size $\eta$: they all converge to $z^*$ with the rate $\eta \tmu_M + O(\eta^2)$, see 
\ifextended%
    \appendixref{apx:algo_rates} 
\else%
    the Appendix
\fi
for a review.%
\footnote{
    More precisely, 
    gradient methods converge locally (for small enough step-sizes)
    if $\min_{\lambda \in \spectrum(M)} \Re(\lambda) > 0$, and only~if $\min_{\lambda \in \spectrum(M) \setminus \{0\}} \Re(\lambda) \geq 0$. Throughout this paper we will assume for simplicity $M$ invertible so that those two quantities are equal.
}

While the condition $\tmu_M > 0$ is very general and tight, it is not obvious how to control or interpret it. Our purpose with this paper is to explore what this condition entails and to emphasize some of its surprising consequences. 

The rest of the paper is organized as follows. In \autoref{sec:quadr} we show that, generically, ${ \tilde \mu_M>0 }$ as long as the problem has \emph{partial curvature}, i.e.\ the diagonal blocks of $M$ are non-zero. 
In \autoref{sec:quantitative}, we study more precisely the case of games 
where the ``interaction'' component dominates the ``potential'' component, i.e.,
when $\nabla_{xx}^2 f(z^*), \nabla_{yy}^2 f(z^*) \ll \nabla_{xy}^2 f(z^*)$. For such games, in a certain random setting,  we show that $\tilde \mu_M$ 
scales as
the \emph{average} of the eigenvalues of the potential part, instead of the minimum as an analogy with minimization problems would suggest.
In \autoref{sec:particleMNE} we consider the computation of mixed Nash equilibria of continuous games using particle methods, a setting where the convergence under partial curvature has a striking consequence.
By optimizing over both the weights and supports of the mixed strategies,
one obtains dynamics that generically converge faster than fixed-support methods, even when the latter are using the optimal supports.

For ease of exposition, in most of the paper we focus on the local convergence of GF. We discuss the convergence of discrete time algorithms in \autoref{subsec:quantitative:DT_algos}.

\subsection{Related work} \label{subsec:intro:related_work}

Throughout this paper, by ``convergence'' we mean convergence of the last iterate to a (local) NE, while 
a different line of work
also considers 
convergence of the averaged iterate 
\citep{nemirovski_prox-method_2004},
and another considers convergence to any critical point \citep{abernethy_last-iterate_2019}.
Also to prevent possible confusion, let us mention a related but distinct line of work \citep{du_linear_2019, doan_convergence_2022} that considers using different step-sizes to update the $x$ and the $y$ variables. 
The resulting two-timescale dynamics can be analyzed globally assuming a min-max analog of the Polyak-Łojasiewicz inequality, using quite different considerations than the ones developed in this paper.

\paragraph{Analysis of GDA and PP for SC-SC, bilinear or convex-SC games.}
The fact that PP converges for 
SC-SC
and bilinear min-max games is well-established since at least \cite{rockafellar_monotone_1976}; for a modern reformulation, see e.g.\ \cite{mokhtari_unified_2020}.
The convergence of simultaneous GDA for SC-SC games and its divergence for bilinear games is also a classical fact, 
see e.g.\ \cite{liang_interaction_2019,lu_osr-resolution_2022}.
The cycling behavior of alternating GDA for bilinear games is proved in \cite[Theorem~4]{bailey2020finite}.

As shown more recently in 
\cite[Appendix~G]{nagarajan_gradient_2017} and
\cite[Theorem~6]{zhang_near-optimal_2022},
convexity in $x$, strong concavity in $y$ and non-degeneracy of the interaction component are in fact sufficient to ensure $\tmu_M > 0$, and so local convergence of GDA.
The setting of the former work
corresponds to ${ \nabla_{xx}^2 f(z^*) = 0 }$, $-\nabla_{yy}^2 f(z^*) \succ 0$ and $\nabla_{xy}^2 f(z^*)$ full-row-rank,
extended in the latter work to $\nabla_{xx}^2 f(z^*) \succeq 0$.
Those two works also provide bounds on 
$\tmu_M$
for those specific cases
in terms of the least eigenvalues of $-\nabla_{yy}^2 f(z^*)$ and of $\nabla^2_{xy} f(z^*) \nabla^2_{xy} f(z^*)^\top$.

\paragraph{Spectral analysis-based convergence analyses.}
Our approach 
to analyzing min-max algorithms is to
directly study the properties of 
the update operator 
(e.g., of $T(z) = z - \begin{pmatrix} \nabla_x f \\ -\nabla_y f \end{pmatrix}(z)$ for simultaneous GDA).
Specifically, by a classical result on discrete-time dynamical systems, the local exponential convergence is characterized by the spectral radius of that operator's Jacobian.
This approach is used for example in 
\cite{gidel_negative_2019} to analyze the local convergence behavior of alternating GDA with negative momentum,
and in \cite{azizian_tight_2020} to derive tight convergence bounds 
-- dependent explicitly on $\eta$ and $\spectrum(M)$ --
for simultaneous GDA and multi-step Extra-Gradient.

\paragraph{Average-case analysis of the (local) convergence rate.}
In \cite{pedregosa2020acceleration,domingo-enrich_averagecase_2021}, the authors
analyze the convergence of gradient methods for affine operator problems -- i.e., for finding $z$ such that $M(z-z^*) = 0$ -- when $M$ is a random normal matrix with a known spectral distribution. They derive average-case optimal methods for this setting.
Our analysis in \sectionref{sec:quantitative} also has an average-case flavor, but our random model is different: we assume that the symmetric $S$ and antisymmetric parts $A$ of $M$ are independent and that $S \ll A$. In particular we do not require $M$ to be normal (i.e., $S$ and $A$ do not commute a priori, which is the typical case in min-max optimization).

\paragraph{Hypocoercivity.}
In the context of partial differential equations (PDEs), the phenomenon that a linear PDE $\partial_t u_t = -L u_t$ may exhibit linear convergence to $0$ even when the generator $L$ is not coercive, is called hypocoercivity and is studied in detail in \cite{villani2009hypocoercivity}.
This is precisely the infinite-dimensional analog of the phenomenon studied in the present paper: 
still denoting by $S, A$ the symmetric resp.\ antisymmetric parts of $M$,
coercivity of $L$ corresponds to $S \succ 0$, while hypocoercivity corresponds to $\tmu_M>0$.
Specifically, \cite{villani2009hypocoercivity} shows how to construct Lyapunov functions to establish hypocoercivity and convergence rates of certain PDEs, by exploiting properties of the iterated commutators of $S$ and $A$.
By contrast our focus is on the finite-dimensional case, where it is easier and more natural to directly study the spectrum of $M$.

\subsection{Notation}

For any matrix $T \in \RR^{d \times d}$ or $\CC^{d \times d}$, 
denote by $\spectrum(T) \subset \CC$ its spectrum, i.e., the set of its eigenvalues,
and by $\rho(T) = \max_{\lambda \in \spectrum(T)} \abs{\lambda}$ its spectral radius.
Recall that the spectral radius is distinct from the operator norm, which is the largest singular value, although they coincide for Hermitian matrices.
Denote eigenspaces as $E_\lambda(T) = \left\lbrace z \in \CC^d;~ Tz = \lambda z \right\rbrace$
and let $\mathrm{Eigvecs}(T) = \bigcup_{\lambda \in \spectrum(T)} E_\lambda(T)$ the set of all (complex) eigenvectors.
$\norm{\cdot}$ denotes Euclidean or Hermitian norm.
For a collection of square matrices (resp.\ scalars) $(C_k)_k$, $\Diag((C_k))$ denotes the (block-)diagonal matrix with blocks (resp.\ coefficients) $(C_k)_k$.



\section{Characterization of local convergence} \label{sec:quadr}

Consider the general smooth min-max game of Equation~\eqref{eq:intro:smooth-minmax} and assume $M$ is invertible.
We decompose $M$ into its symmetric and antisymmetric parts as
\begin{equation}
\label{eq:quadr:def_M_S_A}
    M = 
    \begin{bmatrix}
        \nabla_{xx}^2 f & \nabla_{xy}^2 f \\
        -\nabla_{xy}^2 f^\top & -\nabla_{yy}^2 f
    \end{bmatrix}(z^*)
    \eqqcolon \begin{bmatrix}
        Q & P \\ 
        -P^\top & R
    \end{bmatrix},
    \qquad
    S \coloneqq \begin{bmatrix}
        Q & 0 \\ 
        0^\top & R
    \end{bmatrix},
    \qquad
    A \coloneqq \begin{bmatrix}
        0 & P \\ 
        -P^\top & 0
    \end{bmatrix}.
\end{equation}
By the second-order optimality condition in the definition of NE, $Q \succeq 0$ and $R \succeq 0$.

Following \cite{letcher_differentiable_2019}, 
$S$ and $A$ can be thought of intuitively as the ``potential'' resp.\ ``interaction'' (or Hamiltonian) components of the two-player zero-sum game \eqref{eq:intro:smooth-minmax}.
Indeed, consider the quadratic game ~
$
    \min_{x \in \RR^n} \max_{y \in \RR^m}
    \left\{ 
        \frac{1}{2} 
        \begin{pmatrix} x \\ y \end{pmatrix}^\top 
        \nabla^2 f(z^*)
        \begin{pmatrix} x \\ y \end{pmatrix} 
        ~=~
        \frac{1}{2} x^\top Q x - \frac{1}{2} y^\top R y + x^\top P y
    \right\}
$,
which is essentially sufficient to understand
the local behavior of GF for \eqref{eq:intro:smooth-minmax} around $z^*$.
Then we can interpret the terms in $Q$ and $R$ as quadratic potentials to be optimized independently by each player, with the bilinear term in $P$ capturing all of the interaction between the players.

\paragraph{Partial curvature generically suffices.}
Let us recall that $\tmu_M = \min_{\lambda \in \spectrum(M)} \Re(\lambda)$ governs the local convergence of GF around the local NE $z^*$.
It is a general fact that this quantity is nonnegative (the proof is included below). In the following theorem, 
we give necessary and sufficient conditions for it to be positive, in terms of geometric conditions on $Q, R$ and $P$.

\begin{theorem} \label{thm:quadr:eigvals}
    The following conditions are equivalent:
	\begin{enumerate}[label=$(\roman*)$,itemsep=0mm]
		\item $\tmu_M > 0$.
		\item $\mathrm{Eigvecs}(A) \cap \Ker S = \{0\}$.
		\item For 
		any eigenvector $x$ of $P P^\top$
		(i.e., left-singular vector of $P$), 
		$x \not\in \Ker Q$ or $P^\top x \not\in \Ker R$.
		\item For 
		any eigenvector $y$ of $P^\top P$
		(i.e., right-singular vector of $P$),
		$Py \not\in \Ker Q$ or $y \not\in \Ker R$.
	\end{enumerate}
\end{theorem}

As a consequence of \autoref{thm:quadr:eigvals}, the condition $\tmu_M>0$ holds generically in the following sense: for any fixed ${S\neq 0}$, the set of matrices $P$ such that $\tmu_M>0$ is dense and open in $\RR^{n \times m}$. In particular, this property holds with probability $1$ if $P$ is drawn from an absolutely continuous distribution and is independent from $Q$ and $R$, as in the experiment of \figureref{fig:intro:regbilinear}.

\begin{proof}
	Let $\lambda \in \spectrum(M)$ and $z \in \CC^{n+m}$ non-zero such that $M z = (S+A) z = \lambda z$.
    Since $S$ and $A$ are real and symmetric resp.\ antisymmetric,
	\begin{align}
		\overline{\olz^\top S z}
		&= z^\top S \olz
		= (z^\top S \olz)^\top 
		= \olz^\top S z,
		~~\text{so}~~
		\olz^\top S z \in \RR
		 \\
		\text{and}~~~
		\overline{\olz^\top A z}
		&= z^\top A \olz
		= (z^\top A \olz)^\top
		= -\olz^\top A z,
		~~\text{so}~~
		\olz^\top A z \in i\RR.
	\end{align}
	So by taking the real part in $\olz^\top (S+A) z = \lambda \norm{z}^2$,
	\begin{equation}
		\Re(\lambda) \norm{z}^2
		= \olz^\top S z
		= \Re(z)^\top S \Re(z) + \Im(z)^\top S \Im(z)
		\geq 0
	\end{equation}
    since $S = \Diag(Q, R) \succeq 0$.
    This shows that $\tmu_M \geq 0$.

	Now let us show the equivalence of the conditions.
	\begin{itemize}
		\item[$(ii) \implies (i)$:~]
		By contraposition, suppose there exists
		$\lambda \in \spectrum(M)$ with $\Re(\lambda)=0$,
		and let $z$ non-zero such that $Mz = \lambda z$.
		Then 
		$\olz^\top S z = \Re(\lambda) \norm{z}^2 = 0$, 
		so $z \in \Ker S$.
		So $M z = Az = \lambda z$, and $z \in \left( E_\lambda(A) \cap \Ker S \right) \setminus \{0\}$.
		\item[$(i) \implies (ii)$:~]
		By contraposition, suppose there exists $\lambda \in \spectrum(A)$ and a non-zero $z \in E_\lambda(A) \cap \Ker S$. Since $A$ is antisymmetric, then
		$\Re(\lambda) = 0$. On the other hand, $Mz = (S+A)z = Az = \lambda z$, i.e., $\lambda \in \spectrum(M)$.
        So $\tmu_M \leq \Re(\lambda) =0$.
		\item[$(i) \implies (iii), (iv)$:~]
		By contraposition, suppose there exists an eigenvector $x$ of $PP^\top$ such that $x \in \Ker Q$ and $P^\top x \in \Ker R$, and denote $\sigma \in \RR$ such that $PP^\top x = \sigma^2 x$ (since $PP^\top \succeq 0$).
		Then
		\begin{align}
			M \begin{pmatrix}
				i\sigma~ x \\ P^\top x
			\end{pmatrix}
			= \begin{bmatrix}
				Q & P \\
				-P^\top & R
			\end{bmatrix}
			\begin{pmatrix}
				i\sigma~ x \\ P^\top x
			\end{pmatrix}
			= \begin{pmatrix}
				\sigma^2 x \\ -i \sigma P^\top x
			\end{pmatrix}
			= -i \sigma \begin{pmatrix}
				i\sigma~ x \\ P^\top x
			\end{pmatrix}
		\end{align}
		and so $-i\sigma \in \spectrum(M)$.
		This shows $(i) \implies (iii)$, and $(i) \implies (iv)$ follows analogously.
		\item[$(iii), (iv) \implies (ii)$:~]
        By contraposition, suppose there exists $\lambda = i\sigma \in \spectrum(A)$ and a non-zero 
        $z = (x, y) \in E_\lambda(A) \cap \Ker S$.
		Expanding the blocks in $Az = \lambda z$,
		\begin{align}
			\begin{cases}
				Py = i\sigma x \\
				-P^\top x = i\sigma y
			\end{cases}
			~~~\text{and so}~~~
			\begin{cases}
				P^\top P y = \sigma^2 y \\
				P P^\top x = \sigma^2 x.
			\end{cases}
		\end{align}
		Moreover, this implies that $x = 0 \iff y = 0$, and since $z \neq 0$, then both $x \neq 0$ and $y \neq 0$.
		So $x$ is an eigenvector of $PP^\top$ and, since $Sz=0$, $x \in \Ker Q$ and $P^\top x = -i\sigma y \in \Ker R$, which contradicts $(iii)$.
		Likewise, $y$ is an eigenvector of $P^\top P$ and $y \in \Ker R$ and $P y = i\sigma x \in \Ker Q$, which contradicts $(iv)$.
        \qedhere
	\end{itemize}
\end{proof}

\paragraph{Geometric interpretation using real vectors.}
We draw the attention of the reader to the fact that $(ii)$ involves complex eigenvectors. For a rephrasing in terms of real objects, note that
if $z \in E_{i\sigma}(A)$ with $\sigma \in \mathbb{R}$, then 
$A \begin{bmatrix} \Re(z) & \Im(z) \end{bmatrix}
= \begin{bmatrix} \Re(z) & \Im(z) \end{bmatrix} 
\begin{bmatrix}
    0 & \sigma \\
    -\sigma & 0
\end{bmatrix}$;
geometrically, $F_{i\sigma}(A) = \mathrm{span}(\Re(z),\Im(z))$ is a ``rotation plane'' of the GF for the bilinear game
$\min_x \max_y x^\top P y$,
in the sense that the projection of GF on $F_{i\sigma}(A)$ is a circular motion with constant speed $\sigma$.
Condition $(ii)$ expresses that for each such $F_{i\sigma}(A)$, 
there exists an eigenspace $E_\mu(S)$ of $S$ (for a $\mu > 0$) that is not orthogonal to it.
This causes the GF for 
$\min_x \max_y \frac{1}{2} x^\top Q x - \frac{1}{2} y^\top R y + x^\top P y$
projected on $F_{i\sigma}(A)$ to spiral down to
$0$ instead of cycling around it.

One may naturally wonder whether a notion of non-orthogonality between the potential ($S$) and interaction components ($A$) can be used to bound the convergence quantitatively; this is developed in the next section in the particular case $S\ll A$.



\section{Convergence rate when interaction dominates (\texorpdfstring{$S \ll A$}{S<<A})} \label{sec:quantitative}

Let us now discuss the case of games with a small symmetric part, that is, whose Jacobian at optimum is $M_\alpha = A+\alpha S$ for some 
symmetric $S = \begin{bmatrix}
    Q & 0 \\ 
    0^\top & R
\end{bmatrix}$,
antisymmetric $A = \begin{bmatrix}
    0 & P \\ 
    -P^\top & 0
\end{bmatrix}$
and some small $\alpha>0$. 
In this section we assume $n=m$
(the general case is technically more challenging as \autoref{prop:quantitative:ord1} requires $A$ to have distinct eigenvalues, which requires $\abs{n-m} \leq 1$).

\subsection{Convergence rate of Gradient Flow} \label{subsec:quantitative:GF}

As discussed previously, the normalized local exponential convergence rate $r/\eta$ of gradient methods in the asymptotic regime $\eta \to 0$ -- or equivalently, the convergence rate of GF -- for a game with Jacobian at optimum $M_\alpha$ is equal to $\tmu_{M_\alpha}$. 
We can estimate this quantity using the standard formula for the asymptotic expansion of the eigenvalues of a perturbed matrix, which takes an interesting form in our context.

\begin{proposition} \label{prop:quantitative:ord1}
    Suppose that $P$ is full-rank and has distinct singular values, and let $P = U \Sigma V^\top = \sum_{j=1}^n \sigma_j u_j v_j^\top$ be its singular value decomposition. Then it holds
    \begin{equation}
        \tmu_{M_\alpha} = \frac{1}{2} \alpha \Big(\min_{1\leq j\leq n} u_j^\top Q u_j + v_j^\top Rv_j\Big) + O(\alpha^3).
    \end{equation}
\end{proposition}

This expansion explains in particular why the normalized convergence rate 
$r/\eta \sim \tmu_{M_\alpha}$
is approximately proportional to $\alpha$ in \figureref{fig:intro:regbilinear}.
Interestingly, the error term is $O(\alpha^3)$, which suggests that the approximation can be reasonably accurate even for quite large values of $\alpha$, as illustrated in 
\autoref{fig:intro:regbilinear:spectrum}.

\begin{proof}
    Here $M_0=A$ has distinct eigenvalues $\left\lbrace i s \sigma_j, s \in \{-1,1\}, 1 \leq j \leq n \right\rbrace$, with unit-norm eigenvectors
    ${ A \begin{pmatrix} -i s u_j/\sqrt{2} \\ v_j/\sqrt{2} \end{pmatrix} = i s \sigma_j \begin{pmatrix} -i s u_j/\sqrt{2} \\ v_j/\sqrt{2} \end{pmatrix} }$.
    By the calculation of the eigenvalue derivatives from \cite{tao_when_2008}, we obtain the following expansion for $\spectrum(M_\alpha)$:
    \begin{multline} 
        \spectrum(M_\alpha)
        = \Bigg\lbrace
            i s \sigma_j + \frac{1}{2} \alpha \left( u_j^\top Q u_j + v_j^\top R v_j \right) 
            + \frac{1}{4} \alpha^2 \sum_{(s', j') \neq (s, j)} \frac{1}{is \sigma_j - is' \sigma_{j'}} \left( s s' u_{j'}^\top Q u_j + v_{j'}^\top R v_j \right)^2 \\
            + O(\alpha^3),~~
            s \in \{-1,1\}, 1 \leq j \leq n
        \Bigg\rbrace.
    \label{eq:quantitative:ord1}
    \end{multline}
    Now the zeroth- and second-order terms are all in $i\RR$, hence the announced expansion for $\tmu_{M_\alpha}$.
\end{proof}

\paragraph{Estimate of the leading term under a probabilistic model.}

Assuming the singular vectors $(u_1, ..., u_n)$, $(v_1, ..., v_n)$ of $P$ are distributed uniformly at random -- which is the case for example if $P$ has i.i.d.\ Gaussian entries by rotational invariance%
\footnote{
    If $P = U \Sigma V^\top$ has i.i.d.\ Gaussian entries, then $P$ has the same law as $\tU^\top P \tV$, for any $\tU, \tV \in \OOO_n$ the set of $n \times n$ orthonormal matrices.  So $(U, V)$ has the same law as $(\tU^\top U, \tV^\top V)$. This shows that $(U, V)$ is distributed according to the Haar measure on the product group $\OOO_n \times \OOO_n$, that is, $U$ and $V$ are independently distributed uniformly on $\OOO_n$.
} --,
the leading term in the expansion of $\tmu_{M_\alpha}$ can be estimated in expectation as follows.
\ifextended%
    The proof is placed in \appendixref{apx:RMT_proofs}, where we also include a high-probability version of the estimate as \autoref{prop:RMTinprob:spreadout}.
\else%
    The proof is placed in the Appendix, where we also include a high-probability version of the estimate.
\fi

\begin{restatable}{proposition}{propRMTspreadout}
\label{prop:quantitative:RMT_spreadout}
    Suppose $Q, R$ are fixed and 
    $U, V$ are independently distributed uniformly on the set of $n \times n$ orthonormal matrices.
    Then
    \begin{equation}
        \frac{\trace(S)}{n} \left( 1 - 2\frac{\norm{S}_F}{\trace(S)} \sqrt{\log n} \right)
        \leq
        \EE \left[ \min_{1\leq j\leq n} u_j^\top Q u_j + v_j^\top Rv_j \right]
        \leq 
        \frac{\trace(S)}{n}
    \end{equation}
    where $\norm{\cdot}_F$ denotes the Frobenius norm.
    In particular, provided that
    $\frac{\trace(S)}{\norm{S}_F} \geq 2 \sqrt{\log n} (1+c)$ for some fixed $c>0$,
    we have $\EE \left[ \min_{1\leq j\leq n} u_j^\top Q u_j + v_j^\top Rv_j \right] \asymp \frac{\trace(S)}{n}$ as $n \to \infty$.
\end{restatable}

Note that $\frac{\trace(S)}{\norm{S}_F}$, which always lies in the interval $[1,\sqrt{2n}]$, is a measure of the effective sparsity of the spectrum of $S$
(larger meaning less sparse), so
the condition $\frac{\trace(S)}{\norm{S}_F} \geq 2 \sqrt{\log n} (1+c)$ means the spectrum of $S$ is well spread-out.
So the proposition shows that
\textbf{the exponential convergence rate depends on the \emph{average} of the eigenvalues of $S$,} when $\alpha S \ll A$ and the spectrum of $S$ is well spread-out.
This fact should be contrasted with the case of minimization, where the convergence rate scales as the \emph{minimum} eigenvalue of the Hessian.

Interestingly, when the spectrum of $S$ is 
sparse,
the typical behavior of the leading term 
in the expansion of $\tmu_{M_\alpha}$,
$\min_j u_j^\top Q u_j + v_j^\top R v_j$, is quite different.
In this case, that quantity depends on the geometric mean of the non-zero eigenvalues of $S$, rather than the arithmetic mean,
as formalized in the following proposition.
\ifextended%
    The proof is placed in \appendixref{apx:RMT_proofs}, along with a high-probability version of the estimate in \autoref{prop:RMTinprob:sparse}.
\else%
    The proof is placed in the Appendix, along with a high-probability version of the estimate.
\fi

\begin{restatable}{proposition}{propRMTsparse}
\label{prop:quantitative:RMT_sparse}
    Suppose $Q, R$ are fixed and
    $U, V$ are independently distributed uniformly on the set of $n \times n$ orthonormal matrices.
    Let $s_1 \geq ... \geq s_r > 0 = s_{r+1} = ... = s_{2n}$ the eigenvalues of 
    $S$.
    Then
    \begin{equation}
        \EE \left[ \min_{j \leq n} u_j^\top Q u_j + v_j^\top R v_j \right]
        \geq
        \max_{\SSS \subset \{1,...,r\}}
        \frac{1}{e} \frac{\card{\SSS}}{n} n^{-\frac{2}{\card{\SSS}}} \left[ \prod_{l \in \SSS} s_l \right]^{\frac{1}{\card{\SSS}}}.
    \end{equation}
    In particular, 
    $\EE \left[ \min_{j \leq n} u_j^\top Q u_j + v_j^\top R v_j \right]
    \gtrsim n^{-\frac{2}{r}-1}$
    when $n \to \infty$ and $r$ and $s \in \RR^r$ are fixed.
\end{restatable}

Numerically, the quantity $\min_{j \leq n} u_j^\top Q u_j + v_j^\top R v_j$ indeed scales as $n^{-\frac{2}{r}-1}$ under the conditions of the proposition, suggesting that our lower estimate could be tight.

\subsection{Convergence rate of discrete-time algorithms} \label{subsec:quantitative:DT_algos}

\autoref{prop:quantitative:ord1} gave an expansion of the normalized convergence rate $r/\eta$ of gradient methods (for a game with Jacobian at optimum $M_\alpha$), non-asymptotically in $\alpha$ and in the asymptotic limit $\eta \to 0$.
In this subsection, we give expansions of the convergence rate $r$ that are non-asymptotic in $\alpha$ and $\eta$.

The algorithms we will consider can be written in the form $z^{k+1} = T(z^k)$ with the update operator $T$ dependent only on $\nabla f$ and on step-size $\eta$, and satisfying $z^* = T(z^*)$.
It is well-known that local convergence of such methods is determined by $\rho(\nabla T(z^*))$, where $\nabla T$ is the Jacobian of $T$ and $\rho(\cdot)$ denotes spectral radius. 
Namely, if $\rho(\nabla T(z^*)) < 1$ then the iterates converge locally with $\norm{z^k-z^*} = O \left( (\rho(\nabla T(z^*)) + \eps)^k \right)$, with $\eps>0$ an arbitrarily small constant \citep[Proposition~4.4.1]{bertsekas_nonlinear_1997}.

In \autoref{tab:expans} we give an expansion for $\rho(\nabla T(z^*))$ for three classical gradient methods: simultaneous GDA (Sim-GDA), alternating GDA (Alt-GDA),
and Extra-Gradient (EG).
Let us clarify immediately that throughout this paper, statements made about ``GDA'' without further specification apply to both Sim-GDA and Alt-GDA.
In the second column we denoted by
$g(z) = \begin{pmatrix} \nabla_x f \\ -\nabla_y f \end{pmatrix}(z)$ 
the skewed gradient field of the game,
and in the third column we wrote $M$ instead of $M_\alpha$ for concision.
In the fourth column, $\pm i \sigma_j$ denotes the eigenvalues of $A$ assumed distinct and $w_j$ are the associated eigenvectors
-- equivalently, the singular value decomposition of $P$ is $P = \sum_{j=1}^n \sigma_j u_j v_j^\top$ and for each $j \leq n$, $\olw_j^\top S w_j = \olw_{j+n}^\top S w_{j+n} = \frac{1}{2} \left( u_j^\top Q u_j + v_j^\top R v_j \right)$.
We refer to 
\ifextended%
    \autoref{apx:DT_expansions} 
\else%
    the Appendix
\fi
for the derivation of this table and for explicit bounds on the ``$O(\cdot)$'' terms.

\begin{table}[t]
    \caption{Expansions of $\rho(\nabla T(z^*))$ in $\alpha$ and $\eta$ for classical gradient methods}
    \label{tab:expans}
    \centering
    \begin{tabular}{c l l l}
        \toprule
        Algorithm & $T(z)$ & $\nabla T(z^*)$ & \makecell[l]{$\rho(\nabla T(z^*))^2 = \max_{j \leq 2n} [...]$ \\ $\qquad\qquad\qquad~ + O(\eta \alpha^3 + \eta^2 \alpha^2)$} \\
        \midrule
        Sim-GDA & $z - \eta g(z)$ & $I - \eta M$ & $1 - 2 \alpha \eta \left( \olw_j^\top S w_j \right) + \eta^2 \sigma_j^2$ \\ 
        Alt-GDA & {\small see text} & $I - \eta (I - \frac{\eta}{2} A) M + O(\eta^3)$ & $1 - 2 \alpha \eta \left( \olw_j^\top S w_j \right) + O(\eta^3)$ \\ 
        EG & $z - \eta g(z - \eta g(z))$ & $I - \eta (I - \eta M) M$ & $1 - 2 \alpha \eta \left( \olw_j^\top S w_j \right) - \eta^2 \sigma_j^2 + O(\eta^3)$ \\
        \bottomrule
    \end{tabular}
\end{table}

Informally, as one can directly see from the fourth column, Sim-GDA requires a very small step-size for the first term in $\eta$ to overcome the terms $+\eta^2 \sigma_j^2$, while for EG those terms actually appear with a favorable sign, and Alt-GDA neither benefits nor suffers from those terms.
We also see that Alt-GDA is quite faithful to GF, in that their normalized convergence rates coincide up to terms of order $\eta^3 + \alpha^4$.
All of these insights are in line with common intuition in the min-max optimization literature \cite{bailey2020finite, lu_osr-resolution_2022},
as well as with our numerical experiments for the next section, \autoref{fig:particleMNE:manyalphas}.

\paragraph{A symmetrized formulation of Alt-GDA.}
In order to derive the rate for Alt-GDA, we used the following symmetrized formulation of it: we
let $(x^0, y^\half) \in \RR^d \times \RR^d$ and
\begin{equation}
    \begin{cases}
        \forall k \in \NN,
        & x^{k+1} = x^k - \eta \nabla_x f(x^k, y^{k+\half}) \\
        \forall k \in \NN + \half,
        & y^{k+1} = y^k + \eta \nabla_y f(x^{k+\half}, y^k)
    \end{cases}
    ~~~~\text{and}~~~~
    \begin{cases}
        \forall k \in \NN,
        & x^{k+\half} = \frac{x^{k+1} + x^k}{2} \\
        \forall k \in \NN+\half,
        & y^{k+\half} = \frac{y^{k+1} + y^k}{2}.
    \end{cases}
\end{equation}
That is, $x$ gets updated with the gradient rule at integer time-steps, $y$ at half-integer time-steps, and we define $x^k$, $y^k$ at non-updating time-steps as the average of the preceding and following updating time-steps.
Assuming $\eta \leq \norm{\nabla_{xx}^2 f}_\infty^{-1} \wedge \norm{\nabla_{yy}^2 f}_\infty^{-1}$, we show in
\ifextended%
    \autoref{apx:sym_altGDA} 
\else%
    the Appendix
\fi
that $z^{k+1}$ is indeed entirely determined by $z^k = (x^k, y^k)$ for each $k \in \NN$, and that the associated update operator $T$ satisfies
\begin{equation}
    T(z) = z - \eta g(z) + \frac{\eta^2}{2} A(z) g(z) + O(\eta^3 \norm{g(z)})
    ~~~~\text{where}~~~~
    A(z) = \frac{\nabla g(z) - \nabla g(z)^\top}{2}.
\end{equation}

For comparison, the usual formulation of Alt-GDA considers as the iterates $(\tx^k, \ty^k) = (x^k, y^{k+\half})$.
We emphasize that $(\tx^k, \ty^k)$ and $(x^k, y^k)$ have the same convergence rate if they converge exponentially, as one can check directly from the definition.



\section{Illustration: sparse mixed Nash equilibria of continuous games} \label{sec:particleMNE}

In this section we apply the above considerations to a particular class of min-max problems, which is of its own interest in game theory.
Namely we consider the classical problem of finding the mixed Nash equilibria (MNE) of two-player zero-sum games, that is, given strategy spaces $\XXX$, $\YYY$ and a payoff function $f: \XXX \times \YYY \to \RR$, solving the min-max problem
\begin{equation}
    \min_{\mu \in \PPP(\XXX)} \max_{\nu \in \PPP(\YYY)} 
    \left\{ F(\mu,\nu) = \EE_{x \sim \mu, y \sim \nu} [f(x, y)] \right\}.
\end{equation}
Here $\PPP(\XXX)$ denotes the space of probability measures -- representing mixed strategies -- over $\XXX$.
Let us focus on cases where $\XXX$ and $\YYY$ are continuous sets, say, $\XXX=\YYY=\TT^1$ the one-dimensional Euclidean torus,%
\footnote{
    The choice of $\XXX=\YYY=\TT^1$ is made for simplicity of exposition. The discussion below extends straightforwardly to toruses of any dimension,
    and could be extended to $\XXX, \YYY$ Riemannian manifolds without boundaries at the cost of more technical notation.
}
and $f: \TT^1 \times \TT^1 \to \RR$ is smooth.
The above min-max problem is then infinite-dimensional, and algorithms to solve it explicitly must be based on reparameterized formulations.
More specifically, suppose that the MNE $(\mu^*, \nu^*)$ is unique and ``sparse'', i.e., has finite supports:
$\support(\mu^*) = \left\lbrace x^*_I, 1 \leq I \leq N \right\rbrace$,
$\support(\nu^*) = \left\lbrace y^*_J, 1 \leq J \leq M \right\rbrace$
(this is the case for example if $f$ is a sum of separable functions \citep{stein_separable_2008}).

In this setting there are two natural reparameterizations and associated algorithms:
\begin{enumerate}[leftmargin=5mm,itemsep=0mm]
    \item If the optimal support points $\{x^*_I\}_I$, $\{y^*_J\}_J$ are known, then
    we may reparameterize by 
    $\mu = \sum_{I=1}^N a_I \delta_{x^*_I}$, 
    $\nu = \sum_{J=1}^M b_J \delta_{y^*_J}$
    and optimize over the $a_I, b_J$.
    The problem reduces to a constrained bilinear game
    \begin{equation} \label{eq:particleMNE:MP_pb}
        \min_{a \in \Delta_N} \max_{b \in \Delta_M} 
        \left\lbrace 
            F_1(a,b)
            = \sum_{I=1}^N \sum_{J=1}^M a_I b_J~ f(x^*_I, y^*_J)
            \eqqcolon a^\top P b
        \right\rbrace
    \end{equation}
    where $\Delta_N$
    denotes the standard simplex.
    A classical approach is then to apply Mirror Prox (MP) with entropy link function 
    \citep{nemirovski_prox-method_2004}
    (MP is the Bregman-geometry analog of EG).

    \item If only the number of optimal support points is known, then we may reparameterize by 
    $\mu = \sum_{I=1}^N a_I \delta_{x_I}$,
    $\nu = \sum_{J=1}^M b_J \delta_{y_J}$
    and optimize over both the weights ($a_I, b_J$) and the support points ($x_I, y_J$).
    The problem reduces to
    \begin{equation} \label{eq:particleMNE:CPMP_pb}
        \min_{(a, x) \in \Delta_N \times (\TT^1)^N} ~
        \max_{(b, y) \in \Delta_M \times (\TT^1)^M} ~
        \left\{
            F_2(a, x, b, y)
            = \sum_{I=1}^N \sum_{J=1}^M a_I b_J~ f(x_I, y_J)
        \right\}.
    \end{equation}
    A possible approach is to iteratively update (simultaneously) the $a, b$ using MP steps with step-size $\eta$ and the $x, y$ using EG steps with step-size $\gamma \eta$, for some parameter $\gamma>0$.
    This algorithm is called Conic Particle Mirror Prox in \cite{wang_exponentially_2022}, but for concision we will refer to it simply as ``EG'' in this section.
    Note that the main challenge in that reference is to deal with the case where $N$ and $M$ are unknown,
    but here we assume they are known for the sake of simplicity.
\end{enumerate}

\begin{figure}[t]
    \vspace{-0.5em}
    \begin{subfigure}[t]{0.48\textwidth}
        \centering
    	\includegraphics[width=\textwidth]{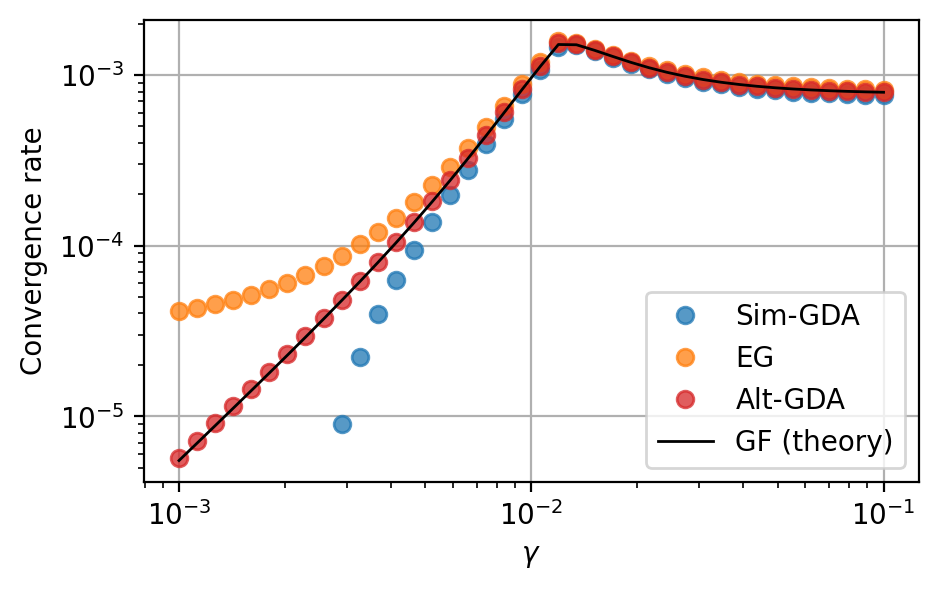}
        \\[-0.5em]
        \caption{
            Convergence rate $r$ vs.~$\gamma$ (fixed $\eta=10^{-2}$)
        }
        \label{fig:particleMNE:rate_vs_alpha}
    \end{subfigure}
    \hfill
    \begin{subfigure}[t]{0.48\textwidth}
        \centering
    	\includegraphics[width=\textwidth]{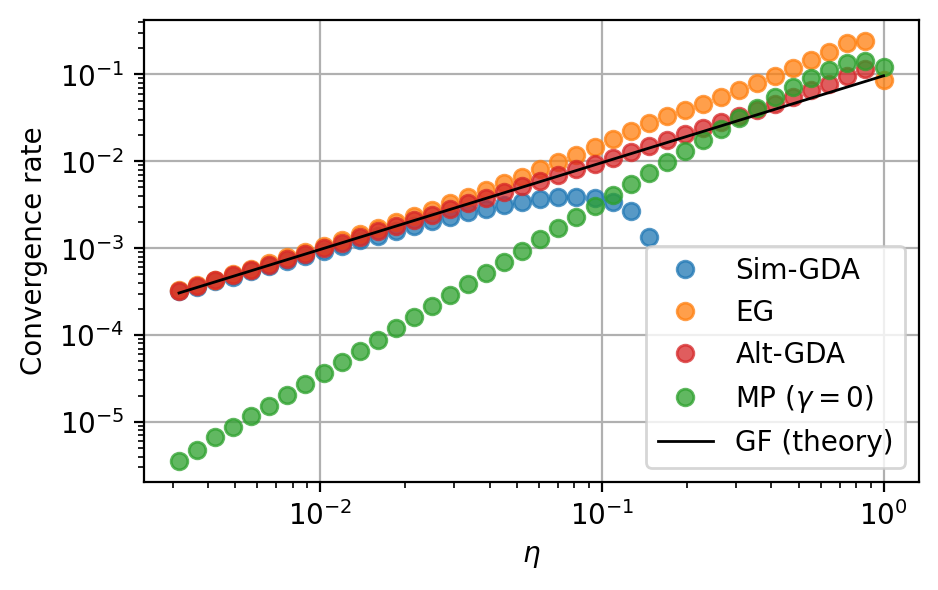}
        \\[-0.5em]
        \caption{
            Convergence rate vs.~$\eta$ (fixed $\gamma=10^{-2}$)
        }
        \label{fig:particleMNE:rate_vs_eta}
    \end{subfigure}
    \caption{Observed local convergence rates $r$ (i.e., $\norm{z^k - z^*} = \Theta((1-r)^k)$) for various conic particle methods using step-size $\eta$ for the weights $(a, b)$ and $\gamma \eta$ for the positions $(x, y)$ and for a fixed random draw of $f$. 
    (left) Fixing $\eta=10^{-2}$, we observe a rate for GF (we plot $\eta\cdot\tilde \mu_\gamma$) scaling as $\gamma^2$ as predicted in~\autoref{prop:particleMNE:scale_gamma2}. Interestingly, Alt-GDA has exactly the same behavior, while the behavior of Sim-GDA and EG differ for small $\gamma$, due to the corrective terms shown in~\autoref{tab:expans}. 
    (right) Fixing $\gamma=10^{-2}$, we observe a rate scaling as $\eta$ for EG and its variants, and as $\eta^2$ for MP; the convergence of EG is mostly faster than MP.}
    \label{fig:particleMNE:manyalphas}
\end{figure}

In \figureref{fig:particleMNE:manyalphas}, we show the dependency on $\eta$ and $\gamma$ of the local convergence rate of MP and EG, as well as that of the analogs of Sim-GDA and Alt-GDA and GF for the problem $\min_{(a,x)} \max_{(b,y)} F_2(a,x,b,y)$ in order to illustrate the insights from \autoref{subsec:quantitative:DT_algos}.
We used a randomly generated payoff function~$f$,
and the convergence is measured by the iterates' $\ell_2$~distance to the solution; see 
\ifextended%
    \appendixref{subsec:details_particleMNE:for_figure} 
\else%
    the Appendix
\fi
for details.

MP (the algorithm based on the reparameterization $F_1$) converges to the MNE for any small enough $\eta$ with an exponential rate at least proportional to $\eta^2$~\cite{wei_linear_2021}, and this scaling is tight numerically (green dots in~\autoref{fig:particleMNE:rate_vs_eta}).
Its explicit variant (Mirror Descent-Ascent) and its continuous-time flow are known to diverge for any $\eta$
\cite{bailey_multiplicative_2018} \cite{mertikopoulos_cycles_2018}.

EG is known to converge locally to the MNE for any small enough $\eta$ (and any $\gamma>0$) with an exponential rate at least proportional to $\eta^2$, despite the non-convexity of $F_2$, under some non-degeneracy assumptions~\cite[Sec.~3.1]{wang_exponentially_2022}.
However numerically the convergence rate of EG typically scales as $\eta$, not $\eta^2$ (orange dots in~\autoref{fig:particleMNE:rate_vs_eta}) -- a fact which we will explain below.
The explicit variant of EG and its continuous-time flow have not previously been analyzed; the discussion below will give a characterization of when they converge locally.

Note that, at least for the particular $f$ and $\gamma$ used for \autoref{fig:particleMNE:rate_vs_eta}, 
EG converges locally faster than MP for the same $\eta$
(in number of iterations, the per-iteration costs differing only by a constant factor), 
even though the former does not use the knowledge of the $\{x^*_I\}_I$, $\{y^*_J\}$!
In other words, even when the optimal support points are known, it is beneficial to use the overparameterized formulation $F_2$ where we also vary the support points.

\paragraph{Overparameterization induces partial curvature.}
Let us inspect the ``Jacobians'' at optimum for the two dynamics, MP vs.\ EG with parameter $\gamma$.
Due to the simplex constraints and the non-Euclidean nature of the updates for $a, b$,
the relevant matrices are $M_{\MP}$ and $M_\gamma$ defined below, in the sense that the exponential convergence rate of each algorithm is $\tmu_M \eta + O(\eta^2)$ 
\ifextended%
    (see \appendixref{subsec:details_particleMNE:MP},~\ref{subsec:details_particleMNE:CPMP}).
\else%
    (see the Appendix).
\fi
Namely,
omitting half of the antisymmetric off-diagonals for readability,
\vspace{-0.25em}  
\begin{equation}
\label{eq:particleMNE:tM_MP}
    M_{\MP} 
    = \begin{bmatrix}
        \bmzero & \Pi_a D_a P D_b \Pi_b^\top \\
        -(*)^\top & \bmzero
    \end{bmatrix}
\vspace{-0.25em}  
\end{equation}
where $D_a = \Diag(\sqrt{a^*})$ (square roots being taken component-wise) and
$\Pi_a \in \RR^{(N-1) \times N}$ is any matrix such that $\Pi_a \Pi_a^\top = I_{N-1}$ and $\Pi_a^\top \Pi_a = I_N - \sqrt{a^*} \sqrt{a^*}^\top$,
and likewise for $D_b$ and $\Pi_b \in \RR^{(M-1) \times M}$;
and for EG, 
\vspace{-0.5em}  
\begin{equation}
\label{eq:particleMNE:tM_CPMP}
    M_{\gamma} = \begin{bmatrix}
        \bmzero & \bmzero & \Pi_a D_a P D_b \Pi_b^\top & \sqrt{\gamma} \ \Pi_a D_a [\partial_y P] D_b \\
        \bmzero & \gamma \Diag(\partial_{xx}^2 P b^*) & \sqrt{\gamma} D_a [\partial_x P] D_b \Pi_b^\top & \gamma D_a [\partial_{xy}^2 P] D_b \\
        -(*)^\top & -\left( * \right)^\top & \bmzero & \bmzero \\
        -\left( * \right)^\top & -\left( * \right)^\top & \bmzero & -\gamma \Diag(\partial_{yy}^2 P^\top a^*)
    \end{bmatrix}
\vspace{-0.1em}  
\end{equation}
where 
$\left[ \partial_x P \right]_{IJ} = \partial_x f(x^*_I, y^*_J)$, and likewise for $\partial_y P$, $\partial_{xx}^2 P$, $\partial_{yy}^2 P$, $\partial_{xy}^2 P$.

For MP, it is clear that the equivalent conditions of \theoremref{thm:quadr:eigvals} are violated for any payoff function $f$, and so $\tmu_{M_{\MP}} = 0$.
For EG, depending on $f$ and $\gamma$, they may or may not be violated.
For all of the random payoff functions we considered in our experiments, we observed that $\tmu_{M_{\gamma}} > 0$, suggesting that the conditions hold generically.
They are violated for certain $f$'s however, as shown in
\ifextended%
    \appendixref{subsec:details_particleMNE:counter_example}, 
\else%
    the Appendix,
\fi
so that the scaling in $\eta^2$ of the convergence rate proved in~\cite{wang_exponentially_2022} is tight in the worst case.
 
More precisely as we show in the next proposition, assuming that the blocks of $M_{\gamma}$ are in general position, we expect $\tmu_{M_{\gamma}}$ to scale as $\gamma^2$, which is indeed observed in the numerical experiment reported in \autoref{fig:particleMNE:rate_vs_alpha}.
The proof, placed in 
\ifextended%
    \autoref{subsec:details_particleMNE:proof_scale_gamma2}, 
\else%
    the Appendix,
\fi
relies on the same tools as \autoref{prop:quantitative:ord1}, that is, on the asymptotic expansions of the eigenvalues of perturbed matrices.

\begin{restatable}{proposition}{particleMNEscaleGammaSq}
\label{prop:particleMNE:scale_gamma2}
    Let $S_2$ symmetric and $A_0, A_1, A_2$ antisymmetric real matrices of the form
    {\small{
    \begin{equation*}
        S_2 = \left[
            \begin{array}{c c|c c}
                \bmzero & & & \\
                & * & & \\
                \hline
                & & \bmzero & \\
                & & & *
            \end{array}
        \right],~
        A_0 = \left[
            \begin{array}{c c|c c}
                & & * & \bmzero \\
                & & \bmzero & \bmzero \\
                \hline
                * & \bmzero & & \\
                \bmzero & \bmzero & &
            \end{array}
        \right],~
        A_1 = \left[
            \begin{array}{c c|c c}
                & & \bmzero & * \\
                & & * & \bmzero \\
                \hline
                \bmzero & * & & \\
                * & \bmzero & &
            \end{array}
        \right],~
        A_2 = \left[
            \begin{array}{c c|c c}
                & & \bmzero & \bmzero \\
                & & \bmzero & * \\
                \hline
                \bmzero & \bmzero & & \\
                \bmzero & * & &
            \end{array}
        \right]
    \end{equation*}
    }}%
    and $M_\gamma = \gamma S_2 + A_0 + \sqrt{\gamma} A_1 + \gamma A_2$ for all $\gamma>0$.
    Then $\tmu_{M_\gamma} = O(\gamma^2)$ as $\gamma \to 0$.
%
\end{restatable}



\section{Conclusion} \label{sec:ccl}

We have investigated the local convergence of gradient methods for min-max games and found that they converge generically under partial curvature. In  more specific settings, we have obtained quantitative estimates of the local convergence rate which exhibit the \emph{average} of the eigenvalues of $S$ as the driving quantity for typical problems.
For the computation of mixed Nash equilibria of continuous games, this leads to a behavior of conic particle gradient methods that is more favorable than that described by the worst-case bounds.

More generally, our analysis leads to the following insights: (i) worst-case bounds might be looser in min-max optimization (compared to minimization) as they fail to capture the interplay between interaction and potential parts; (ii) the addition of new degrees of freedom with curvature typically accelerates the local convergence, as we illustrated in \autoref{sec:particleMNE}.

We note that the phenomenon described in this paper is fundamentally a consequence of
the fact that the skewed gradient field's Jacobian at optimum has a positive-semidefinite symmetric part.
This property is satisfied at local Nash equilibria of min-max games, i.e., of two-player zero-sum differentiable games, but is not true for general differentiable games.



%


\nocite{*}
\printbibliography
\addcontentsline{toc}{section}{\refname} 

\ifextended%
    \newpage
    \appendix
    

\section{Local convergence rates of classical gradient methods} \label{apx:algo_rates}

Consider a game $\min_{x \in \RR^n} \max_{y \in \RR^m} f(x, y)$ with local NE $z^*$ and let
$
    M = \begin{bmatrix}
        \nabla_{xx}^2 f & \nabla_{xy}^2 f \\
        -\nabla_{xy}^2 f^\top~ & -\nabla_{yy}^2 f
    \end{bmatrix}(z^*)
$,
assumed invertible for simplicity.
Gradient-based iterative methods can be written in the form $z^{k+1} = T(z^k)$ with the update operator $T$ dependent only on $\nabla f$ and on step-size $\eta$, and satisfying $z^* = T(z^*)$.
It is well-known that local convergence is determined by $\rho(\nabla T(z^*))$ where $\nabla T$ is the Jacobian of $T$. Namely, if $\rho(\nabla T(z^*)) < 1$ then the iterates converge locally with $\norm{z^k-z^*} = O \left( \rho(\nabla T(z^*)) + \eps)^k \right)$ with $\eps>0$ an arbitrarily small constant \citep[Proposition~4.4.1]{bertsekas_nonlinear_1997}.

\tableref{tab:rates} summarizes the local convergence rates of three classical gradient methods, up to third order in $\eta$. 
The last two columns of the table are valid for $\eta \leq \eta_{\max}$ (given by the second column), and the $O(\cdot)$'s hide only universal constants.
In the last column, we denoted
$\tmu = \tmu_M$ for concision,
and $L = \rho(M)$.
``$k$-EG'' stands for $k$-step Extra-Gradient.
The table was extracted from the proofs of \cite[Appendix~E]{azizian_tight_2020}.


As shown in the fourth column, the quantity $\tmu_M$ appears naturally in the bounds on the local convergence rate
$1-\rho(\nabla T(z^*))$,
with all methods benefitting from $\tmu_M$ being larger.
Moreover, the rate depends linearly on $\tmu_M$ in the asymptotic regime of small $\eta$, as can be seen in the third column, since $1 - \rho(\nabla T(z^*))^2 \sim 2\eta \min_{\lambda \in \spectrum(M)} \Re(\lambda) = 2\eta \tmu_M$.%
\footnote{
    $1 - \rho(\nabla T(z^*))^2$ is the exponential convergence rate for $\norm{z^k-z^*}^2$, and $1 - \rho(\nabla T(z^*))$ is the one for $\norm{z^k-z^*}$. Hence those two quantities differ by a factor $2$ at first order in $\eta$. Equivalently, the additional factor $2$ comes from the fact that
    $\rho^2 = \left[ 1-(1-\rho) \right]^2 \approx 1 - 2(1-\rho)$ when $\rho$ is close to $1$.
}
This asymptotic equivalence corresponds to the fact that GF, the continuous-time flow of all of the classical gradient methods, has local exponential convergence rate equal to $\tmu_M$.
(Of course using a large $\eta$, or following a different continuous-time flow than GF, may lead to faster convergence, but it may require knowledge of problem parameters or the use of more complex adaptive schemes.)

Interestingly, 
as can be seen in 
the third column, EG may have a slower convergence rate than GDA if
$\abs{\Im(\lambda)} < \Re(\lambda)$
for $\lambda = \argmin_{\spectrum(M)} \Re(\cdot)$, 
when $\eta$ is small.

\begin{table}[ht]
  \caption{Moduli of eigenvalues of $\nabla T(z^*)$ for classical gradient methods}
  \label{tab:rates}
  \centering
  \begin{tabular}{c c l l}
    \toprule
    Algorithm & $\eta_{\max}$ & $\left\lbrace \abs{\nu}^2, \nu \in \spectrum(\nabla T(z^*)) \right\rbrace$ & Upper bound on $\rho(\nabla T(z^*))^2$ \\
    \midrule
    Sim-GDA & $\infty$ & $\left\lbrace 1 - 2\eta \Re(\lambda) + \eta^2 \abs{\lambda}^2, \lambda \in \spectrum(M) \right\rbrace$ & $1 - 2\eta \tmu + \eta^2 L^2$ \\ 
    \makecell{$k$-EG \\ ($k \geq 2$)} & \makecell{$(1-c)/L$ \\ ($\forall c>0$)} & \makecell[l]{$\Big\lbrace 1 - 2\eta \Re(\lambda) - \eta^2 \left( \abs{\lambda}^2 - 4 \Re(\lambda)^2 \right)$ \\ \qquad~~~~ $+ O(\frac{1}{1-c} \eta^3 \abs{\lambda}^3), \lambda \in \spectrum(M) \Big\rbrace$} & $1 - \frac{1}{(2-c)^2} \max\left( 2\eta \tmu, \eta^2 L^2 \right)$ \\
    PP & $\infty$ & $\left\lbrace \frac{1}{1+2\eta \Re(\lambda) + \eta^2 \abs{\lambda^2}}, \lambda \in \spectrum(M) \right\rbrace$ & $1 - \max\left( \frac{2\eta \tmu}{1+2\eta \tmu}, \frac{\eta^2 L^2}{1+\eta^2 L^2} \right)$ \\ 
    \bottomrule
  \end{tabular}
\end{table}



\section{Details for \autoref{subsec:quantitative:GF}}
\label{apx:RMT_proofs}

\subsection{Proof of \texorpdfstring{\propositionref{prop:quantitative:RMT_spreadout}}{Proposition~\ref{prop:quantitative:RMT_spreadout}}}
\label{subsec:RMT_proofs:RMT_spreadout}

For ease of reference, we restate the proposition below.
\propRMTspreadout*

For the sake of concision, let
$W = \begin{bmatrix} U \\ V \end{bmatrix}$
and $w_j = W_{\bullet j} = \begin{pmatrix} u_j \\ v_j \end{pmatrix}$,
and (in this appendix only) 
$\mu \coloneqq \min_{j \leq n} u_j^\top Q u_j + v_j^\top R v_j 
= \min_{j \leq n} w_j^\top S w_j$
the quantity which we want to estimate.
Furthermore, denote $s_1 \geq ... \geq s_r > 0 = s_{r+1} = ... = s_{2n}$ the eigenvalues of $S$,
and assume $S$ is diagonal w.l.o.g.
Also let, for all $j \leq n$, $N_1, ..., N_n, M_1, ..., M_n \sim \chi^2_n$ i.i.d.\ and independent of $U$ and $V$, and pose 
$a_j = \begin{pmatrix} \sqrt{N_j}~ u_j \\ \sqrt{M_j}~ v_j \end{pmatrix}$.

Note that for each $j$, $a_j \sim \NNN(0, I_{2n})$.
Indeed, $\sqrt{N_j}~ u_j \sim \NNN(0, I_n)$ since it is isotropic and its norm has the correct distribution, and likewise for $\sqrt{M_j}~ v_j$;
so $a_j \sim \NNN(0, I_{2n})$ as the concatenation of two independent standard Gaussian vectors.

\begin{lemma} \label{lm:RMT:reduce_to_gaussians}
    We have
    $\frac{1}{n} \EE\left[ \min_{j \leq n} a_j^\top S a_j \right] \leq \EE \mu \leq \frac{\trace(S)}{n}$.
\end{lemma}

\begin{proof}
    Let $J \coloneqq \argmin_{j \leq n} w_j^\top S w_j$.
    Then $a_J^\top S a_J = N_J u_J^\top Q u_J + M_J v_J^\top R v_J$.
    Now $J$ is a deterministic function of the $(w_j)_{j \leq n}$, and $(N_j, M_j)_{j \leq n}$ are i.i.d.\ and independent of $(w_j)_{j \leq n}$, so we have the conditional expectation $\EE[N_J | (w_j)_{j \leq n}] = \EE N_1 = \EE \chi^2_n = n$ and likewise for $M_J$.
    Thus, 
    \begin{equation}
        \EE \left[ a_J^\top S a_J \middle| (w_j)_{j \leq n} \right]
        = \EE \left[ N_J \middle| (w_j)_{j \leq n} \right]
        u_J^\top Q u_J 
        + \EE \left[ M_J \middle| (w_j)_{j \leq n} \right]
        v_J^\top R v_J
        = n \mu,
    \end{equation}
    and the first inequality follows by taking total expectations.

    For the second inequality,
    $\EE \mu 
    \leq \EE u_1^\top Q u_1 + v_1^\top R v_1$.
    Now letting $a^{(1)}_1 = \sqrt{N_1}~ u_1 \sim \NNN(0,I_n)$,
    $N_1 u_1^\top Q u_1 = (a^{(1)}_1)^\top Q a^{(1)}$,
    so taking expectations on both sides,
    $n \EE u_1^\top Q u_1 = \trace(Q)$.
    So
    $\EE \mu
    \leq \frac{\trace(Q)+\trace(R)}{n} = \frac{\trace(S)}{n}$.
\end{proof}

Let for concision $\zeta_j = a_j^\top S a_j$ for each $j \leq n$.
\begin{lemma} \label{lm:RMT:MGF_zetaj}
    The moment-generating function of $\zeta_1 (\sim \zeta_2 \sim ... \sim \zeta_n)$ is
    \begin{equation}
        \EE e^{t \zeta_1}
        = \prod_{l=1}^r (1-2 s_l t)^{-1/2}
        ~~~\text{for all}~~~
        t<\frac{1}{2 \max_l s_l}.
    \end{equation}
\end{lemma}

\begin{proof}
    Since $\zeta_1 = \sum_{l=1}^{2n} s_l a_1[l]^2$ and $a_1[l]^2 \sim \chi^2$ i.i.d., we have
    $\EE e^{t a_1[l]^2} = (1-2t)^{-1/2}$ for all $t<\frac{1}{2}$ and
    $
        \EE e^{t \zeta_1} = \prod_{l=1}^r \EE e^{t s_l a_1[l]^2}
        = \prod_{l=1}^r (1-2 s_l t)^{-1/2}
    $ 
    for all $t<\frac{1}{2 \max_l s_l}$.
\end{proof}

We now lower-bound the expectation of $\min_{j \leq n} \zeta_j$
using a Chernoff bound-type argument,
which we note does not require independence.

\begin{lemma} \label{lm:RMT:spreadout_lb_zetaj}
    We have
    \begin{equation}
        \EE \min_{j \leq n} \zeta_j
        \geq \trace(S) \left(1 - 2\frac{\norm{S}_F}{\trace(S)} \sqrt{\log n} \right)
    \end{equation}
    where $\norm{\cdot}_F$ denotes Frobenius norm, i.e., $\ell_2$-norm of the vector of eigenvalues.
\end{lemma}

\begin{proof}
    By Jensen's inequality 
    and monotonicity of $\exp(\cdot)$,
    \begin{equation}
        \forall t>0,~
        \exp\left( t~ \EE \max_{j \leq n} -\zeta_j \right)
        \leq \EE \exp\left( t \max_{j \leq n} -\zeta_j \right) = \EE \max_{j \leq n} \exp \left( -t \zeta_j \right).
    \end{equation}
    So, taking $\log$ and optimizing the bound,
    \begin{align}
        \EE \max_{j \leq n} (-\zeta_j) 
        &\leq \inf_{t>0} \frac{1}{t} \log \EE 
        \max_{j \leq n} \exp\left( t (-\zeta_j) \right) \\
        &\leq \inf_{t>0} \frac{1}{t} \log\left[ n \EE e^{t (-\zeta_1)} \right]
        = \inf_{t>0} \frac{1}{t} \log\left[
            n
            \prod_{l=1}^r (1+2 s_l t)^{-1/2}
        \right] \\
        &= \inf_{t>0} \frac{1}{t} \left( \log(n) -\frac{1}{2} \sum_{l=1}^r \log(1+2 s_l t) \right)
        \eqqcolon \inf_{t>0} g(t).
    \label{eq:RMT:chernoff_bound_type_ineq}
    \end{align}
    
    By calculating we find that
    $g'(t)>0 \iff \sum_{l=1}^r \left[ \log(1+2s_l t) - \frac{2s_l t}{1+2s_l t} \right] > 2 \log n$.
    With the case $r \gg \log n$ in mind, let us evaluate at $t^{(1)}$ defined by
    $\sum_{l=1}^r 2 (s_l t)^2 = 2 \log n$, i.e.,
    $t^{(1)} = \sqrt{\log n} / \norm{S}_F$.
    (This choice is obtained by Taylor expansion for $t\to 0$ of the condition $g'(t)=0$.)
    Using that $\log(1+y) \geq y-\frac{1}{2}y^2$ for $y \geq 0$, we get
    \begin{align*}
        g(t^{(1)}) &= \frac{\norm{S}_F}{\sqrt{\log n}} \left( \log n - \frac{1}{2} \sum_{l=1}^r \log \left(1+2s_l \frac{\sqrt{\log n}}{\norm{S}_F} \right) \right) \\
        &\leq \frac{\norm{S}_F}{\sqrt{\log n}} \left( \log n - \frac{1}{2} \sum_{l=1}^r \left( 2s_l \frac{\sqrt{\log n}}{\norm{S}_F} - \frac{1}{2} 4s_l^2 \frac{\log n}{\norm{S}_F^2} \right) \right) \\
        &= \frac{\norm{S}_F}{\sqrt{\log n}} \left( \log n - \frac{\trace(S) \sqrt{\log n}}{\norm{S}_F} + \log n \right) \\
        &= -\trace(S) + 2 \norm{S}_F \sqrt{\log n}.
    \end{align*}
    Thus,
    $
        \EE \min_{j \leq n} \zeta_j
        \geq -g(t^{(1)}) 
        \geq \trace(S) - 2 \norm{S}_F \sqrt{\log n}
        = \trace(S) \left(1 - 2\frac{\norm{S}_F}{\trace(S)} \sqrt{\log n} \right)
    $.
\end{proof}

The upper bound of \propositionref{prop:quantitative:RMT_spreadout} is shown in \lemmaref{lm:RMT:reduce_to_gaussians},
and the lower bound follows immediately from substituting \lemmaref{lm:RMT:spreadout_lb_zetaj} into \lemmaref{lm:RMT:reduce_to_gaussians}.

\subsection{Proof of \texorpdfstring{\propositionref{prop:quantitative:RMT_sparse}}{Proposition~\ref{prop:quantitative:RMT_sparse}}}

For ease of reference, we restate the proposition below.
\propRMTsparse*

We are still exactly in the same setting as for \propositionref{prop:quantitative:RMT_spreadout}, so all the lemmas of \autoref{subsec:RMT_proofs:RMT_spreadout} apply.
We also reuse notations from that subsection:
$\mu = \min_{j \leq n} u_j^\top Q u_j + v_j^\top R v_j$,
$a_j = \begin{pmatrix} \sqrt{N_j}~ u_j \\ \sqrt{M_j}~ v_j \end{pmatrix}$ where
$N_1, ..., N_n, M_1, ..., M_n \sim \chi^2_n$ i.i.d.\ and independent of $U$ and $V$,
and $\zeta_j = a_j^\top S a_j$ for $j \leq n$.
We also assume $S$ diagonal w.l.o.g.

Recall from \lemmaref{lm:RMT:reduce_to_gaussians} that $\EE \mu \geq \frac{1}{n} \EE \min_{j \leq n} \zeta_j$.
The proposition then follows immediately from the following lemma.
\begin{lemma} \label{lm:RMT:sparse_lb_zetaj}
    %
    Let any $\SSS \subset \{1,...,r\}$, denote
    $r_\SSS = \lvert \SSS \rvert$ and $G_\SSS = \left[ \prod_{l\in \SSS} s_l \right]^{1/r_\SSS}$.
    Then we have 
    $\EE \min_{j \leq n} \zeta_j \geq \frac{1}{e} ~r_\SSS~ n^{-\frac{2}{r_\SSS}} G_\SSS$.
\end{lemma}

\begin{proof}
    Recall from \eqref{eq:RMT:chernoff_bound_type_ineq} that
    \begin{equation}
        \EE \max_{j \leq n} (-\zeta_j) 
        \leq \inf_{t>0} \frac{1}{t} \left( \log(n) -\frac{1}{2} \sum_{l=1}^r \log(1+2 s_l t) \right)
        \eqqcolon \inf_{t>0} g(t)
    \end{equation}
    and that 
    $g'(t)>0 \iff \sum_{l=1}^r \left[ \log(1+2s_l t) - \frac{2s_l t}{1+2s_l t} \right] > 2 \log n$.

    With the case $r \ll \log n$ in mind, let us evaluate at $t^{(2)}$ defined by
    $\sum_{l=1}^r \left[ \log(1+2s_l t) - 1 \right] = 2 \log n$,
    i.e., $\prod_{l=1}^r (1+2s_l t) = e^{2 \log n + r} = n^2 e^r$ --
    more precisely, let $t^{(2)}$ be the smallest positive such $t$ (since there may be several solutions to that polynomial equation).
    Then $g(t^{(2)}) = \frac{1}{t^{(2)}} (-r/2)$.
    
    Let us further upper-bound $t^{(2)}$ by some $u>0$.
    Since $t^{(2)}$ is defined as the smallest positive root of the polynomial
    $P(t) = \prod_{l=1}^r (1+2s_l t) - n^2 e^r$,
    and since $P(0) <0$,
    it suffices to pick any $u>0$ that $P(u)\geq 0$.
    Since $P(t) \geq \left[ \prod_{l=1}^r 2s_l \right] t^r - n^2 e^r$ for any $t>0$,%
    \footnote{
        And this approximation of $P(t)$ intuitively makes sense to do since our choice of $t^{(2)}$ was guided by the Ansatz that $\frac{2s_l t}{1+2s_l t} \approx 1$, i.e.\ $s_l t \gg 1$, for all $l \leq r$.
    }
    we may choose
    $u = (n^2 e^r)^{\frac{1}{r}} / \left[ \prod_{l=1}^r 2s_l \right]^{\frac{1}{r}}$.
    
    Thus we have
    \begin{align*}
        -\EE \min_{j \leq n} \zeta_j
        &\leq g(t^{(2)}) = \frac{1}{t^{(2)}} (-r/2)
        \leq \frac{1}{u} (-r/2)
        = \frac{\left[ \prod_{l=1}^r 2s_l \right]^{1/r}}{(n^2 e^r)^{1/r}} (-r/2) \\
        \EE \min_{j \leq n} \zeta_j
        &\geq r \frac{\left[ \prod_{l=1}^r s_l \right]^{1/r}}{e \cdot n^{2/r}}.
    \end{align*}
    This shows the lemma in the case $\SSS = \{1, ..., r\}$.

    For the case of arbitrary $\SSS \subset \{1, ..., r\}$,
    drop some terms in the Chernoff bound-type inequality; that is, upper-bound $-\sum_{l \not\in \SSS} \log(1+2s_l t)$ by zero:
    \begin{equation}
        -\EE \min_{j \leq n} \zeta_j
        \leq \inf_{t>0} \frac{1}{t} \left( \log(n) -\frac{1}{2} \sum_{l=1}^r \log(1+2 s_l t) \right)
        \leq \inf_{t>0} \underbrace{
            \frac{1}{t} \left( \log(n) -\frac{1}{2} \sum_{l \in \SSS} \log(1+2 s_l t) \right)
        }_{\eqqcolon g_\SSS(t)}.
    \end{equation}
    Thereafter, going through exactly the same considerations as above, we obtain the inequality where we restricted attention to the components $l \in \SSS$.
\end{proof}

\subsection{High-probability bounds for the spread-out spectrum case}
\label{subsec:RMT_proofs:RMTinprob_spreadout}

In this subsection we provide a high-probability counterpart to the expectation estimate 
from \autoref{prop:quantitative:RMT_spreadout},
where we showed
$\EE \left[ \min_{j \leq n} u_j^\top Q u_j + v_j^\top R v_j \right] \sim \frac{\trace(S)}{n}$
when $\frac{\trace(S)}{\norm{S}_F} \geq \Omega(\sqrt{\log n})$.

\begin{proposition} \label{prop:RMTinprob:spreadout}
    Suppose $Q, R$ are fixed and 
    $U, V$ are independently distributed uniformly on the set of $n \times n$ orthonormal matrices.
    Then,
    denoting $\mu = \min_{j \leq n} u_j^\top Q u_j + v_j^\top R v_j$,
    \begin{align*}
        \forall 0 \leq \gamma \leq 1,~~
        \PP\left( 
            \mu \geq \frac{\trace(S)}{2n} (1-\gamma)
        \right)
        & ~ \geq ~
        1 - n e^{-\frac{\trace(S)^2}{4 \norm{S}_F^2} \gamma^2} - 2 e^{-\frac{n}{8}} \\
        \text{and}~~~~
        \forall \gamma \geq 0,~~
        \PP\left(
            \mu \leq \frac{4 \trace(S)}{n} (1+\gamma)
        \right)
        & ~ \geq ~
        \begin{cases}
            1 - e^{-\frac{\trace(S)^2}{8 \norm{S}_F^2} \gamma^2} - 2 e^{-\frac{n}{8}} 
            & \text{if}~ \gamma \leq \frac{\norm{S}_F^2}{\trace(S) \spnorm{S}} \\
            1 - e^{-\frac{\trace(S)}{8 \spnorm{S}} \gamma} - 2 e^{-\frac{n}{8}}
            & \text{otherwise}
        \end{cases}
    \end{align*}
    where $\norm{\cdot}_F$ denotes the Frobenius norm
    and $\spnorm{\cdot}$ denotes the operator norm.
\end{proposition}

The remainder of this subsection is dedicated to proving the above proposition.
Let as in \autoref{subsec:RMT_proofs:RMT_spreadout}
$a_j = \begin{pmatrix} \sqrt{N_j}~ u_j \\ \sqrt{M_j}~ v_j \end{pmatrix}$ where
$N_1, ..., N_n, M_1, ..., M_n \sim \chi^2_n$ i.i.d.\ and independent of $U$ and $V$,
and $\zeta_j = a_j^\top S a_j$ for $j \leq n$.
Also denote $s_1 \geq ... \geq s_r > 0 = s_{r+1} = ... s_{2n}$ the eigenvalues of $S$, and assume $S$ diagonal w.l.o.g.

\begin{lemma} \label{lm:RMTinprob:inv_chi2_concentration_bounds}
    For each $j$,
    $\PP\left( \frac{1}{2n} \leq \frac{1}{N_j} \right) \geq 1-e^{-\frac{n}{8}}$
    and 
    $\PP\left( \frac{1}{N_j} \leq \frac{4}{n} \right) \geq 1-e^{-\frac{n}{8}}$.
\end{lemma}

\begin{proof}
    Since $N_j \sim \chi^2_n$,
    we have the classical concentration bounds
    \citep[Equations~(4.3), (4.4)]{laurent2000adaptive}
    \begin{equation} \label{eq:RMTinprob:chi2_concentration}
        \forall t > 0,~
        \PP\left( N_j - n \leq 2 \sqrt{n t} + 2 t \right) \geq 1- e^{-t}
        ~~~~\text{and}~~~~
        \PP\left( N_j - n \geq -2 \sqrt{n t} \right) \geq 1- e^{-t}.
    \end{equation}
    Evaluating the first inequality at $t = \frac{n}{8}$ yields
    $N_j \leq n + \frac{n}{\sqrt{2}} + \frac{n}{4} \leq 2n$
    with probability $\geq 1 - e^{-\frac{n}{8}}$.
    Evaluating the second inequality at $t = \frac{n}{8}$ yields
    $N_j \geq n - \frac{n}{\sqrt{2}} \geq \frac{n}{4}$
    with probability $\geq 1 - e^{-\frac{n}{8}}$.
\end{proof}

Letting the random variable $J = \argmin_{j \leq n} u_j^\top Q u_j + v_j^\top R v_j$ that is only dependent on $U$ and $V$, so independent of the $N_j$, $M_j$, 
we have
\begin{align*}
    \mu = \frac{1}{N_J} \cdot \sqrt{N_J} u_J^\top Q \sqrt{N_J} u_J + \frac{1}{M_J} \cdot \sqrt{M_J} v_J^\top R \sqrt{M_J} v_J
    &\geq 
    \left( \frac{1}{N_J} \wedge \frac{1}{M_J} \right) a_J^\top S a_J \\
    &\geq 
    \left( \frac{1}{N_J} \wedge \frac{1}{M_J} \right) \min_{j \leq n} a_j^\top S a_j
\end{align*}
and so by union bound, using that 
$N_J \sim N_1$ and $M_J \sim M_1$
by independence,
\begin{equation}
    \forall c \geq 0,~
    \PP\left( \mu \geq c \right)
    \geq \PP\left( \min_{j \leq n} a_j^\top S a_j \geq 2n c \right) - 2 e^{-\frac{n}{8}}.
\label{eq:RMTinprob:reduce_to_gaussians_lb}
\end{equation}

Let for concision $\zeta_j = a_j^\top S a_j$ for each $j \leq n$.
Recall from \autoref{lm:RMT:MGF_zetaj} that the moment-generating function of $\zeta_1 (\sim \zeta_2 \sim ... \sim \zeta_n)$ is
\begin{equation}
    \EE e^{t \zeta_1}
    = \prod_{l=1}^r (1-2 s_l t)^{-1/2}
    ~~~\text{for all}~~~
    t<\frac{1}{2 \max_l s_l}.
\end{equation}
We can now use union bound and Chernoff's bound to lower-bound $\min_{j \leq n} \zeta_j$ with high probability.
The derivation is essentially an instantiation of the general subexponential tail bound \cite[Sec.~2.1.3]{wainwright2019high} using our precise knowledge of the moment-generating function of the $\zeta_j$.

\begin{lemma} \label{lm:RMTinprob:spreadout_lb_zetaj}
    For any $0 \leq \gamma \leq 1$, we have
    \begin{equation}
        \PP\left( \min_{j \leq n} \zeta_j \geq \frac{\trace(S)}{2n} (1-\gamma) \right)
        \geq 1 - n e^{-\frac{\trace(S)^2}{4 \norm{S}_F^2} \gamma^2}
    \end{equation}
    where $\norm{\cdot}_F$ denotes Frobenius norm, i.e., $\ell_2$-norm of the vector of eigenvalues.
\end{lemma}

\begin{proof}
    By union bound, since the $\zeta_j$ are identically distributed, 
    $
        \PP\left( \min_j \zeta_j \leq x \right)
        \leq n \PP\left( \zeta_1 \leq x \right)
    $.
    By Markov's inequality,
    \begin{align*}
        \PP\left(\zeta_1 \leq x \right)
        &= \PP\left(-\zeta_1 \geq -x \right)
        \leq \inf_{t>0} \frac{1}{e^{-tx}} \EE e^{-t \zeta_1}
        = \inf_{t>0} e^{tx} \prod_{l=1}^r (1 + 2s_l t)^{-1/2} \\
        \log \PP\left(\zeta_1 \leq x \right)
        &\leq 
        \inf_{t>0}~
        tx - \frac{1}{2} \sum_{l=1}^r \log (1 + 2s_l t)
        \eqqcolon \inf_{t>0} g_x(t).
    \end{align*}
    By calculating we find that
    $g_x'(t) > 0 \iff 2x > \sum_{l=1}^r \frac{2s_l t}{1+2s_l t}$.
    With the case ${x \approx \EE \zeta_1 = \trace(S) \ll r}$ in mind,
    let us evaluate at $t^{(1)}_x$ defined by
    $2x = \sum_{l=1}^r 2s_l (1-2s_l t) = 2 \trace(S) - 4 \norm{S}_F^2 t$, i.e., 
    $t^{(1)}_x = \frac{\trace(S)-x}{2\norm{S}_F^2}$.
    Assume henceforth that $x \leq \trace(S)$ so that $t^{(1)}_x \geq 0$.
    Using that $\log(1+y) \geq y - \frac{1}{2} y^2$ for $y \geq 0$, we get
    \begin{align*}
        g_x(t^{(1)}_x)
        &= \frac{\trace(S) - x}{2 \norm{S}_F^2} x 
        - \frac{1}{2} \sum_{l=1}^r \log \left( 1 + 2s_l \frac{\trace(S)-x}{2\norm{S}_F^2} \right) \\
        &\leq \frac{\trace(S) x - x^2}{2 \norm{S}_F^2} 
        - \frac{1}{2} \sum_{l=1}^r 2s_l \frac{\trace(S)-x}{2\norm{S}_F^2} 
        + \frac{1}{4} \sum_{l=1}^r 4 s_l^2 \left(\frac{\trace(S)-x}{2\norm{S}_F^2}\right)^2 \\
        &= \frac{\trace(S) x - x^2}{2 \norm{S}_F^2} 
        - \trace(S) \frac{\trace(S)-x}{2\norm{S}_F^2} 
        + \norm{S}_F^2 \left(\frac{\trace(S)-x}{2\norm{S}_F^2}\right)^2 \\
        &= \frac{1}{2\norm{S}_F^2} \left[
            \trace(S) x - x^2 - \trace(S)^2 + \trace(S) x + \frac{1}{2} \left( \trace(S) - x \right)^2 
        \right] \\
        &= -\frac{1}{4 \norm{S}_F^2} [\trace(S) - x]^2.
    \end{align*}
    Hence we have shown
    \begin{align*}
        \forall x \leq \trace(S),~
        \PP\left( \min_j \zeta_j \leq x \right)
        \leq n \PP\left( \zeta_1 \leq x \right)
        \leq n \exp\left( -\frac{1}{4 \norm{S}_F^2} [\trace(S) - x]^2 \right)
    \end{align*}
    and the announced bound follows by a change of variables 
    $\gamma = 1 - \frac{x}{\trace(S)}$.
\end{proof}


Combining \eqref{eq:RMTinprob:reduce_to_gaussians_lb} with \autoref{lm:RMTinprob:spreadout_lb_zetaj} yields the high-probability lower bound 
\begin{align*}
    \forall 0 \leq \gamma \leq 1,~~~~
    \PP\left( \mu \geq \frac{\trace(S)}{2n} (1-\gamma) \right)
    &\geq 1 - n e^{-\frac{\trace(S)^2}{4 \norm{S}_F^2} \gamma^2} - 2 e^{-\frac{n}{8}}
\end{align*}
which is exactly the lower bound of \autoref{prop:RMTinprob:spreadout}.

For a high-probability upper bound on $\mu = \min_{j \leq n} w_j^\top S w_j$, in the regime of interest in this subsection it is sufficient to start from
\begin{equation}
    \mu \leq w_1^\top S w_1
    = \frac{1}{N_1} \cdot \sqrt{N_1} u_1^\top Q \sqrt{N_1} u_1 + \frac{1}{M_1} \cdot \sqrt{M_1} v_1^\top R \sqrt{M_1} v_1
    \leq \left( \frac{1}{N_1} \vee \frac{1}{M_1} \right) a_1^\top S a_1
\end{equation}
so that, by union bound and \autoref{lm:RMTinprob:inv_chi2_concentration_bounds},
\begin{equation}
    \forall C \geq 0,~
    \PP(\mu \leq C)
    \geq \PP\left( a_1^\top S a_1 \leq \frac{n}{4} C \right) - 2 e^{-\frac{n}{8}}.
\label{eq:RMTinprob:reduce_to_gaussians_ub}
\end{equation}
Hence it is sufficient to upper-bound $\zeta_1 = a_1^\top S a_1$ with high probability.

\begin{lemma} \label{lm:RMTinprob:spreadout_ub_zetaj}
    For any $\eps \geq 0$, we have
    \begin{equation}
        \PP\left( \zeta_1 \leq \trace(S) + \eps \right)
        \geq \begin{cases}
            1 - e^{-\frac{\eps^2}{8 \norm{S}_F^2}} & \text{if}~ 0 \leq \eps \leq \frac{\norm{S}_F^2}{\max_l s_l} \\
            1 - e^{-\frac{\eps}{8 \max_l s_l}} & \text{otherwise}.
        \end{cases}
    \end{equation}
\end{lemma}

\begin{proof}
    $\zeta_1 = \sum_{l=1}^r s_l a_1[l]^2$ is subexponential with parameters $\left( \sqrt{\sum_{l=1}^r (2 \cdot s_l)^2}, \max_{l \leq r} 4 \cdot s_l \right)$ as a linear combinations of independent $\chi^2$ random variables $a_j[l]^2$ which are subexponential with parameters $(2, 4)$ \cite[Sec.~2.1.3]{wainwright2019high}.
    The announced subexponential tail bound follows by a direct application of \cite[Prop.~2.9]{wainwright2019high}.
\end{proof}

Combining \eqref{eq:RMTinprob:reduce_to_gaussians_ub} with \autoref{lm:RMTinprob:spreadout_ub_zetaj} yields the high-probability upper bound 
\begin{align*}
    \forall \eps \geq 0,~~~~
    \PP\left( \mu \leq 4 \frac{\trace(S) + \eps}{n} \right) 
    &\geq \begin{cases}
        1 - e^{-\frac{\eps^2}{8 \norm{S}_F^2}} - 2 e^{-\frac{n}{8}} 
        & \text{if}~ \eps \leq \frac{\norm{S}_F^2}{\max_l s_l} \\
        1 - e^{-\frac{\eps}{8 \max_l s_l}} - 2 e^{-\frac{n}{8}}
        & \text{otherwise}
    \end{cases} \\
    \text{and so} \quad
    \forall \gamma \geq 0,~~~~
    \PP\left( \mu \leq \frac{4 \trace(S)}{n} (1+\gamma) \right) 
    &\geq \begin{cases}
        1 - e^{-\frac{\trace(S)^2}{8 \norm{S}_F^2} \gamma^2} - 2 e^{-\frac{n}{8}} 
        & \text{if}~ \gamma \leq \frac{\norm{S}_F^2}{\trace(S) \left[ \max_l s_l \right]} \\
        1 - e^{-\frac{\trace(S) \gamma}{8 \max_l s_l}} - 2 e^{-\frac{n}{8}}
        & \text{otherwise}
    \end{cases}
\end{align*}
by the change of variables $\trace(S) \gamma = \eps$, 
which is exactly the upper bound of \autoref{prop:RMTinprob:spreadout}.

\subsection{High-probability bounds for the sparse spectrum case}

In this subsection we provide a high-probability counterpart to the expectation lower estimate 
from \autoref{prop:quantitative:RMT_sparse},
where we showed
$\EE \left[ \min_{j \leq n} u_j^\top Q u_j + v_j^\top R v_j \right] \gtrsim n^{-\frac{2}{r}-1}$
when $n \to \infty$ with $r = \rank(S)$ and $s \in \RR^r$ fixed.

\begin{proposition} \label{prop:RMTinprob:sparse}
    Suppose $Q, R$ are fixed and 
    $U, V$ are independently distributed uniformly on the set of $n \times n$ orthonormal matrices.
    Let $s_1 \geq ... \geq s_r > 0 = s_{r+1} = ... = s_{2n}$ the eigenvalues of $S$.
    Then 
    denoting $\mu = \min_{j \leq n} u_j^\top Q u_j + v_j^\top R v_j$,
    for any $\SSS \subset \{1, ..., r\}$,
    \begin{equation}
        \forall 0 \leq \gamma \leq 1,~~
        \PP\left( 
            \mu \geq 
            \frac{1}{2e}
            \frac{\card{\SSS}}{n} n^{-\frac{2}{\card{\SSS}}}
            \left[\prod_{l \in \SSS} s_l \right]^{\frac{1}{\card{\SSS}}} 
            \cdot (1-\gamma)
        \right)
        ~ \geq ~
        1 - e^{-\frac{\card{\SSS}}{2} \gamma}
        - 2 e^{-\frac{n}{8}}.
    \end{equation}
\end{proposition}

The remainder of this subsection is dedicated to proving the above proposition.
As we are exactly in the same setting as for \propositionref{prop:RMTinprob:spreadout}, all the lemmas of \autoref{subsec:RMT_proofs:RMTinprob_spreadout} apply.
We also reuse notations from that subsection: 
$a_j = \begin{pmatrix} \sqrt{N_j}~ u_j \\ \sqrt{M_j}~ v_j \end{pmatrix}$ where
$N_1, ..., N_n, M_1, ..., M_n \sim \chi^2_n$ i.i.d.\ and independent of $U$ and $V$,
and $\zeta_j = a_j^\top S a_j$ for $j \leq n$.
We also assume $S$ diagonal w.l.o.g.

Given 
\eqref{eq:RMTinprob:reduce_to_gaussians_lb},
in order to prove the proposition it suffices to derive a different lower tail bound for $\min_{j \leq n} \zeta_j$ that is more adapted to the sparse spectrum case.

\begin{lemma} \label{lm:RMTinprob:sparse_lb_zetaj}
    Let any $\SSS \subset \{1,...,r\}$,
    $r_\SSS = \lvert \SSS \rvert$ and $G_\SSS = \left[ \prod_{l\in \SSS} s_l \right]^{1/r_\SSS}$. 
    Then
    \begin{equation}
        \forall \gamma \in \RR,~~
        \PP \left( 
            \min_{j \leq n} \zeta_j 
            \geq \frac{r_\SSS}{e} n^{-\frac{2}{r_\SSS}} \left[\prod_{l \in \SSS} s_l \right]^{\frac{1}{r_\SSS}} \cdot (1-\gamma)
        \right)
        \geq 1 - e^{-\frac{r_\SSS}{2} \gamma}.
    \end{equation}
\end{lemma}


\begin{proof}
    Let any $x \in \RR$.
    By Markov's inequality,
    \begin{equation}
        \PP\left( \min_j \zeta_j \leq x \right)
        = \PP\left( \max_j (-\zeta_j) \geq -x \right)
        \leq \inf_{t>0} \frac{1}{e^{-tx}} \EE e^{t \max_j (-\zeta_j)}
        = \inf_{t>0} e^{tx} \cdot \EE \max_j e^{t (-\zeta_j)}.
    \end{equation}
    Now for any $\alpha \geq 1$, by monotonicity and concavity of $y \mapsto y^{1/\alpha}$,
    \begin{align*}
        \EE \max_{j \leq n} e^{t (-\zeta_j)}
        = \EE\left[ \left( \max_{j \leq n} e^{\alpha t (-\zeta_j)} \right)^{1/\alpha} \right]
        \leq \EE\left[ \left( \sum_{j \leq n} e^{\alpha t (-\zeta_j)} \right)^{1/\alpha} \right] 
        &\leq \left[ \EE \sum_{j \leq n} e^{\alpha t (-\zeta_j)} \right]^{1/\alpha} \\
        &= \left[ n \EE e^{\alpha t (-\zeta_1)} \right]^{1/\alpha}.
    \end{align*}
    Hence, taking the infimum over $\alpha$
    and substituting into the above inequality,%
    \footnote{In fact one can check, by inverting $\inf_t$ and $\inf_\alpha$ and optimizing first over $t$, that the infimum over the joint $(t, \alpha)$ is always attained at some $(t^*, 1)$; that is, the true optimum is always at $\alpha=1$. What introducing the extra degree of freedom $\alpha$ affords us is a way to obtain a tractable upper bound, whereas upper-bounding the $\inf_t$ with $\alpha=1$ fixed seems more difficult.}
    \begin{align*}
        \PP\left( \min_j \zeta_j \leq x \right)
        &\leq \inf_{t>0} e^{tx} \cdot \inf_{\alpha \geq 1} \left[ n \EE e^{\alpha t(-\zeta_1)} \right]^{1/\alpha} \\
        \log \PP\left( \min_j \zeta_j \leq x \right)
        &\leq \inf_{t>0} tx + \inf_{\alpha \geq 1} \frac{1}{\alpha} \left( \log n + \log \EE e^{\alpha t(-\zeta_1)} \right) \\
        &= \inf_{t>0} t \left[ x + \inf_{\beta \geq t} \frac{1}{\beta} \left( \log n + \log \EE e^{\beta (-\zeta_1)} \right) \right].
    \end{align*}
    Let, as in the proof of \autoref{lm:RMT:sparse_lb_zetaj},
    \begin{equation}
        g(\beta) \coloneqq 
        \frac{1}{\beta} \left( \log n + \log \EE e^{-\beta \zeta_1} \right)
        = \frac{1}{\beta} \left( \log n - \frac{1}{2} \sum_{l=1}^r \log(1+2 s_l \beta) \right)
        \quad \text{and} \quad
        u \coloneqq 
        \frac{e}{2} n^{\frac{2}{r}} \left[ \prod_{l=1}^r s_l \right]^{-\frac{1}{r}}.
    \end{equation}
    One can check, by lower-bounding $\log(1+2s_l u)$ by $\log(2s_l u)$ for each term,
    that
    $g(u) \leq \frac{1}{u} (-r/2)$.
    Hence,
    \begin{align*}
        \log \PP\left( \min_j \zeta_j \leq x \right)
        &\leq \inf_{t>0} t \left[ x + \inf_{\beta \geq t} g(\beta) \right]
        \leq \inf_{0 < t \leq u} t \left( x + g(u) \right) \\
        &\leq u (x + g(u))
        \leq u x - \frac{r}{2}
        = \frac{r}{2} \cdot \left( \frac{u}{r/2} x - 1 \right).
    \end{align*}
    Hence we have shown that for all $x \in \RR$,
    \begin{equation}
        \PP\left( \min_j \zeta_j \leq x \right) \leq \exp\left(
            - \frac{r}{2} \left( 1 - \frac{u}{r/2} x \right)
        \right).
    \end{equation}
    The inequality of the lemma with $\SSS = \{1, ..., r\}$ follows
    by the change of variable 
    $\gamma = 1 - \frac{u}{r/2} x$
    and substituting the expression of $u$.

    The inequality for arbitrary $\SSS \subset \{1, ..., r\}$
    follows by exactly the same considerations as above applied to
    \begin{equation}
        g_\SSS(\beta) = 
        \frac{1}{\beta} \left( \log n - \frac{1}{2} \sum_{l \in \SSS} \log(1+2 s_l \beta) \right)
        \quad \text{and} \quad
        u_\SSS =
        \frac{e}{2} n^{\frac{2}{r_{\SSS}}} \left[ \prod_{l \in \SSS} s_l \right]^{-\frac{1}{r_{\SSS}}}
    \end{equation}
    instead of $g$ and $u$, noting that $g(\beta) \leq g_\SSS(\beta)$ for any $\beta>0$.
\end{proof}

Combining \eqref{eq:RMTinprob:reduce_to_gaussians_lb} with \autoref{lm:RMTinprob:sparse_lb_zetaj} yields the high-probability lower bound
\begin{equation}
    \forall \gamma \in \RR,~~
    \PP\left( 
        \mu \geq 
        \frac{1}{2n}
        \frac{\card{\SSS}}{e} n^{-\frac{2}{\card{\SSS}}} \left[\prod_{l \in \SSS} s_l \right]^{\frac{1}{\card{\SSS}}} \cdot (1-\gamma)
    \right)
    \geq 
    1 - e^{-\frac{\card{\SSS}}{2} \gamma}
    - 2 e^{-\frac{n}{8}}
\end{equation}
for any $\SSS \subset \{1, ..., r\}$,
which is exactly the inequality of \autoref{prop:RMTinprob:spreadout}.



\section{A symmetric expression for Alt-GDA} \label{apx:sym_altGDA}

This appendix contains results stated in \autoref{subsec:quantitative:DT_algos} and used in \autoref{subsec:DT_expan:altGDA} and that may be of independent interest.

As announced in \autoref{subsec:quantitative:DT_algos}, we analyze a symmetrized formulation of Alt-GDA, 
whose definition we recall here for ease of reference:
let $(x^0, y^\half) \in \RR^d \times \RR^d$ and
\begin{equation} \label{eq:sym_altGDA:altGDA_upd}
    \begin{cases}
        \forall k \in \NN,
        & x^{k+1} = x^k - \eta \nabla_x f(x^k, y^{k+\half}) \\
        \forall k \in \NN + \half,
        & y^{k+1} = y^k + \eta \nabla_y f(x^{k+\half}, y^k)
    \end{cases}
    ~~~~\text{and}~~~~
    \begin{cases}
        \forall k \in \NN,
        & x^{k+\half} = \frac{x^{k+1} + x^k}{2} \\
        \forall k \in \NN+\half,
        & y^{k+\half} = \frac{y^{k+1} + y^k}{2}.
    \end{cases}
\end{equation}
Note that this symmetrized formulation is indeed equivalent to standard one (used e.g.\ in~\cite{gidel_negative_2019, bailey2020finite, zhang_near-optimal_2022, grimmer2022limiting}), via the correspondence
$\forall k \in \NN,$ $\tx^k = x^k$ and $\ty^k = y^{k+\half}$.
We emphasize that, as one can easily check from the definition, the symmetrized iterates $(z^k)_k = (x^k, y^k)_k$ converge if and only if the classical ones $(\tz^k)_k = (x^k, y^{k+\half})_{k \in \NN}$ do, and if they converge exponentially, then they have the same convergence rate.

The following proposition shows that the evolution of the symmetrized Alt-GDA iterates can be expressed approximately in terms only of the (skewed) gradient field $g(z)$.

\begin{proposition} \label{prop:sym_altGDA:altGDA_upd_withO}
    Assume $\eta \leq \norm{\nabla_{xx}^2 f}_\infty^{-1} \wedge \norm{\nabla_{yy}^2 f}_\infty^{-1}$.
    Then the symmetrized Alt-GDA iterates defined by 
    \eqref{eq:sym_altGDA:altGDA_upd}
    satisfy
    \begin{equation}
        \forall k \in \NN \cup (\NN + \half),~
        z^{k+1} = z^k - \eta g(z^k) + \frac{\eta^2}{2} A(z^k) g(z^k)
        + O\left( \eta^3 \norm{g(z^k)} \right)
    \end{equation}
    where 
    $A(z) = 
    \begin{bmatrix}
        0 & \nabla_{xy}^2 f \\
        -\nabla_{xy}^2 f^\top & 0
    \end{bmatrix}(z)$.
    More precisely the ``$O(\cdot)$'' term is a vector with norm at most
    $10 \eta^3 \norm{g(z^k)} \left( 
        \norm{\nabla^3 f}_\infty 
        (1 + \norm{\nabla_{xy}^2 f}_\infty^2 \eta^2) 
        + \norm{\nabla^2 f}_\infty \norm{g(z^k)} 
    \right)$.
\end{proposition}

\begin{remark}[High-resolution ODE for the symmetrized iterates]
    The symmetrized update operator with the $O(\eta^3)$ term neglected,
    $z \mapsto z - \eta g(z) + \frac{\eta^2}{2} A(z) g(z)$, leads to a different high-resolution ODE than the one derived by \cite[Eq.~(12)]{grimmer2022limiting}
    for the standard formulation.
    The one derived in that reference is
    \begin{equation}
        \frac{d\tz}{dt} = -g(\tz) - \frac{\eta}{2} \begin{bmatrix}
            Q & P \\
            P^\top & -R
        \end{bmatrix}(\tz)
        g(\tz)
        \qquad ~~\text{where}~~
        \begin{bmatrix}
            Q & P \\
            -P^\top & R
        \end{bmatrix}(z)
        \coloneqq
        \nabla g(z)
    \end{equation}
    while the high-resolution ODE corresponding to the symmetrized update $z^{k+1} \approx T(z^k)$ can be shown -- simply by applying the high-resolution ODE for explicit Euler steps \cite[Eq.~(20)]{lu_osr-resolution_2022} to the vector field $-g(z) + \frac{\eta}{2} A(z) g(z)$ -- to be
    \begin{align*}
        \frac{dz}{dt} = -g(z) - \frac{\eta}{2} M(z) g(z) + \frac{\eta}{2} A(z) g(z)
        &= -g(z) - \frac{\eta}{2} 
        \begin{bmatrix}
            Q & 0 \\
            0 & R
        \end{bmatrix}(z)
        g(z).
    \end{align*}
    There is no contradiction here, as the classical Alt-GDA iterates $\tz^k$ are not identical pointwise to the symmetrized ones $z^k$; instead the correspondence involves a time-shift of $+\half$ for one of the variables.
    Informally, the former high-resolution ODE tracks the evolution of the Alt-GDA iterates using a piecewise-constant (forward in time) interpolation scheme, with the arbitrary choice of measuring at time-steps where $y$ was updated more recently than $x$, and the latter tracks the same iterates but using a piecewise-linear (trapezoidal) interpolation scheme. 
\end{remark}

\begin{proof}
    As usual in this paper we let for concision $z^k = (x^k, y^k)$ for all $k \in \NN \cup (\NN + \half)$.
    Denote $L_{xx} = \norm{\nabla_{xx}^2 f}_\infty$, $L_{yy} = \norm{\nabla_{yy}^2 f}_\infty$, $L_{xy} = \norm{\nabla_{xy}^2 f}_\infty$, $L_2 = L_{xx} \vee L_{yy} \vee L_{xy}$ and $L_3 = \norm{\nabla^3 f}_\infty$.
    Throughout this proof we will write $\beps$ or $\beps'$ to denote elements of $[-1, 1]$ or of the unit ball of $\RR^n$ or $\RR^m$ or $\RR^{n+m}$, and that may change from line to line.

    Let $k \in \NN$. By Taylor expansion with remainder in Lagrange form,
    \begin{align*}
        x^{k+1} -x^k &= - \eta \nabla_x f(x^k, y^{k+\half}) \\
        &= - \eta \left[ \nabla_x f(z^k) + \nabla_{xy}^2 f(z^k) \cdot (y^{k+\half} - y^k) + \beps \left( \frac{L_3}{2} \norm{y^{k+\half}-y^k}^2 \right) \right].
    \end{align*}
    Now 
    \begin{align*}
        y^{k+\half} - y^k
        &= y^{k+\half} - \frac{y^{k+\half} + y^{k-\half}}{2} = \frac{y^{k+\half} - y^{k-\half}}{2} \\
        \text{and}~~~~
        y^{k+\half} - y^{k-\half}
        &= \eta \nabla_y f(x^k, y^{k-\half}) 
        = \eta \nabla_y f(z^k) + \beps \left( \eta L_{yy} \norm{y^k - y^{k-\half}} \right) \\
        &= \eta \nabla_y f(z^k) + \beps \left( \eta L_{yy} \frac{1}{2} \norm{y^{k+\half} - y^{k-\half}} \right) \\
        \text{since}~~~~
        y^{k-\half} - y^k &= y^{k-\half} - \frac{y^{k+\half} + y^{k-\half}}{2}
        = -\frac{y^{k+\half}-y^{k-\half}}{2}
        = -\left( y^{k+\half} - y^k \right).
    \end{align*}
    Moreover, the above computation also shows that
    \begin{align}
        \norm{y^{k+\half}-y^{k-\half}} 
        &\leq \eta \norm{\nabla_y f(z^k)} + \eta L_{yy} \frac{1}{2} \norm{y^{k+\half}-y^{k-\half}} \\
        \text{and so}~~~~
        \norm{y^{k+\half}-y^{k-\half}} &\leq \frac{1}{1 - \eta L_{yy}/2} \eta \norm{\nabla_y f(z^k)}
    \label{eq:sym_altGDA:altGDA_upd:control_y}
    \end{align}
    for any $\eta \leq 2/L_{yy}$.
    In particular, since we assume $\eta \leq 1/L_{yy}$ then
    $\norm{y^{k+\half}-y^{k-\half}} \leq 2 \eta \norm{\nabla_y f(z^k)}$.
    Thus, using this last fact to express the error terms in terms of $\norm{g(z^k)} = \norm{\nabla f(z^k)}$, by substituting into the Taylor expansion of $x^{k+1}-x^k$ we obtain
    \begin{multline}
        \forall k \in \NN,~
        x^{k+1} - x^k = - \eta \nabla_x f(z^k) - \frac{1}{2} \eta^2 \nabla_{xy}^2 f(z^k) \cdot \nabla_y f(z^k) \\
        \qquad\qquad
        + \beps \left( \eta^2 L_{xy} L_{yy} \frac{1}{2} \norm{y^{k+\half}-y^{k-\half}} \right)
        + \beps' \left( \eta \frac{L_3}{8} \norm{y^{k+\half}-y^{k-\half}}^2 \right) \\
        = - \eta \nabla_x f(z^k) - \frac{1}{2} \eta^2 \nabla_{xy}^2 f(z^k) \cdot \nabla_y f(z^k)
        + \eta^3 \beps \left( 
            L_{xy} L_{yy} \norm{g(z^k)} 
            + \frac{L_3}{2} \norm{g(z^k)}^2
        \right).
    \label{eq:sym_altGDA:altGDA_upd:x_updating}
    \end{multline}
    By the symmetric calculations, we have 
    \begin{multline}
        \forall k \in \NN+\half,~
        y^{k+1} - y^k = \eta \nabla_y f(z^k) - \frac{1}{2} \eta^2 \nabla_{yx}^2 f(z^k) \nabla_x f(z^k) \\
        + \eta^3 \beps \left( 
            L_{xy} L_{yy} \norm{g(z^k)} 
            + \frac{L_3}{2} \norm{g(z^k)}^2
        \right).
    \label{eq:sym_altGDA:altGDA_upd:y_updating}
    \end{multline}
    
    Now let us compute an expansion for the increments between non-updating time-steps. Let $k \in \NN$ and let us compute $y^{k+1}-y^k$:
    \begin{align}
        & y^{k+1} - y^k = \frac{y^{k+\frac{3}{2}} + y^{k+\half}}{2} - \frac{y^{k+\half}+y^{k-\half}}{2} 
        = \frac{1}{2} \left( y^{k+\frac{3}{2}} - y^{k-\half} \right) \\
        &= \frac{1}{2} \left( \eta \nabla_y f\left( x^{k+1}, y^{k+\half} \right) + \eta \nabla_y f\left( x^k, y^{k-\half} \right) \right) \\
        &= \frac{1}{2} \eta \left[ 
            2 \nabla_y f(z^k)
            + \nabla_{yx}^2 f(z^k) \cdot (x^{k+1}-x^k)
            + \nabla_{yy}^2 f(z^k) \cdot (y^{k+\half} - y^k)
            + \nabla_{yy}^2 f(z^k) \cdot (y^{k-\half} - y^k)
        \right] \\
        &\qquad\qquad
        + \beps \left(
            \eta \frac{L_3}{2} \norm{x^{k+1}-x^k}^2 + \eta \frac{L_3}{2} \norm{y^{k+\half}-y^k}^2
        \right) 
        + \beps' \left(
            \eta \frac{L_3}{2} \norm{y^k-y^{k-\half}}^2
        \right) \\
        &= \eta \nabla_y f(z^k) 
        + \frac{1}{2} \eta \nabla_{yx}^2 f(z^k) \cdot (x^{k+1}-x^k) 
        + \frac{1}{2} \eta \nabla_{yy}^2 f(z^k) \cdot \underbrace{
            \left( y^{k+\half} + y^{k-\half} - 2 y^k \right)
        }_{=0} \\
        &\qquad\qquad
        + \beps \left( 
            \eta \frac{L_3}{2} \norm{x^{k+1}-x^k}^2 + \eta \frac{L_3}{4} \norm{y^{k+\half}-y^{k-\half}}^2 
        \right) \\
        &= \eta \nabla_y f(z^k) 
        + \frac{1}{2} \eta \nabla_{yx}^2 f(z^k) \cdot (x^{k+1}-x^k) 
        + \beps \left( 
            \eta^3 L_3 \left( 2 + L_{xy}^2 \eta^2 \right) \norm{g(z^k)}^2
        \right).
    \end{align}
    Here in order to bound the error term in $\norm{x^{k+1}-x^k}$ we used that
    \begin{align}
        \norm{x^{k+1}-x^k} &= \eta \norm{\nabla_x f(x^k, y^{k+\half})} 
        \leq \eta \norm{\nabla_x f(z^k)} + \eta L_{xy} \norm{y^{k+\half} - y^k} \\
        &= \eta \norm{\nabla_x f(z^k)} + \frac{1}{2} \eta L_{xy} \norm{y^{k+\half} - y^{k-\half}}
        \leq \eta (1 + L_{xy} \eta) \norm{g(z^k)} \\
        \norm{x^{k+1}-x^k} &\leq \eta^2 \left( 2 + 2 L_{xy}^2 \eta^2 \right) \norm{g(z^k)}^2
    \end{align}
    since
    $\norm{y^{k+\half}-y^{k-\half}} \leq 2 \eta \norm{\nabla_y f(z^k)}$
    as shown previously \eqref{eq:sym_altGDA:altGDA_upd:control_y}.
    Hence, substituting the expansion of $x^{k+1}-x^k$ from the previous paragraph, we have
    \begin{multline}
        \forall k \in \NN,~
        y^{k+1} - y^k
        = \eta \nabla_y f(z^k) 
        + \frac{1}{2} \eta \nabla_{yx}^2 f(z^k) \cdot \left[
            - \eta \nabla_x f(z^k) - \frac{1}{2} \eta^2 \nabla_{xy}^2 f(z^k) \cdot \nabla_y f(z^k) 
        \right] \\
        \qquad\qquad
        + \eta^3 \beps \left( 
            L_3 \left( 2 + L_{xy}^2 \eta^2 \right) \norm{g(z^k)}^2
            + 
            L_{xy} L_{yy} \norm{g(z^k)} + \frac{L_3}{2} \norm{g(z^k)}^2
        \right) \\
        = \eta \nabla_y f(z^k) 
        - \frac{1}{2} \eta^2 \nabla_{yx}^2 f(z^k) \cdot \nabla_x f(z^k)
        + \eta^3 \beps \left(
            L_3 \left( \frac{5}{2} + L_{xy}^2 \eta^2 \right) \norm{g(z^k)}^2
            + \frac{5}{4} L_2^2 \norm{g(z^k)}
        \right).
    \label{eq:sym_altGDA:altGDA_upd:x_nonupdating}
    \end{multline}
    By the symmetric calculations, we have
    \begin{multline}
        \forall k \in \NN+\half,~
        x^{k+1} - x^k
        = -\eta \nabla_x f(z^k) - \frac{1}{2} \eta^2 \nabla_{xy}^2 f(z^k) \cdot \nabla_y f(z^k) \\
        + \eta^3 \beps \left( 
            L_3 \left( \frac{5}{2} + L_{xy}^2 \eta^2 \right) \norm{g(z^k)}^2
            + \frac{5}{4} L_2^2 \norm{g(z^k)}
        \right).
    \label{eq:sym_altGDA:altGDA_upd:y_nonupdating}
    \end{multline}
    
    The announced expansion for $z^{k+1} - z^k$, both at time-steps $k \in \NN$ and $k \in \NN + \half$, 
    follows immediately from 
    \eqref{eq:sym_altGDA:altGDA_upd:x_updating},
    \eqref{eq:sym_altGDA:altGDA_upd:y_updating},
    \eqref{eq:sym_altGDA:altGDA_upd:x_nonupdating} 
    and \eqref{eq:sym_altGDA:altGDA_upd:y_nonupdating}.
\end{proof}

The above \autoref{prop:sym_altGDA:altGDA_upd_withO} is already a good indication of the fact that local convergence of Alt-GDA can be analyzed via the spectrum of $\nabla\left[ \id - \eta g + \frac{\eta^2}{2} A \cdot g \right](z^*) = I-\eta M + \frac{\eta^2}{2} A M$.
In fact, since the error term in the above proposition is $O(\eta^3 \norm{g(z^k)})$, one can follow the potential-based approach developed in \cite[Sec.~4]{lu_osr-resolution_2022} (via continuous-time)
or \cite[Sec.~H]{wang_exponentially_2022} (for a special case),
to lower-bound the decrease of $\norm{g(z^k)}$ or $\norm{z^k-z^*}$ at each iteration.
However that approach may require $\eta$ to be smaller than actually necessary.
The following proposition offers a more direct path to linking convergence of Alt-GDA to spectral properties of the aforementioned matrix.

\begin{proposition} \label{prop:sym_altGDA:altGDA_upd_withexplicitThalf}
    Assume $\eta \leq \norm{\nabla_{xx}^2 f}_\infty^{-1} \wedge \norm{\nabla_{yy}^2 f}_\infty^{-1}$.
    Consider $z^k = (x^k, y^k)_{k \in \NN \cup (\NN + \half)}$ the symmetrized Alt-GDA iterates defined by 
    \eqref{eq:sym_altGDA:altGDA_upd}.
    
    For any $k \in \NN \cup (\NN + \half)$, $z^{k+\half}$ is entirely determined by $z^k$.
    More precisely, there exist well-defined operators $T^\half_{xy}, T^\half_{yx}$ such that
    $\begin{cases}
        \forall k \in \NN, & z^{k+\half} = T^\half_{xy}(z^k) \\
        \forall k \in \NN + \half, & z^{k+\half} = T^\half_{yx}(z^k)
    \end{cases}$.
    
    Moreover, for $z^*$ such that $g(z^*) = 0$, it holds
    $T^\half_{xy}(z^*) = T^\half_{yx}(z^*) = z^*$ and
    \begin{equation}
        \nabla \left[ T^\half_{xy} \circ T^\half_{yx} \right](z^*)
        ,~~
        \nabla \left[ T^\half_{yx} \circ T^\half_{xy} \right](z^*)
        = I - \eta M + \frac{\eta^2}{2} A M
        + O(\eta^3).
    \end{equation}
    More precisely the ``$O(\cdot)$'' term is a matrix with operator norm -- and so spectral radius -- at most
    $2 \eta^3 \spnorm{A} \left( \spnorm{A} \vee \spnorm{S} \right)^2$, where $\spnorm{\cdot}$ denotes operator norm.
\end{proposition}

Thanks to the above proposition, one can analyze the local convergence of the integer-time-step iterates $(z^k)_{k \in \NN}$, say, by applying the usual analysis to its well-defined update operator $T^\half_{yx} \circ T^\half_{xy}$;
that is to say, by bounding the spectral radius of its Jacobian at $z^*$.

For comparison, the update operator $\tT_{xy}$ for the standard formulation is
given by
$\tz^{k+1} = \begin{pmatrix}
    \tx^k - \nabla_x f(\tx^k, \ty^k) \\
    \ty^k + \nabla_y f(\tx^{k+1}, \ty^k)
\end{pmatrix}$.
Denoting
$M = \begin{bmatrix}
    Q & P \\ -P^\top & R
\end{bmatrix}$
as usual the Jacobian of $g$ at a fixed point $z^*$,
the Jacobian of $\tT_{xy}$ at $z^*$ writes
\cite[Sec.~A.3]{zhang_near-optimal_2022}
\begin{equation}
    \nabla \tT_{xy}(z^*) 
    = \begin{bmatrix}
        I & 0 \\
        \eta P^\top & I - \eta R
    \end{bmatrix}
    \begin{bmatrix}
        I - \eta Q & -\eta P \\
        0 & I
    \end{bmatrix} 
    = \begin{bmatrix}
        I - \eta Q & -\eta P \\
        \eta P^\top - \eta^2 P^\top Q & I-\eta R - \eta^2 P^\top P
    \end{bmatrix}
\end{equation}
which cannot be written only in terms of the symmetric and antisymmetric parts of $M$.
Interestingly, in order to heuristically obtain an expression for a ``Jacobian'' that is symmetric in the $x$ and $y$ players, a natural idea is to simply consider the average
\begin{equation}
    \frac{\nabla T_{xy}(z^*) + \nabla T_{yx}(z^*)}{2}
    = \begin{bmatrix}
        I - \eta Q - \frac{\eta^2}{2} P P^\top & -\eta P + \frac{\eta^2}{2} P R \\
        \eta P^\top - \frac{\eta^2}{2} P^\top Q & I - \eta R - \frac{\eta^2}{2} P^\top P
    \end{bmatrix}
    = I - \eta M + \frac{\eta^2}{2} A M,
\end{equation}
which yields the same matrix as in \autoref{prop:sym_altGDA:altGDA_upd_withexplicitThalf} (ignoring the $O(\eta^3)$ term).

\begin{proof}
    Denote $L_{xx} = \norm{\nabla_{xx}^2 f}_\infty$, $L_{yy} = \norm{\nabla_{yy}^2 f}_\infty$, $L_{xy} = \norm{\nabla_{xy}^2 f}_\infty$, $L_2 = L_{xx} \vee L_{yy} \vee L_{xy}$ and $L_3 = \norm{\nabla^3 f}_\infty$.
    
    Let $k \in \NN$. 
    We have by definition 
    \begin{align}
        \begin{cases}
            x^{k+\half} = \frac{x^k - \eta \nabla_x f(x^k, y^{k+\half}) + x^k}{2} \\
            y^{k+\half} = y^{k-\half} + \eta \nabla_y f(x^k, y^{k-\half}) \\
            y^{k-\half} = y^k - (y^{k+\half} - y^k)
        \end{cases}
        & \text{and so}~~
        \begin{cases}
            x^{k+\half} = x^k - \frac{1}{2} \eta \nabla_x f(x^k, y^{k+\half}) \\
            y^{k+\half} = 2y^k - y^{k+\half} + \eta \nabla_y f(x^k, 2y^k - y^{k+\half})
        \end{cases} \\
        & \text{i.e.,}~~
        \begin{cases}
            x^{k+\half} = x^k - \frac{1}{2} \eta \nabla_x f(x^k, y^{k+\half}) \\
            y^{k+\half} = y^k + \frac{1}{2} \eta \nabla_y f(x^k, 2y^k - y^{k+\half}).
        \end{cases}
    \end{align}
    It remains to check that the second equation determines $y^{k+\half}$ completely for any given $x^k$ and $y^k$. And indeed, thanks to our assumption that $\frac{\eta}{2} L_{yy} < 1$, the mapping $y \mapsto y^k + \frac{1}{2} \eta \nabla_y f(x^k, 2y^k - y)$ is a contraction, so in particular it has a unique fixed point.
    This shows that $z^{k+\half}$ is entirely determined by $z^k$ for $k \in \NN$, and the case $k \in \NN+\half$ is similar. 
    
    Moreover the above also shows that $T^\half_{xy}$ is characterized by
    \begin{equation}
        (x^{k+\half}, y^{k+\half}) = T^\half_{xy}(x^k, y^k) 
        \iff
        \begin{cases}
            x^{k+\half} &= x^k - \frac{1}{2} \eta \nabla_x f(x^k, y^{k+\half}) \\
            y^{k+\half} &= y^k + \frac{1}{2} \eta \nabla_y f(x^k, 2y^k - y^{k+\half})
        \end{cases}
    \end{equation}
    and symmetrically for $T^\half_{yx}$.
    So by substituting, we indeed find that $T^\half_{xy}(x^*, y^*) = T^\half_{yx}(x^*, y^*) = (x^*, y^*)$.
    In order to compute the Jacobian, let any small $\delta = (\delta_x, \delta_y)$, 
    and consider
    $T^\half_{xy}(z^* + \delta) - T^\half_{xy}(z^*) \eqqcolon \Delta = (\Delta_x, \Delta_y)$.
    Then, denoting $y'$ the unique solution of $y' = y^* + \delta_y + \frac{1}{2} \eta \nabla_y f(x^* + \delta_x, 2y^* + 2 \delta_y - y')$, we have
    \begin{align*}
        y' - y^* &= \delta_y + \frac{\eta}{2} \left[
            \nabla_y f(x^* + \delta_x, 2y^* + 2 \delta_y - y')
            - \nabla_y f(z^*)
        \right] \\
        \norm{y'-y^*} &\leq \norm{\delta_y} + \frac{\eta}{2} \left[ L_{xy} \norm{\delta_x} + L_{yy} \left( 2 \norm{\delta_y} + \norm{y^* - y'} \right) \right] \\
        \norm{y'-y^*} &\leq \frac{1}{1 - \eta L_{yy}/2} \cdot \left(
            \norm{\delta_y} +
            \frac{\eta}{2} \left[ L_{xy} \norm{\delta_x} + L_{yy} 2 \norm{\delta_y} \right]
        \right)
        = O\left( \norm{\delta_x} + \norm{\delta_y} \right)
    \end{align*}
    since we assume that $\eta L_{yy}/2 \leq 1/2$,
    and so
    \begin{align*}
        y' - y^*
        &= \delta_y + \frac{\eta}{2} \left[
            \nabla_{yx}^2 f(z^*) \delta_x
            + \nabla_{yy}^2 f(z^*) (2\delta_y + y^* - y')
            + O\left( \norm{\delta_x}^2 + \norm{\delta_y}^2 \right)
        \right] \\
        y' - y^* &= \left[ I + \frac{\eta}{2} \nabla_{yy}^2 f(z^*) \right]^{-1} \left(
            \delta_y +
            \frac{\eta}{2} \left[
                \nabla_{yx}^2 f(z^*) \delta_x
                + 2 \nabla_{yy}^2 f(z^*) \delta_y
            \right]
        \right)
        + O\left( \norm{\delta_x}^2 + \norm{\delta_y}^2 \right).
    \end{align*}
    This directly gives an expansion for $\Delta_y = y' - y^*$, and we have
    \begin{align*}
        \Delta_x &= \delta_x - \frac{1}{2} \eta \left[ \nabla_x f(x^*+\delta_x, y') - \nabla_x f(z^*) \right] \\
        &= \delta_x - \frac{1}{2} \eta \left[ 
            \nabla_{xx}^2 f(z^*) \delta_x 
            + \nabla_{xy}^2 f(z^*) (y'-y^*)
        \right]
        + O(\norm{\delta_x}^2 + \norm{\delta_y}^2) \\
        &= \delta_x - \frac{\eta}{2} \nabla_{xx}^2 f(z^*) \delta_x
        - \frac{\eta}{2} \nabla_{xy}^2 f(z^*) 
        \left[ I + \frac{\eta}{2} \nabla_{yy}^2 f(z^*) \right]^{-1} \left(
            \delta_y +
            \frac{\eta}{2} \left[
                \nabla_{yx}^2 f(z^*) \delta_x
                + 2 \nabla_{yy}^2 f(z^*) \delta_y
            \right]
        \right) \\
        &\qquad\qquad
        + O(\norm{\delta_x}^2 + \norm{\delta_y}^2) \\
        &= \left\{
            I - \frac{\eta}{2} \nabla_{xx}^2 f(z^*)
            - \frac{\eta^2}{4} \nabla_{xy}^2 f(z^*) \left[ I + \frac{\eta}{2} \nabla_{yy}^2 f(z^*) \right]^{-1} \nabla_{yx}^2 f(z^*)
        \right\} \delta_x \\
        &\qquad\qquad
        + \left\{
            -\frac{\eta}{2} \nabla_{xy}^2 f(z^*) \left[ I + \frac{\eta}{2} \nabla_{yy}^2 f(z^*) \right]^{-1}
            \left( I + \eta \nabla_{yy}^2 f(z^*) \right)
        \right\} \delta_y 
        + O(\norm{\delta_x}^2 + \norm{\delta_y}^2).
    \end{align*}
    Writing for concision
    $
        \begin{bmatrix}
            \nabla_{xx}^2 f & \nabla_{xy}^2 f \\ 
            \nabla_{yx}^2 f & \nabla_{yy}^2 f
        \end{bmatrix}(z^*)
        = \begin{bmatrix}
            Q & P \\
            P^\top & -R
        \end{bmatrix}
    $
    as usual in this paper,
    the above expansions write
    \begin{align*}
        \Delta_x &= \left\{
            I - \frac{\eta}{2} Q
            - \frac{\eta^2}{4} P \left[ I - \frac{\eta}{2} R \right]^{-1} P^\top
        \right\} \delta_x
        + \left\{
            -\frac{\eta}{2} P \left[ I - \frac{\eta}{2} R \right]^{-1} (I - \eta R)
        \right\} \delta_y 
        + O(\norm{\delta_x}^2 + \norm{\delta_y}^2) \\
        \Delta_y &= \left\{
            \left[ I - \frac{\eta}{2} R \right]^{-1} \frac{\eta}{2} P^\top
        \right\} \delta_x
        + \left\{ 
            \left[ I - \frac{\eta}{2} R \right]^{-1} (I - \eta R)
        \right\} \delta_y
        + O\left( \norm{\delta_x}^2 + \norm{\delta_y}^2 \right)
    \end{align*}
    which implies by definition of the Jacobian that
    \begin{equation}
        \nabla T^{\half}_{xy}(z^*) = \begin{bmatrix}
            I - \frac{\eta}{2} Q - \frac{\eta^2}{4} P \left[ I - \frac{\eta}{2} R \right]^{-1} P^\top
            & -\frac{\eta}{2} P \left[ I - \frac{\eta}{2} R \right]^{-1} (I - \eta R) \\
            \left[ I - \frac{\eta}{2} R \right]^{-1} \frac{\eta}{2} P^\top
            & \left[ I - \frac{\eta}{2} R \right]^{-1} (I - \eta R)
        \end{bmatrix}.
    \end{equation}
    By symmetry, the Jacobian for $T^{\half}_{yx}$ can be obtained from the above expression by swapping $\eta$ for $-\eta$, $Q$ for $-R$, $P^\top$ for $P$, and swapping the lines resp.\ columns of the block matrix, yielding
    \begin{equation}
        \nabla T^{\half}_{yx}(z^*) = \begin{bmatrix}
            \left[ I - \frac{\eta}{2} Q \right]^{-1} (I - \eta Q)
            & -\left[ I - \frac{\eta}{2} Q \right]^{-1} \frac{\eta}{2} P \\
            \frac{\eta}{2} P^\top \left[ I - \frac{\eta}{2} Q \right]^{-1} (I - \eta Q) 
            & I - \frac{\eta}{2} R - \frac{\eta^2}{4} P^\top \left[ I - \frac{\eta}{2} Q \right]^{-1} P
        \end{bmatrix}.
    \end{equation}

    By straightforward but tedious computations,
    and noting that $\left[ I - \frac{\eta}{2} R \right]^{-1}$ commutes with $R$,
    we arrive at the following 
    expression for
    $\nabla \left[ T^\half_{xy} \circ T^\half_{yx} \right](z^*)
    = \nabla T^\half_{xy}(z^*) \cdot \nabla T^\half_{yx}(z^*)$
    (since $z^*$ is a fixed point of both operators):
    \begin{align*}
        & \nabla \left[ T^\half_{xy} \circ T^\half_{yx} \right](z^*) \\
        &= \begin{bmatrix}
            I - \eta Q - \frac{\eta^2}{2} P P^\top \left[ I - \frac{\eta}{2} Q \right]^{-1} (I-\eta Q)
            & -\eta P + \frac{\eta^2}{2} P R + \frac{\eta^3}{4} P P^\top \left[ I - \frac{\eta}{2} Q \right]^{-1} P \\
            \eta P^\top \left[ I - \frac{\eta}{2} Q \right]^{-1} (I-\eta Q)
            & I - \eta R - \frac{\eta^2}{2} P^\top \left[ I - \frac{\eta}{2} Q \right]^{-1} P
        \end{bmatrix} \\
        &= \begin{bmatrix}
            I - \eta Q - \frac{\eta^2}{2} P P^\top
            & -\eta P + \frac{\eta^2}{2} P R \\
            \eta P^\top - \frac{\eta^2}{2} P^\top Q
            & I - \eta R - \frac{\eta^2}{2} P^\top P
        \end{bmatrix}
        + \bE 
        = I - \eta M + \frac{\eta^2}{2} A M + \bE
    \end{align*}
    where the error term $\bE$ is
    \begin{equation}
        \bE = \begin{bmatrix}
            -\frac{\eta^2}{2} P P^\top \left\{
                \left[ I - \frac{\eta}{2} Q \right]^{-1} (I-\eta Q)
                - I
            \right\}
            & \frac{\eta^3}{4} P P^\top \left[ I - \frac{\eta}{2} Q \right]^{-1} P \\
            \eta P^\top \left\{
                \left[ I - \frac{\eta}{2} Q \right]^{-1} (I-\eta Q)
                - \left( I - \frac{\eta}{2} Q \right)
            \right\} 
            & -\frac{\eta^2}{2} P^\top \left\{ \left[ I - \frac{\eta}{2} Q \right]^{-1} - I \right\} P
        \end{bmatrix}
        = O(\eta^3)
    \end{equation}
    thanks to our assumption that $\frac{\eta}{2} L_{xx} < 1$.
    More precisely, by bounding the operator norm of each block and summing the bounds, and using that $\spnorm{\left[I - \frac{\eta}{2} Q \right]^{-1}} \leq \sum_{k=0}^\infty \spnorm{\frac{\eta}{2} Q}^k = \frac{1}{1-\eta \spnorm{Q}/2} \leq \frac{1}{1-\eta L_{xx}/2} \leq 2$, one can check that
    \begin{align*}
        \spnorm{R}
        &\leq
        \spnorm{\left[I - \frac{\eta}{2} Q \right]^{-1}} \left(
            \frac{\eta^2}{2} \spnorm{P}^2 \cdot \frac{\eta}{2} \spnorm{Q}
            + \frac{\eta^3}{4} \spnorm{P}^3
            + \eta \spnorm{P} \cdot \frac{\eta^2}{4} \spnorm{Q}^2
            + \frac{\eta^2}{2} \spnorm{P}^2 \cdot \frac{\eta}{2} \spnorm{Q}
        \right) \\
        &\leq 
        \frac{1}{1-\eta L_{xx}/2} \cdot \eta^3 \spnorm{P} (\spnorm{P} \vee \spnorm{Q})^2
        \leq 2 \eta^3 \spnorm{P} (\spnorm{P} \vee \spnorm{Q})^2.
    \end{align*}
    Note that by definition of $A = \begin{bmatrix} 0 & P \\ -P^\top & 0 \end{bmatrix}$ and $S = \begin{bmatrix} Q & 0 \\ 0 & R \end{bmatrix}$, we have $\spnorm{A} = \spnorm{P}$ and $\spnorm{S} = \spnorm{Q} \vee \spnorm{R}$.
    
    This proves the claimed expansion for $\nabla \left[ T^\half_{xy} \circ T^\half_{yx} \right](z^*)$, and the expansion for $\nabla \left[ T^\half_{yx} \circ T^\half_{xy} \right](z^*)$ follows by symmetry.
\end{proof}



\section{Details for \autoref{subsec:quantitative:DT_algos}}
\label{apx:DT_expansions}

In this section we prove the expansions of the convergence rates of discrete-time algorithms (with non-asymptotic bounds in $\eta$ and $\alpha$) reported in \autoref{tab:expans}.
That is, we derive approximate expressions for $\rho(\nabla T(z^*))$ the spectral radius of the update operator's Jacobian at optimum, for various update rules (GDA, EG, etc.), when the skewed gradient field's Jacobian $M$ has a small symmetric part.

We will repeatedly make use of the following estimate for the spectrum of a perturbed normal matrix.
The expansion itself is a special case of \cite[Eqs.~(5) and (7)]{tao_when_2008} for $M_0$ normal, explicitly pointed out in that reference.
Deriving the explicit bound on the error term involves rather tedious calculations however, so the proof is deferred to \autoref{subsec:DT_expan:proof_gen_expan_formula}.

\begin{restatable}{proposition}{DTexpanGenExpanFormula}
\label{prop:DT_expan:gen_expan_formula}
    Let $M_0, M_1, M_2 \in \RR^{d \times d}$ or $\CC^{d \times d}$ and
    $M_\alpha = M_0 + \alpha M_1 + \frac{\alpha^2}{2} M_2$.
    Assume $M_0$ is normal and has distinct eigenvalues; denote its eigenvalue decomposition as $M_0 w_j = \lambda_j w_j$ with $(w_j)_j$ unitary,
    and let $\gamma_0 = \min_{k \neq j} \abs{\lambda_k - \lambda_j}$.
    Then for all $\alpha$ 
    such that
    $\spnorm{\alpha M_1 + \frac{\alpha^2}{2} M_2} \leq \frac{\gamma_0}{4 \sqrt{2d}}$,
    \begin{equation}
        \spectrum(M_\alpha) = \left\lbrace
            \lambda_j + \alpha \olw_j^\top M_1 w_j
            + \frac{\alpha^2}{2} \olw_j^\top M_2 w_j
            + \alpha^2 \sum_{k \neq j} \frac{ \left( \olw_j^\top M_1 w_k \right) \left( \olw_k^\top M_1 w_j\right) }{\lambda_j-\lambda_k}
            + \br_j,~
            1 \leq j \leq d
        \right\rbrace
    \end{equation}
    where $\spnorm{\cdot}$ denotes operator norm and for each $j$,
    ${
        \abs{\br_j} \leq 
        \alpha^3 \cdot 8 d \gamma_0^{-1} \spnorm{M_1} \left(
            \spnorm{M_2} + 4 d \gamma_0^{-1} \spnorm{M_1}^2
        \right)
    }$.
\end{restatable}

In the remainder of this section, consider $S$ real symmetric and $A$ real antisymmetric in $\RR^{d \times d}$, and let $M = M_\alpha = \alpha S + A$.
Assume $A$ has simple eigenvalues and denote its eigenvalue decomposition as $A w_j = i \sigma_j w_j$ with $\sigma_j \in \RR$.
This corresponds exactly to \autoref{sec:quantitative} with the correspondence $d \coloneqq n + n$ and 
$\spectrum(A) = \left\{ i \sigma_j, j \leq d \right\} \coloneqq \left\{ i s \tsigma_j, s \in \{-1, 1\}, j \leq n \right\}$ with $(\tsigma_j)_{j \leq n}$ the singular values of $P$.
Furthermore, let for concision
$\gamma_A = \min_{k \neq j} \abs{\sigma_k - \sigma_j}$
and
$L_S = \spnorm{S}$,
$L_A = \spnorm{A}$.

Throughout the derivations of this section, we will write $\beps, \beps', \beps_j$ or $\beps'_j$ to denote elements of $[-1, 1]$ or of the unit ball of $\RR^d$, and that may change from line to line.
Likewise, $\bzeta, \bzeta', \bzeta_j, \bzeta'_j$ will denote elements of $\{ z \in \CC; \abs{z} \leq 1 \}$ or of the unit ball of $\CC^d$, and that may change from line to line.

\subsection{Simultaneous GDA} \label{subsec:DT_expan:simGDA}

The Sim-GDA update rule is
$z^{k+1} = T(z^k) = z^k - \eta g(z^k)$,
so the update operator's Jacobian at optimum is
\begin{equation}
    \nabla T(z^*) = I - \eta M = \left( I - \eta A \right) - \alpha \eta S.
\end{equation}
Observe that $I-\eta A$ is normal with eigenvalue/eigenvector pairs $(1 - i\eta \sigma_j, w_j)$.
Hence let us apply \autoref{prop:DT_expan:gen_expan_formula} with $M_0 = I-\eta A$ and $M_1 = -\eta S$ (and $M_2=0$).
We get that
\begin{multline*}
    \forall \alpha \leq \frac{1}{\eta L_S} \frac{\eta \gamma_A}{4 \sqrt{2d}} = \frac{\gamma_A}{L_S \cdot 4 \sqrt{2d}},~~ \\
    \spectrum(\nabla T(z^*)) = \left\lbrace
        1 - \eta i \sigma_j 
        - \alpha \eta \olw_j^\top S w_j
        + \alpha^2 \sum_{k \neq j} \frac{ \eta^2 \abs{\olw_j^\top S w_k}^2 }{-i \eta \sigma_j + i \eta \sigma_k}
        + \alpha^3 \eta \cdot \bzeta_j \left( 32 d^2 \gamma_A^{-2} L_S^3 \right),
        1 \leq j \leq d
    \right\rbrace.
\end{multline*}
It seems unlikely that the term in $\alpha^2$ will ever lead to friendly expressions, 
so we will bound it uniformly; let us nonetheless note that it is pure imaginary.
Namely, since
\begin{equation} \label{eq:DT_expansions:bound_alpha2term_nod}
    \abs{ \alpha^2 \sum_{k \neq j} \frac{ \eta^2 \abs{\olw_j^\top S w_k}^2 }{-i \eta \sigma_j + i \eta \sigma_k} }
    \leq \alpha^2 \eta \sum_{k \neq j} \frac{\abs{\olw_k^\top S w_j}^2}{\abs{\sigma_k-\sigma_j}}
    \leq \alpha^2 \eta \gamma_A^{-1} \underbrace{ \sum_k \abs{\olw_k^\top S w_j}^2 }_{= \norm{S w_j}^2 \leq L_S^2},
\end{equation}
then we may write that term as
$i \alpha^2 \eta \cdot \beps_j \left( \gamma_A^{-1} L_s^2 \right)$.
For the spectral radius we get
\begin{align*}
    & \rho(\nabla T(z^*))^2 
    = \max_{\lambda \in \spectrum(\nabla T(z^*))} \abs{\lambda}^2 \\
    &= \max_{j \leq d} \abs{1 - \alpha \eta \olw_j^\top S w_j + \alpha^3 \eta \cdot \beps'_j 32 d^2 \gamma_A^{-2} L_S^3}^2
    + \abs{\eta \sigma_j + \alpha^2 \eta \cdot \beps_j \gamma_A^{-1} L_S^2 + \alpha^3 \eta \cdot \beps''_j 32 d^2 \gamma_A^{-2} L_S^3}^2 \\
    &= \max_{j \leq d} 1 - 2 \alpha \eta \left( \olw_j^\top S w_j \right) + \eta^2 \sigma_j^2 + O(\alpha^3 \eta + \alpha^2 \eta^2)
\end{align*}
and more precisely one can check that the $O(\cdot)$ term is absolutely bounded by
$\alpha^3 \eta \cdot 128 d^2 \gamma_A^{-2} L_S^3 \left( 1 + 5 \eta L_A \right) + \alpha^2 \eta^2 \cdot 2 \gamma_A^{-1} L_S^2 L_A$,
for all $\alpha \leq \frac{\gamma_A}{L_S \cdot 4 \sqrt{2d}}$.

\subsection{Proximal Point} \label{subsec:DT_expan:PP}

The PP update rule is $z^{k+1} = T(z^k)$ with $T^{-1}(z) = z + \eta g(z)$ so the update operator's Jacobian at optimum is
$\nabla T(z^*) = \left( I + \eta M \right)^{-1}$.
So
$\spectrum(\nabla T(z^*)) = \left\lbrace \lambda^{-1}, \lambda \in \spectrum(I+\eta M) \right\rbrace$ and so, by the exact same calculations as for Sim-GDA with $\eta$ replaced by $-\eta$,
\begin{equation}
    \rho(\nabla T(z^*))^2 = \left[
        \min_{j \leq d}
        1 + 2 \alpha \eta \left( \olw_j^\top S w_j \right) + \eta^2 \sigma_j^2 
        + O(\alpha^3 \eta + \alpha^2 \eta^2)
    \right]^{-1}
\label{eq:DT_expansions:PP}
\end{equation}
for all $\alpha \leq \frac{\gamma_A}{L_S \cdot 4\sqrt{2d}}$,
and we have the same absolute bound on the $O(\cdot)$ term as for Sim-GDA.

\subsection{Alternating GDA} \label{subsec:DT_expan:altGDA}

As we show in \autoref{apx:sym_altGDA},
for any $\eta \leq \norm{\nabla_{xx}^2 f}_\infty^{-1} \wedge \norm{\nabla_{yy}^2 f}_\infty^{-1}$,
the (symmetrized) Alt-GDA update operator's Jacobian at optimum is
\begin{gather*}
    \nabla T(z^*) 
    = \nabla \olT(z^*) + \bE 
    ~~~~\text{for some}~~ \bE \in \RR^{d \times d} ~~\text{with}~~ \spnorm{\bE} \leq 2 \eta^3 L_A (L_A \vee \alpha L_S)^2, \\
    \text{where}~~
    \nabla \olT(z^*) = I - \eta M + \frac{\eta^2}{2} A M
    = \left( I - \eta A + \frac{\eta^2}{2} A^2 \right)
    + \alpha \left( -\eta S + \frac{\eta^2}{2} A S \right).
\end{gather*}
Observe that $I - \eta A + \frac{\eta^2}{2} A^2$ is normal with eigenvalue/eigenvector pairs $\left( 1 - i\eta \sigma_j - \frac{\eta^2}{2} \sigma_j^2, w_j \right)$.
Hence let us apply \autoref{prop:DT_expan:gen_expan_formula} with $M_0 = I-\eta A + \frac{\eta^2}{2} A^2$ and $M_1 = -\eta S + \frac{\eta^2}{2} A S$ (and $M_2=0$).
Further observe that
\begin{equation}
    \olw_j^\top A S w_k
    = \olw_j^\top \left( \sum_l i \sigma_l w_l \olw_l^\top \right) S w_k
    = i \sigma_j \cdot \olw_j^\top S w_k.
\end{equation}
So we get that, for all $\alpha \leq \frac{\gamma_A}{4 \sqrt{2d} \cdot L_s \left( 1+\frac{\eta}{2} L_A \right)}$,
\begin{align*}
    &\spectrum(\nabla \olT(z^*)) \\
    &= \Bigg\lbrace
        1 - \eta i\sigma_j - \frac{\eta^2}{2} \sigma_j^2
        + \alpha \olw_j^\top \left( -\eta S + \frac{\eta^2}{2} AS \right) w_j
        + \alpha^2 \sum_{k \neq j} \frac{
            \eta^2
            \left( \olw_j^\top \left( S - \frac{\eta}{2} A S \right) w_k \right)
            \left( \olw_k^\top \left( S - \frac{\eta}{2} A S \right) w_j \right)
        }{
            -i\eta \sigma_j + i \eta \sigma_k + \frac{\eta^2}{2} (-\sigma_j^2 + \sigma_k^2)
        } \\
        &\qquad\qquad
        + \alpha^3 \eta \cdot \bzeta_j\left( 32 d^2 \gamma_A^{-2} L_S^3 \left(1+\frac{\eta}{2} L_A \right)^3 \right),
        1 \leq j \leq d
    \Bigg\rbrace \\
    &= \Bigg\lbrace
        1 - \eta i\sigma_j - \frac{\eta^2}{2} \sigma_j^2
        - \alpha \eta \left( \olw_j^\top S w_j \right)
        + \alpha \frac{\eta^2}{2} i \sigma_j \cdot \left( \olw_j^\top S w_j \right) 
        + i \cdot \alpha^2 \eta \cdot \beps_j\left( \gamma_A^{-1} L_S^2 \right)
        \\
        &\qquad\qquad
        + \alpha^3 \eta \cdot \bzeta_j\left( 32 d^2 \gamma_A^{-2} L_S^3 \left(1+\frac{\eta}{2} L_A \right)^3 \right)
        + \alpha^2 \eta^3 \cdot \bzeta'_j \left( \gamma_A^{-1} L_S^2 L_A^2 /4 \right),
        1 \leq j \leq d
    \Bigg\rbrace.
\end{align*}
Here for the second equality, we computed the coefficient for the term in $\alpha^2$ as
\begin{align*}
    &
    \sum_{k \neq j} \frac{
        \eta^2
        \left( \olw_j^\top \left( S - \frac{\eta}{2} A S \right) w_k \right)
        \left( \olw_k^\top \left( S - \frac{\eta}{2} A S \right) w_j \right)
    }{
        -i\eta \sigma_j + i \eta \sigma_k + \frac{\eta^2}{2} (-\sigma_j^2 + \sigma_k^2)
    } \\
    &= \sum_{k \neq j} \eta \frac{
        \left( \olw_j^\top S w_k - \frac{\eta}{2} i \sigma_j \olw_j^\top S w_k \right)
        \left( \olw_k^\top S w_j - \frac{\eta}{2} i \sigma_k \olw_k^\top S w_j \right)
    }{
        i (\sigma_k - \sigma_j) 
        + \frac{\eta}{2} (\sigma_k^2 - \sigma_j^2)
    }
    = \sum_{k \neq j} \eta \abs{\olw_j^\top S w_k}^2 
    \frac{
        \left( 1 - \frac{\eta}{2} i \sigma_j \right)
        \left( 1 - \frac{\eta}{2} i \sigma_k \right)
    }{
        i (\sigma_k - \sigma_j)
        \left( 1 - i \frac{\eta}{2} (\sigma_k + \sigma_j) \right)
    } \\
    &= \sum_{k \neq j}
    \frac{ \eta \abs{\olw_j^\top S w_k}^2 }{i (\sigma_k-\sigma_j)}
    \left( 1 + \frac{ -\frac{\eta^2}{4} \sigma_j \sigma_k }{1-i\frac{\eta}{2} (\sigma_k+\sigma_j)} \right) 
    = i \eta \cdot \beps_j( \gamma_A^{-1} L_S^2 )
    + \eta^3 \cdot \bzeta'_j\left( \gamma_A^{-1} L_S^2 \cdot L_A^2/4 \right)
\end{align*}
where the last equality follows from the same bound as in~\eqref{eq:DT_expansions:bound_alpha2term_nod}.
For the spectral radius we get
\begin{align*}
    \rho(\nabla \olT(z^*))^2 &= \max_{j \leq d} \abs{1 - \frac{\eta^2}{2} \sigma_j^2 - \alpha \eta \left( \olw_j^\top S w_j \right)}^2
    + \abs{-\eta \sigma_j + \alpha \frac{\eta^2}{2} \sigma_j \cdot \left( \olw_j^\top S w_j \right) + O(\alpha^2 \eta)}^2 \\
    &\qquad\qquad + O(\alpha^3 \eta + \alpha^2 \eta^3) \\
    &= \max_{j \leq d} 1 + \frac{\eta^4}{4} \sigma_j^4 - \eta^2 \sigma_j^2 - 2 \alpha \eta \left( \olw_j^\top S w_j \right) + \alpha \eta^3 \sigma_j^2 \left( \olw_j^\top S w_j \right) \\
    &\qquad\qquad + \eta^2 \sigma_j^2 - \alpha \eta^3 \sigma_j^2 \left( \olw_j^\top S w_j \right) 
    + O(\alpha^3 \eta + \alpha^2 \eta^2) \\
    &= \max_{j \leq d} 1 - 2 \alpha \eta \left( \olw_j^\top S w_j \right) + \frac{\eta^4}{4} \sigma_j^4 
    + O(\alpha^3 \eta + \alpha^2 \eta^2).
\end{align*}
The two cancellations that occur in the last line (of $\pm \eta^2 \sigma_j^2$ and of $\pm \alpha \eta^3 \sigma_j^2 \left( \olw_j^\top S w_j \right)$) are consistent with the intuition that Alt-GDA is a good symplectic integrator.
More precisely, one can check that the final $O(\cdot)$ term is absolutely bounded by
$
    \alpha^3 \eta \cdot 512 d^2 \gamma_A^{-2} L_S^3 (1+\eta L_A)^6
    + \alpha^2 \eta^2 \cdot 4 d \gamma_A^{-1} L_S^2 L_A (1+\eta L_A)^4
$,
for all $\alpha \leq \frac{\gamma_A}{4 \sqrt{2d} \cdot L_s \left( 1+\frac{\eta}{2} L_A \right)}$.

Finally, by smoothness of the eigenvalues of perturbed matrices (\autoref{lm:DT_expan:restate_taoeigenvalues2008}),
\begin{equation}
    \rho(\nabla T(z^*))^2 = \rho(\nabla \olT(z^*) + \bE)^2
    = \rho(\nabla \olT(z^*)^2 + O(\eta^3).
\end{equation}
One could further derive explicit bounds on the $O(\eta^3)$ term in that last expression
by controlling the distance to normality of $\nabla \olT(z^*)$ -- since it is continuous in $\alpha$ and normal for $\alpha=0$ -- and adapting the proof of \autoref{lm:DT_expan:control_chi} to matrices that are close to normal.

\subsection{Extra-Gradient} \label{subsec:DT_expan:EG}

The EG update rule is $z^{k+1} = T(z^k) = z^k - \eta g(z^k - \eta g(z^k))$
so the update operator's Jacobian at optimum is
\begin{align*}
    \nabla T(z^*) = I - \eta M (I - \eta M)
    &= I - \eta A - \alpha \eta S + \eta^2 M^2 \\
    &= \left( I - \eta A + \eta^2 A^2 \right) + \alpha \left( -\eta S + \eta^2 (AS + SA) \right) + \frac{1}{2} \alpha^2 \left( 2 \eta^2 S^2 \right).
\end{align*}
Similarly as in the previous subsection, observe that $I - \eta A + \eta^2 A^2$ is normal with eigenvalue/eigenvector pairs $\left( 1 - i\eta \sigma_j - \eta^2 \sigma_j^2, w_j \right)$, and that
\begin{equation}
    \olw_j^\top (A S + S A) w_k
    = i (\sigma_j + \sigma_k) \cdot \olw_j^\top S w_k.
\end{equation}
Let us apply \autoref{prop:DT_expan:gen_expan_formula} to $M_0 = I - \eta A + \eta^2 A^2$, $M_1 = -\eta S + \eta^2 (AS + SA)$ and $M_2 = 2\eta^2 S^2$.
The proposition yields a bound for all 
$\alpha$ such that $\alpha \spnorm{S - \eta (AS + SA) + \alpha \eta S^2} \leq \frac{\gamma_A}{4 \sqrt{2d}}$,
for which a simpler sufficient condition is
$\alpha \leq 1 \wedge \frac{\gamma_A}{4 \sqrt{2d} \cdot L_S (1+\eta L_S (1+2 L_A))}$.
Using that $\spnorm{M_1} \leq \eta L_S (1+2\eta L_A)$
to upper-bound the term in $\alpha^3$,
we get that for all such $\alpha$,
\begin{align*}
    \spectrum(\nabla T(z^*)) 
    &= \Bigg\lbrace
        1 - i \eta \sigma_j - \eta^2 \sigma_j^2
        + \alpha \left(
            - \eta \olw_j^\top S w_j + \eta^2 2 i \sigma_j \olw_j^\top S w_j
        \right)
        + \alpha^2 \eta^2 \olw_j^\top S^2 w_j \\
        &\qquad\qquad
        + \alpha^2 \sum_{k \neq j} \frac{
            \eta^2
            \olw_j^\top (S - \eta (AS + SA)) w_k
            \cdot 
            \olw_k^\top (S - \eta (AS + SA)) w_j
        }{
            -i\eta \sigma_j + i \eta \sigma_k + \eta^2 (-\sigma_j^2 + \sigma_k^2)
        } \\
        &\qquad\qquad
        + \alpha^3 \eta \cdot \bzeta_j \cdot 8d \gamma_A^{-1} L_S (1+2\eta L_A) \left( 2\eta L_S^2 + 4 d \gamma_A^{-1} L_S^2 (1+2\eta L_A)^2 \right),
        1 \leq j \leq d
    \Bigg\rbrace \\
    &= \Bigg\lbrace
        1 - i \eta \sigma_j - \eta^2 \sigma_j^2
        - \alpha \eta \left( \olw_j^\top S w_j \right)
        + 2 i \alpha \eta^2 \sigma_j \left( \olw_j^\top S w_j \right) \\
        &\qquad\qquad
        + i \cdot \alpha^2 \eta \cdot \beps_j\left( \gamma_A^{-1} L_S^2 \right)
        + \alpha^2 \eta^2 \cdot \beps'_j\left( 2 \gamma_A^{-1} L_S^2 L_A \right)
        + \alpha^2 \eta^2 \olw_j^\top S^2 w_j \\
        &\qquad\qquad
        + \alpha^3 \eta \cdot \bzeta_j \left( 512 d^2 \gamma_A^{-2} L_S^3 (1+\eta L_A)^3 \right),
        1 \leq j \leq d
    \Bigg\rbrace.
\end{align*}
Here for the second equality we computed the terms in $\alpha^2$ as
\begin{align*}
    & \alpha^2 \sum_{k \neq j} \frac{
        \eta^2
        \olw_j^\top (S - \eta (AS + SA)) w_k
        \cdot 
        \olw_k^\top (S - \eta (AS + SA)) w_j
    }{
        -i\eta \sigma_j + i \eta \sigma_k + \eta^2 (-\sigma_j^2 + \sigma_k^2)
    } \\
    &= \alpha^2 \sum_{k \neq j} \eta \frac{
        \left( \olw_j^\top S w_k \right) \left( 1 - i \eta (\sigma_j + \sigma_k) \right)
        \cdot 
        \left( \olw_k^\top S w_j \right) \left( 1 - i \eta (\sigma_j + \sigma_k) \right)
    }{
        i (\sigma_k-\sigma_j) \left( 1 - i \eta (\sigma_k+\sigma_j) \right)
    } \\
    &= \alpha^2 \sum_{k \neq j} \frac{\eta \abs{\olw_j^\top S w_k}^2}{i (\sigma_k-\sigma_j)} \left( 1 - i \eta (\sigma_k+\sigma_j) \right) 
    = i \cdot \alpha^2 \eta \cdot \beps_j\left( \gamma_A^{-1} L_S^2 \right) + \alpha^2 \eta^2 \cdot \beps'_j\left( \gamma_A^{-1} L_S^2 \cdot 2 L_A \right)
\end{align*}
where the last equality follows from the same bound as in~\eqref{eq:DT_expansions:bound_alpha2term_nod}.
For the spectral radius we get
\begin{align*}
    \rho(\nabla T(z^*))^2 
    &= \max_{j \leq d}
    \abs{1 - \eta^2 \sigma_j^2 - \alpha \eta \left( \olw_j^\top S w_j \right) + O(\alpha^2 \eta^2)}^2
    + \abs{ -\eta \sigma_j + 2\alpha \eta^2 \sigma_j \left( \olw_j^\top S w_j \right) + O(\alpha^2 \eta)}^2 \\
    &\qquad\qquad
    + O(\alpha^3 \eta) \\
    &= \max_{j \leq d}
    1 + \eta^4 \sigma_j^4 - 2 \eta^2 \sigma_j^2 - 2\alpha \eta \left( \olw_j^\top S w_j \right) + 2 \alpha \eta^3 \sigma_j^2 \left( \olw_j^\top S w_j \right) \\
    &\qquad\qquad + \eta^2 \sigma_j^2 - 4 \alpha \eta^3 \sigma_j^2 \left( \olw_j^\top S w_j \right)
    + O(\alpha^3 \eta + \alpha^2 \eta^2) \\
    &= \max_{j \leq d}
    1 - 2 \alpha \eta \left( \olw_j^\top S w_j \right)
    - \eta^2 \sigma_j^2
    - 2 \alpha \eta^3 \sigma_j^2 \left( \olw_j^\top S w_j \right)
    + \eta^4 \sigma_j^4
    + O(\alpha^3 \eta + \alpha^2 \eta^2)
\end{align*}
and more precisely one can check that the final $O(\cdot)$ term is absolutely bounded by 
$\alpha^3 \eta \cdot 2^{12} d^2 \gamma_A^{-2} L_S^3 (1+\eta L_A)^7 + \alpha^2 \eta^2 \cdot 15 d \gamma_A^{-1} L_S^2 L_A (1+\eta L_A)^2$.



\subsection{Proof of \autoref{prop:DT_expan:gen_expan_formula}}
\label{subsec:DT_expan:proof_gen_expan_formula}

The following lemma summarizes the eigenvalue derivative formulas up to order $2$ for perturbed matrices with distinct eigenvalues.

\begin{lemma}[\cite{tao_when_2008}] \label{lm:DT_expan:restate_taoeigenvalues2008}
    Let $M_\alpha \in \CC^{d \times d}$ for $\alpha \in (-1, 1)$ a smooth curve of matrices
    such that $M_0$ has distinct eigenvalues.
    Then there exists an open interval $I \ni 0$ such that $M_\alpha$ has distinct eigenvalues for all $\alpha \in I$; in particular they are diagonalizable.
    Denote their eigenvalue decompositions as
    $M_\alpha = \sum_{k=1}^d \lambda_k(\alpha) v_k(\alpha) \olw_k(\alpha)^\top$; that is, the eigenvalues of $M_\alpha$ are $(\lambda_k(\alpha))_k$, the associated eigenvectors are $(v_k(\alpha))_k$, and $(w_k(\alpha))_k$ is a dual basis to the basis of eigenvectors, i.e., $\olw_k(\alpha)^\top v_j(\alpha) = \ind_{j=k}$.

    The eigenvalues of $M_\alpha$ are smooth over $I$ and their first two derivatives at any $\alpha \in I$ are given by, using $\dot{\bullet}$ to denote differentiation w.r.t.\ $\alpha$
    and leaving the dependency on $\alpha$ implicit,
    \begin{equation}
        \dot{\lambda}_k = w^*_k \dot{M} v_k
        \quad \text{and} \quad
        \ddot{\lambda}_k = w_k^* \ddot{M} v_k + 2 \sum_{j \neq k} \frac{(w^*_k \dot{M} v_j) (w^*_j \dot{M} v_k)}{\lambda_k-\lambda_j}
    \end{equation}
    where we denoted $w^*_k = \olw_k^\top$. 
    Furthermore the eigenvectors $v_k$ and dual basis vectors $w_k$ can also be chosen smooth and their derivatives at any $\alpha \in I$ are given by
    \begin{equation}    
        \dot{v}_k = \sum_{j \neq k} \frac{w^*_j \dot{M}v_k}{\lambda_k - \lambda_j} v_j + c_k v_k
        \quad \text{and} \quad
        \dot{w}^*_k = \sum_{j \neq k} \frac{w^*_k \dot{M} v_j}{\lambda_k - \lambda_j} w^*_j - c_k w^*_k
    \end{equation}
    for some scalars $c_k$ that reflect the normalization of the eigenvectors.
\end{lemma}

By applying a Taylor expansion with remainder in Lagrange form to the eigenvalues $\lambda_k(\alpha)$ of the matrix $M_\alpha = M_0 + \alpha M_1 + \frac{\alpha^2}{2} M_2$ of \autoref{prop:DT_expan:gen_expan_formula} -- since the eigenvalues are smooth by the above lemma --,
we have that for all $\alpha$ in some neighborhood of zero,
\begin{equation} \label{eq:DT_expan:gen_expan_TaylorLagrange}
    \lambda_k(\alpha) = \lambda_k(0) + \alpha \dot{\lambda}_k(0) + \frac{\alpha^2}{2} \ddot{\lambda}_k(0) + \frac{\alpha^3}{6} \dddot{\lambda}_k(\xi)
    ~~~~\text{for some}~~ 0 < \xi < \alpha.
\end{equation}
By substituting the expressions from \autoref{lm:DT_expan:restate_taoeigenvalues2008} for the first two eigenvalue derivatives, we already get the terms in $\alpha$ and $\alpha^2$ in \autoref{prop:DT_expan:gen_expan_formula}.
In order to control the last term in $\alpha^3$, we need to compute the third eigenvalue derivatives $\dddot{\lambda}_k(\xi)$.
Note that $M_\xi$ is never assumed normal for any $\xi >0$, which is why we do not use normality in \autoref{lm:DT_expan:restate_taoeigenvalues2008} nor in \autoref{lm:DT_expan:third_eigval_deriv} below.

\begin{lemma} \label{lm:DT_expan:third_eigval_deriv}
    Under the conditions of \autoref{lm:DT_expan:restate_taoeigenvalues2008},
    \begin{multline*}
        \dddot{\lambda}_k =
        w^*_k \dddot{M} v_k
        + 3 \sum_{j \neq k} \frac{1}{\lambda_k-\lambda_j}
        \left[
            (w^*_k \dot{M} v_j) (w^*_j \ddot{M} v_k) 
            + (w^*_k \ddot{M} v_j) (w^*_j \dot{M} v_k)
        \right] \\
        + 6 \sum_{j, l \neq k} \frac{(w^*_k \dot{M} v_j) (w^*_j \dot{M} v_l) (w^*_l \dot{M} v_k)}{(\lambda_k-\lambda_j) (\lambda_k-\lambda_l)}
        - 6 \sum_{j \neq k} \frac{(w^*_k \dot{M} v_j) (w^*_j \dot{M} v_k)}{(\lambda_k-\lambda_j)^2} w^*_k \dot{M} v_k.
    \end{multline*}
    In particular, 
    \begin{equation}
        \max_k \abs{ \dddot{\lambda}_k } \leq 
        \chi \spnorm{\dddot{M}} 
        + 6 d \gamma^{-1} \chi^2 \spnorm{\dot{M}} \spnorm{\ddot{M}}
        + 6 d^2 \gamma^{-2} \chi^3 \spnorm{\dot{M}}^3
    \end{equation}
    where $\gamma = \min_{j \neq k} \abs{\lambda_j - \lambda_k}$, 
    $\chi = \max_k \norm{v_k} \norm{w_k}$ 
    and $\spnorm{\cdot}$ denotes operator norm.
\end{lemma}

\begin{proof}
    By differentiating the identity $M v_k = \lambda_k v_k$, we have that
    \begin{align*}
        \dot{M} v_k + M \dot{v}_k &= \dot{\lambda}_k v_k + \lambda_k \dot{v}_k \\
        \ddot{M} v_k + 2 \dot{M} \dot{v}_k + M \ddot{v}_k &= \ddot{\lambda}_k v_k + 2 \dot{\lambda}_k \dot{v}_k + \lambda \ddot{v}_k 
        \\
        \dddot{M} v_k + 3 \ddot{M} \dot{v}_k + 3 \dot{M} \ddot{v}_k + M \dddot{v}_k &= \dddot{\lambda}_k v_k + 3 \ddot{\lambda}_k \dot{v}_k + 3 \dot{\lambda}_k \ddot{v}_k + \lambda_k \dddot{v}_k.
    \end{align*}
    Also note that $w^*_k M = \lambda_k w^*_k$.
    By multiplying the third identity by $w^*_k$ on the left, we have that
    \begin{align*}
        w^*_k \dddot{M} v_k + 3 w^*_k \ddot{M} \dot{v}_k + 3 w^*_k \dot{M} \ddot{v}_k + w^*_k M \dddot{v}_k 
        &= \dddot{\lambda}_k + 3 \ddot{\lambda}_k w^*_k \dot{v}_k + 3 \dot{\lambda}_k w^*_k \ddot{v}_k + \lambda_k w^*_k \dddot{v}_k \\
        w^*_k \dddot{M} v_k + 3 w^*_k \ddot{M} \dot{v}_k + 3 w^*_k \dot{M} \ddot{v}_k
        &= \dddot{\lambda}_k + 3 \ddot{\lambda}_k w^*_k \dot{v}_k + 3 \dot{\lambda}_k w^*_k \ddot{v}_k
    \end{align*}
    since $(w^*_k M) \dddot{v}_k = \lambda_k w^*_k \dddot{v}_k$.
    Now let us compute $\ddot{v}_k$. By multiplying the identity for the second derivatives by $w^*_j$ on the left for any $j \neq k$, we have that
    \begin{align*}
        w^*_j \ddot{M} v_k + 2 w^*_j \dot{M} \dot{v}_k + \underbrace{w^*_j M}_{=\lambda_j w^*_j} \ddot{v}_k &= \ddot{\lambda}_k w^*_j v_k + 2 \dot{\lambda}_k w^*_j \dot{v}_k + \lambda_k w^*_j \ddot{v}_k \\
        w^*_j \ddot{v}_k &= \frac{1}{\lambda_j - \lambda_k} \left( 2 \dot{\lambda}_k w^*_j \dot{v}_k - w^*_j \ddot{M} v_k - 2 w^*_j \dot{M} \dot{v}_k \right).
    \end{align*}
    Hence we can compute
    \begin{multline*}
        w^*_k \dot{M} \ddot{v}_k - \dot{\lambda}_k w^*_k \ddot{v}_k
        = w^*_k \dot{M} \sum_j (w^*_j \ddot{v}_k) v_j
        - \dot{\lambda}_k w^*_k \ddot{v}_k \\
        = \sum_{j \neq k} \left( w^*_k \dot{M} v_j \right) \frac{1}{\lambda_j - \lambda_k} \left( 2 \dot{\lambda}_k w^*_j \dot{v}_k - w^*_j \ddot{M} v_k - 2 w^*_j \dot{M} \dot{v}_k \right)
        + \underbrace{
            (w^*_k \dot{M} v_k) (w^*_k \ddot{v}_k)
            - \dot{\lambda}_k w^*_k \ddot{v}_k
        }_{=0}
    \end{multline*}
    since $\dot{\lambda}_k = w^*_k \dot{M} v_k$ by \autoref{lm:DT_expan:restate_taoeigenvalues2008}.
    Substituting back into the identity for the third derivative left-multiplied by $w^*_k$, we find that
    \begin{equation*}
        \dddot{\lambda}_k = w^*_k \dddot{M} v_k + 3 w^*_k \ddot{M} \dot{v}_k - 3 \ddot{\lambda}_k w^*_k \dot{v}_k 
        + 3 \sum_{j \neq k} \frac{w^*_k \dot{M} v_j}{\lambda_j - \lambda_k} \left( 2 \dot{\lambda}_k w^*_j \dot{v}_k - w^*_j \ddot{M} v_k - 2 w^*_j \dot{M} \dot{v}_k \right).
    \end{equation*}
    The claimed expression for $\dddot{\lambda}_k$ will follow by substituting the expressions for $\dot{v}_k$, $\dot{\lambda}_k$ and $\ddot{\lambda}_k$ from \autoref{lm:DT_expan:restate_taoeigenvalues2008} and simplifying.
    Namely, letting for concision $\delta_{jk} = \lambda_j - \lambda_k$ and $M_{jk} = w^*_j M v_k$, $\dot{M}_{jk} = w^*_j \dot{M} v_k$ and similarly for $\ddot{M}_{jk}$, $\dddot{M}_{jk}$ for any $j, k \leq d$, we have
    \begin{equation}
        \dot{\lambda}_k = \dot{M}_{kk}, \qquad
        \dot{v}_k = \sum_{j \neq k} \frac{\dot{M}_{jk}}{\delta_{kj}} v_j + c_k v_k, \qquad
        \ddot{\lambda}_k = \ddot{M}_{kk} + 2 \sum_{j \neq k} \frac{\dot{M}_{kj} \dot{M_{jk}}}{\delta_{kj}}
    \end{equation}
    and in particular $w^*_j \dot{v}_k = \frac{\dot{M}_{jk}}{\delta_{kj}}$ for all $j \neq k$ and so
    \begin{align*}
        \dddot{\lambda}_k &= \dddot{M}_{kk} + 3 \sum_{j \neq k} \frac{\ddot{M}_{kj} \dot{M}_{jk}}{\delta_{kj}} + 3 \ddot{M}_{kk} c_k
        - 3 \left( \ddot{M}_{kk} + 2 \sum_{j \neq k} \frac{\dot{M}_{kj} \dot{M}_{jk}}{\delta_{kj}} \right) c_k \\
        &\qquad - 3 \sum_{j \neq k} \frac{\dot{M}_{kj} \ddot{M}_{jk}}{\delta_{jk}}
        + 6 \sum_{j \neq k} \frac{\dot{M}_{kj}}{\delta_{jk}} \left( 
            \dot{M}_{kk} \frac{\dot{M}_{jk}}{\delta_{kj}} - \sum_{l \neq k} \frac{\dot{M}_{jl} \dot{M}_{lk}}{\delta_{kl}} - \dot{M}_{jk} c_k
        \right) \\
        &= \dddot{M}_{kk} + 3 \sum_{j \neq k} \frac{\ddot{M}_{kj} \dot{M}_{jk} + \dot{M}_{kj} \ddot{M}_{jk}}{\delta_{kj}}
        + 6 \sum_{j \neq k} \frac{\dot{M}_{kj}}{\delta_{jk}} \left( 
            \dot{M}_{kk} \frac{\dot{M}_{jk}}{\delta_{kj}} - \sum_{l \neq k} \frac{\dot{M}_{jl} \dot{M}_{lk}}{\delta_{kl}}
        \right).
    \end{align*}
    In order to simplify the last term, we simply write it as
    \begin{align*}
        \sum_{j \neq k} \frac{\dot{M}_{kj}}{\delta_{jk}} \left( 
            \dot{M}_{kk} \frac{\dot{M}_{jk}}{\delta_{kj}} - \sum_{l \neq k} \frac{\dot{M}_{jl} \dot{M}_{lk}}{\delta_{kl}}
        \right)
        &= 
        \sum_{j, l \neq k} \frac{\dot{M}_{kj} \dot{M}_{jl} \dot{M}_{lk}}{\delta_{kj} \delta_{kl}} - \sum_{j \neq k} \frac{\dot{M}_{kj} \dot{M}_{jk}}{\delta_{kj}^2} \dot{M}_{kk}.
    \end{align*}
    By substituting, we obtain the announced expression for $\dddot{\lambda}_k$.
\end{proof}

In order to bound the term in $\dddot{\lambda}_k(\xi)$ in~\eqref{eq:DT_expan:gen_expan_TaylorLagrange}, we want to control 
\begin{equation}
    \gamma(\alpha) \coloneqq \min_{j \neq k} \abs{\lambda_j(\alpha) - \lambda_k(\alpha)}
    ~~~~\text{and}~~~~
    \chi(\alpha) \coloneqq \max_k \norm{v_k(\alpha)} \norm{w_k(\alpha)}
\end{equation}
the uniform eigengap resp.\ maximal eigenvalue condition number,
uniformly in $\alpha$ for $\alpha$ in a neighborhood of zero.
Now assuming $M_0$ has distinct eigenvalues then $\gamma(0) > 0$, and assuming $M_0$ is normal -- as is the case for \autoref{prop:DT_expan:gen_expan_formula} -- then $\chi(0) = 1$.
So we want to bound the perturbation of the eigenvalues and eigenvectors uniformly.

For the eigenvalues, thanks to the assumption that $M_0$ is normal, we easily get the following bound.
\begin{lemma} \label{lm:DT_expan:control_gamma}
    Under the conditions of \autoref{lm:DT_expan:restate_taoeigenvalues2008}, if additionally $M_0$ is normal, then it holds
    \begin{equation}
        \max_k \abs{\lambda_k(\alpha) - \lambda_k(0)} \leq \spnorm{M_\alpha - M_0}
    \end{equation}
    for all $\alpha \in I'$ the maximal open interval containing $0$ such that 
    $\forall \alpha \in I', \spnorm{M_\alpha - M_0} 
    < \frac{1}{2} \gamma(0)$.

    In particular, we have for all $\alpha \in I'$
    \begin{equation}
        \gamma(\alpha) \geq \gamma(0) - 2 \spnorm{M_\alpha - M_0} > 0.
    \end{equation}
\end{lemma}

Note that $I' \subset I$ the interval from \autoref{lm:DT_expan:restate_taoeigenvalues2008}, since $\gamma(\alpha) > 0$ if and only if $M_\alpha$ has distinct eigenvalues.

\begin{proof}
    Let $\lambda$ any eigenvalue of $M_\alpha$ and $v$ such that $M_\alpha v = \lambda v$ and $\norm{v}=1$. 
    Since $M_0$ is normal, it is orthonormally diagonalizable, so we may write it as $M_0 = U_0 \Lambda_0 \olU_0^\top$ with $U_0$ unitary and $\Lambda_0 = \Diag((\lambda_k(0)))$. 
    Then, letting $\tv = \olU_0^\top v$ and $\Delta = M_\alpha - M_0$ and $\widetilde{\Delta} = \olU_0^\top \Delta U_0$,
    \begin{gather*}
        M_\alpha v = (M_0 + \Delta) v
        = U_0^\top (\Lambda_0 + \widetilde{\Delta}) \olU_0^\top v
        = \lambda v,
        ~~\text{i.e.,}~~
        (\Lambda_0 + \widetilde{\Delta}) \tv = \lambda \tv \\
        \text{and so}~~
        \min_j \abs{\lambda - \lambda_{0j}}^2 \leq
        \sum_j \abs{\tv[j]}^2 \abs{\lambda - \lambda_{0j}}^2
        = \norm{(\Lambda_0 + \widetilde{\Delta}) \tv}^2
        \leq \norm{\widetilde{\Delta} \tv}^2 \leq \spnorm{\Delta}^2
    \end{gather*}
    since $U_0$ is unitary and $\norm{v}=1$.

    This shows that for any $\lambda \in \spectrum(M_\alpha)$, there exists a $\lambda_k(0) \in \spectrum(M_0)$ which is close to it, i.e.,
    \begin{equation}
        \forall k,~ \exists j~ \text{s.t}~~
        \abs{\lambda_k(\alpha) - \lambda_j(0)} \leq \spnorm{M_\alpha - M_0}.
    \end{equation}
    A fortiori, we can ensure that the smooth parametrization of the eigenvalues $(\lambda_k(\alpha))_k$ satisfies the inequality announced in the lemma, by restraining $\alpha$ to some $I'$ small enough so that $\argmin_j \abs{\lambda_k(\alpha) - \lambda_j(0)} = k$ for all $\alpha \in I'$.
    More explicitly, this can be achieved by choosing $I'$ such that
    $\sup_{\alpha \in I'} \spnorm{M_\alpha - M_0} 
    < \frac{1}{2} \min_{j \neq k} \abs{\lambda_j(0) - \lambda_k(0)}
    = \frac{1}{2} \gamma(0)$.
    Hence the choice of $I'$ announced in the lemma.
\end{proof}

For the eigenvectors, also using the assumption that $M_0$ is normal, we get the following bound.
\begin{lemma} \label{lm:DT_expan:control_chi}
    Under the conditions of \autoref{lm:DT_expan:restate_taoeigenvalues2008}, assume additionally that $M_0$ is normal, 
    and choose the normalization of the eigenvectors $(v_k(\alpha))_k$ such that $\norm{v_k(\alpha)}=1$ and $v_k^*(0) \cdot v_k(\alpha) \in \RR_+$ for all $k$. 
    Then
    \begin{align*}
        \forall k,~
        \norm{v_k(\alpha) - v_k(0)} &\leq \frac{2 \sqrt{2} \spnorm{M_\alpha - M_0}}{\gamma(0)} \\
        \text{and}~~~~
        \norm{w^*_k(\alpha) - v^*_k(0)} &\leq \sqrt{d} \norm{w_k(\alpha)} \frac{2 \sqrt{2} \spnorm{M_\alpha - M_0}}{\gamma(0)}
    \end{align*}
    for all $\alpha \in I'$ the interval from \autoref{lm:DT_expan:control_gamma}.

    In particular, we have
    \begin{equation}
        \chi(\alpha) 
        \leq 
        \frac{\gamma(0)}{\gamma(0) - 2 \sqrt{2 d} \spnorm{M_\alpha - M_0}}
    \end{equation}
    for all $\alpha \in I''$ the maximal open interval containing $0$ such that
    $\forall \alpha \in I'', \spnorm{M_\alpha - M_0} < \frac{1}{2 \sqrt{2d}} \gamma(0)$.
\end{lemma}

\begin{proof}
    We will write for concision $v_k = v_k(\alpha)$ and $v_{0k} = v_k(0) = w_k(0)$ since $M_0$ is normal, and $\lambda_k = \lambda_k(\alpha)$, $\lambda_{0k} = \lambda_k(0)$.
    Fix $\alpha \in I'$ and let $\Delta = M_\alpha - M_0$; by \autoref{lm:DT_expan:control_gamma} we have $\max_k \abs{\lambda_k - \lambda_{0k}} \leq \spnorm{\Delta} \leq \frac{1}{2} \gamma(0)$.

    Fix $k$. Subtracting $M_0 v_{0k} = \lambda_{0k} v_{0k}$ from $(M_0 + \Delta) v_k = \lambda_k v_k$, we have
    \begin{align*}
        M_0 (v_k - v_{0k}) + \Delta v_k 
        = \lambda_k v_k - \lambda_{0k} v_{0k}
        &= \lambda_{0k} (v_k-v_{0k}) + (\lambda_k - \lambda_{0k}) v_{0k} \\
        (M_0 - \lambda_{0k} I) (v_k-v_{0k}) &= -\Delta v_k + (\lambda_k - \lambda_{0k}) v_{0k} \\
        \sum_{j \neq k} \abs{\lambda_{0j}-\lambda_{0k}}^2 \abs{v_{0j}^* v_k}^2 = 
        \norm{(M_0 - \lambda_{0k} I) (v_k-v_{0k})}^2
        &\leq 
        \spnorm{\Delta}^2 ( \norm{v_k} + 1)^2
        = 4 \spnorm{\Delta}^2
    \end{align*}
    where $v_{0j}^* \coloneqq \olv_{0j}^\top$,
    using the unitary basis $(v_{0j})_j$ to compute the norm on the left-hand side.
    The left-hand side is further lower-bounded by $\gamma(0)^2 \sum_{j \neq k} \abs{v_{0j}^* v_k}^2 = \gamma(0)^2  \left( 1 - \abs{v_{0k}^* v_k}^2 \right)$,
    so
    \begin{equation}
        1 - \abs{v_{0k}^* v_k}^2
        \leq \frac{4 \spnorm{\Delta}^2}{\gamma(0)^2}.
    \end{equation}
    Consequently, by our choice of normalization: $v_k^*(0) v_k(\alpha) \in \RR_+$,
    we have
    \begin{equation}
        \norm{v_k-v_{0k}}^2 = 2 - 2 \Re(v_{0k}^* v_k)
        = 2 - 2 \abs{v_{0k}^* v_k} 
        \leq 2 \left( 1 - \sqrt{1 - \frac{4 \spnorm{\Delta}^2}{\gamma(0)^2}} \right)
        \leq \frac{8 \spnorm{\Delta}^2}{\gamma(0)^2}
    \end{equation}
    using that $\forall y \in [0,1],~ 1 - \sqrt{1-y^2} \leq y^2$.
    This shows the first inequality of the lemma.

    Let us now show the control on $w_k \coloneqq w_k^*(\alpha)$ -- the second inequality of the lemma.
    Using the unitary basis $(v_{0j})_j$ to compute the norm, we have 
    $\norm{w^*_k - v^*_{0k}}^2 = \sum_{j \neq k} \abs{w^*_k v_{0j}}^2 + \abs{w^*_k v_{0k} - 1}^2$ 
    and since by definition $w_k^* v_j = \ind_{j=k}$, we can bound each term as
    \begin{align*}
        \forall j \neq k,~
        \abs{w_k^* v_{0j}} &\leq \abs{w_k^* v_j} + \abs{w_k^* (v_{0j} - v_j)}
        \leq 0 + \norm{w_k^*} \frac{2 \sqrt{2} \spnorm{\Delta}}{\gamma(0)} \\
        \text{and}~~~~
        \abs{w^*_k v_{0k} - 1}
        &= \abs{w^*_k v_k - 1 + w^*_k (v_{0k} - v_k)}
        \leq \norm{w^*_k} \norm{v_{0k} - v_k}
        \leq \norm{w_k^*} \frac{2 \sqrt{2} \spnorm{\Delta}}{\gamma(0)}.
    \end{align*}
    So in total,
    $\norm{w^*_k - v^*_{0k}}^2 \leq d \cdot \norm{w_k^*}^2 \frac{8 \spnorm{\Delta}^2}{\gamma(0)^2}$,
    as announced.

    The second part of the lemma follows by noting that $\norm{v_k} \norm{w_k} = 1 \cdot \norm{w_k}$ is bounded by
    \begin{align*}
        \norm{w_k} \leq \norm{v_{0k}} + \norm{w_k - v_{0k}} 
        &\leq 1 + \sqrt{d} \norm{w_k} \frac{2 \sqrt{2} \spnorm{\Delta}}{\gamma(0)} \\
        \implies~~
        \norm{w_k} &\leq \frac{1}{1 - \sqrt{d} \cdot \frac{2 \sqrt{2} \spnorm{\Delta}}{\gamma(0)}} = \frac{\gamma(0)}{\gamma(0) - 2 \sqrt{2} \sqrt{d} \spnorm{\Delta}}
    \end{align*}
    for all $\alpha$ such that the denominator in the second line is positive,
    as announced.
\end{proof}

By combining the three above lemmas, we have that under the conditions of \autoref{lm:DT_expan:restate_taoeigenvalues2008}, if additionally $M_0$ is normal, then
for all $\alpha$ in a neighborhood of zero such that
$\spnorm{M_\alpha - M_0} \leq \frac{\gamma(0)}{4 \sqrt{2d}}$, one can check that
$\gamma(\alpha) \geq \frac{1}{2} \gamma(0)$ and $\chi(\alpha) \leq 2$ and so
\begin{equation}
    \frac{1}{6} \max_k \abs{\dddot{\lambda_k}} 
    \leq \frac{1}{3} \spnorm{\dddot{M}} 
    + 8 d \gamma(0)^{-1} \spnorm{\dot{M}} \cdot \spnorm{\ddot{M}} 
    + 32 d^2 \gamma(0)^{-2} \spnorm{\dot{M}}^3.
\end{equation}
In particular, for a matrix $M_\alpha = M_0 + \alpha M_1 + \frac{\alpha^2}{2} M_2$ with $M_0$ normal, the right-hand side translates exactly to the bound on $\abs{\br_j}$ announced in \autoref{prop:DT_expan:gen_expan_formula}.
So that proposition follows directly from \eqref{eq:DT_expan:gen_expan_TaylorLagrange} and the above discussion.



\section{Local convergence of equality-constrained Mirror Flow} \label{apx:MF}

This appendix contains results used in \appendixref{apx:details_particleMNE}.

Let
$\XXX, \YYY$ convex subsets of $\RR^n$ resp.\ $\RR^m$ and $\phi_x: \XXX \to \RR$, $\phi_y: \YYY \to \RR$ strictly convex and differentiable.
Let $\ZZZ = \XXX \times \YYY$ and $\phi: \ZZZ \to \RR$ with $\phi(x, y) = \phi_x(x) + \phi_y(y)$.
Consider a twice continuously differentiable min-max objective $f: \XXX \times \YYY \to \RR$, denote 
$g(z) = \Diag(I_n, -I_m) \nabla f(z)$ and $M(z) = \nabla g(z)$.
Throughout this appendix 
we make the following assumption.
\begin{assumption} \label{assum:MF:eq_constr}
    The constraint set is defined by equalities:
    ${ 
        \ZZZ = \left\lbrace z \in \RR^{n+m}; Az = b \right\rbrace = z_b + \Ker A 
    }$,
    where $z_b$ is any solution of $A z_b=b$.
    Furthermore, $\phi$ is 
    strictly convex and
    three times differentiable.
\end{assumption}

\begin{definition}[{\cite[Proposition~1]{amid_reparameterizing_2020}}] \label{def:MF:MF}
    For an initial point $z^0 \in \ZZZ$, the mirror flow (MF) with link function $\phi$ is the unique curve $z(t)$ such that $z(0) = z^0$ and
    \begin{align} \label{eq:MF:MF}
        \frac{dz}{dt} = -\Phi^{-1}_z P_z g(z) 
        \eqqcolon -g_\eff(z),
    \end{align}
    where $\Phi_z$ denotes the Hessian of $\phi$ at $z$ and
    $
        P_z \coloneqq I - A^\top \left[ A \Phi^{-1}_z A^\top \right]^{-1} A \Phi^{-1}_z
    $.
\end{definition}

One can check that $P_z^2 = P_z$ and that 
$P_z \Phi_z = \Phi_z P_z^\top$.
Furthermore, $P_z A^\top = 0$,
i.e.,
$\Ima P_z^\top \subset \Ker A$,
and in fact $\Ima P_z^\top = \Ker A$ since $(\Ima P_z^\top)^\perp = \Ker P_z = \Ima(I-P_z) \subset \Ima A^\top = (\Ker A)^\perp$.
%
In particular
the MF preserves the constraint set, since
${ \frac{d}{dt} (Az-b) = -A P^\top \Phi^{-1}_z g(z) = 0 }$.
As a consequence, $g_\eff$ can be seen as an operator from the affine space $\ZZZ = z_b + \Ker A$ to itself.

\begin{lemma} \label{lm:MF:jac_g_eff}
    The Jacobian of $g_\eff: \RR^{n+m} \to \RR^{n+m}$ at a local NE $z^*$ is equal to
    $
        \Phi^{-1}_{z^*} P_{z^*} \cdot M(z^*)
    $.

    Furthermore, the Jacobian of $g_\eff: \ZZZ \to \ZZZ$ (seen as an operator between affine spaces) at $z^*$ is
    \begin{equation}
        M_\eff(z^*) = \Phi^{-1}_{z^*} P_{z^*} \cdot M(z^*) P_{z^*}^\top
    \end{equation}
    which is a linear operator from $\Ker A$ to itself.
\end{lemma}

\begin{proof}
    In this proof, write for concision $P = P_z$ and $\Phi^{-1} = \Phi^{-1}_z$.
    Using Einstein's summation notation, pose
	\begin{equation}
		P\indices{_i^j} = P = I - A^\top \left[ A \Phi^{-1} A^\top \right]^{-1} A \Phi^{-1} 
        \qquad \text{and} \qquad
		Q^{ij} = Q \coloneqq \Phi^{-1} P.
	\end{equation}
    In particular $g_\eff(z) = Q g$.
    Now
	\begin{equation} \label{eq:MF:nabla_Qg}
		\nabla \{ Q g \}\indices{^i_j} 
        = \frac{\partial [Q g]^i}{\partial z^j}
		= (Q \nabla g)\indices{^i_j} + \frac{\partial Q^{ik}}{\partial z^j} g_k.
	\end{equation}
	Using formula~(59) from \cite{petersen_cookbook_2008}:
    $
        \frac{\partial Y^{-1}}{\partial x} = -Y^{-1} \frac{\partial Y}{\partial x} Y^{-1}
    $,
	we have that
	\begin{align}
		\frac{\partial (\Phi^{-1})^{il}}{\partial z^j} &= -(\Phi^{-1})\indices{^{ia}} K_{abj} (\Phi^{-1})\indices{^{bl}} \\
		\frac{\partial P\indices{_l^k}}{\partial z^j} 
		&= - \left(
			A^\top [A \Phi^{-1} A^\top]^{-1} A
		\right)\indices{_l^s} 
		K_{stj}
		\left(
			\Phi^{-1} A^\top [A \Phi^{-1} A^\top]^{-1} A \Phi^{-1}
		\right)^{tk} \\
		&~~~~ - (A^\top)_{la}
		\frac{\partial \{ [A \Phi^{-1} A^\top]^{-1} \}\indices{^{ab}} }{\partial z^j}
		(A \Phi^{-1})\indices{_b^k}
	\end{align}
    where
    $K = \nabla^3 \phi(z)$,
	and after calculating and simplifying,
	\begin{align}
		\frac{\partial P\indices{_l^k}}{\partial z^j} 
		&= (I-P)\indices{_l^s} K_{stj} (\Phi^{-1} P)^{tk}.
	\end{align}
	So
	\begin{align}
		\frac{\partial Q^{ik}}{\partial z^j}
		&= (\Phi^{-1})^{il} \frac{\partial P\indices{_l^k}}{\partial z^j}
		+ \frac{\partial (\Phi^{-1})^{il}}{\partial z^j} P\indices{_l^k} \\
		&= -(\Phi^{-1} P)^{is} K_{stj} (\Phi^{-1} P)^{tk}
		= -Q^{is} K_{stj} Q^{tk} \\
        \text{and finally}~~~~
		\frac{\partial Q^{ik}}{\partial z^j} g_k
		&= -Q^{is} K_{stj} (Q g)^t.
	\end{align}
    Now $g_\eff(z^*) = 0$, so $Qg=0$ at $z^*$.
	The first part of the lemma follows by substituting into \eqref{eq:MF:nabla_Qg}.

    To check the second part of the lemma, simply use that $\Phi^{-1}_{z^*} P_{z^*} = P_{z^*}^\top \Phi^{-1}_{z^*}$ and that $P_{z^*}^\top$ is a projector onto the kernel of $A$,
    as remarked under \definitionref{def:MF:MF}.
\end{proof}

The next lemma
is a special case of the stable manifold theorem \citep[Section~2.7]{perko_differential_2013}.
\begin{lemma} \label{lm:MF:suffcond_cv_conttime}
    Let $\ZZZ = z_b + T\ZZZ$ an affine subspace of $\RR^d$. Let a continuously differentiable operator $g: \ZZZ \to \ZZZ$ and denote its Jacobian as $\nabla g: \ZZZ \to T\ZZZ \times T\ZZZ$.
    Suppose there exists $z^*$ such that $g(z^*)=0$ and $\min_{\lambda \in \spectrum(\nabla g(z^*))} \Re(\lambda) > 0$.

    Then for any $\eps>0$, there exists $C>0$ and a relative neighborhood of $z^*$ such that, for any $z^0$ in that neighborhood, the flow $\frac{dz}{dt} = -g(z)$ converges exponentially to $z^*$ with
    \begin{equation}
        \norm{z(t)-z^*}^2 \leq C \norm{z^0-z^*}^2 \exp\left( -2\left[ \min_{\lambda \in \spectrum(\nabla g(z^*))} \Re(\lambda) - \eps \right] t \right).
    \end{equation}
\end{lemma}

\begin{proposition} \label{prop:MF:suffcond_MF}
    If $\tM \coloneqq \Phi^{-\frac{1}{2}}_{z^*} P_{z^*} \cdot M(z^*) \cdot P_{z^*}^\top \Phi^{-\frac{1}{2}}_{z^*}$ is invertible and satisfies the equivalent conditions of \theoremref{thm:quadr:eigvals},
    then MF converges locally exponentially to $z^*$
    at a rate $\tmu_{\tM}$.
\end{proposition}

\begin{proof}
    $\tM$ is similar to $M_\eff(z^*)$. In particular they have the same eigenvalues.
    So the proposition follows immediately from the above two lemmas.
\end{proof}

\begin{remark} \label{rk:MF:interpret_tM_R}
    The matrix $\tM$ can be interpreted as 
    \begin{equation}
        \tM =
        \underbrace{
            \Phi^{-\frac{1}{2}}_{z^*} P_{z^*} \Phi^{\frac{1}{2}}_{z^*} 
        }_{\eqqcolon R_{z^*}}
        ~\cdot~
        \Phi^{-\frac{1}{2}}_{z^*} M(z^*) \Phi^{-\frac{1}{2}}_{z^*} 
        ~\cdot~ 
        \underbrace{
            \Phi^{\frac{1}{2}}_{z^*} P_{z^*}^\top \Phi^{-\frac{1}{2}}_{z^*} 
        }_{= R_{z^*}^\top}.
    \end{equation}
    Note that $R_{z^*}$ is an orthogonal projection.
    Furthermore, the central factor consists in a transformation of $M(z^*)$ that is compatible with the block structure:
    \begin{itemize}[leftmargin=4mm]
        \item If we write $M(z^*) = S+A$ with $S$ symmetric and $A$ antisymmetric, then 
        $\Phi^{-\frac{1}{2}}_{z^*} M(z^*) \Phi^{-\frac{1}{2}}_{z^*} 
        = \Phi^{-\frac{1}{2}}_{z^*} S \Phi^{-\frac{1}{2}}_{z^*} 
        + \Phi^{-\frac{1}{2}}_{z^*} A \Phi^{-\frac{1}{2}}_{z^*}$
        and the first term is symmetric and the second term is antisymmetric.
        \item If we write
        $M(z^*) = \begin{bmatrix}
            Q & P \\ -P^\top & R
        \end{bmatrix}$ 
        with $Q, R$ symmetric,
        then \\
        ${
            \Phi^{-\frac{1}{2}}_{z^*} M(z^*) \Phi^{-\frac{1}{2}}_{z^*} 
            = \begin{bmatrix}
                \Phi^{-\frac{1}{2}}_{x^*} & \bmzero \\ \bmzero & \Phi^{-\frac{1}{2}}_{y^*}
            \end{bmatrix}
            \begin{bmatrix}
                Q & P \\ -P^\top & R
            \end{bmatrix}
            \begin{bmatrix}
                \Phi^{-\frac{1}{2}}_{x^*} & \bmzero \\ \bmzero & \Phi^{-\frac{1}{2}}_{y^*}
            \end{bmatrix}
            = \begin{bmatrix}
                \Phi^{-\frac{1}{2}}_{x^*} Q \Phi^{-\frac{1}{2}}_{x^*} 
                & \Phi^{-\frac{1}{2}}_{x^*} P \Phi^{-\frac{1}{2}}_{y^*} \\
                -\left( \Phi^{-\frac{1}{2}}_{x^*} P \Phi^{-\frac{1}{2}}_{y^*} \right)^\top 
                & \Phi^{-\frac{1}{2}}_{y^*} R \Phi^{-\frac{1}{2}}_{y^*}
            \end{bmatrix}
        }$
        is of the same form,
        where $\Phi_{x^*}$ is the Hessian of $\phi_x$ at $x^*$ and similarly for $\Phi_{y^*}$.
    \end{itemize}
\end{remark}

We state an analogous result for Mirror Descent-Ascent (MDA) and its variants, Mirror Prox (MP) and Bregman PP.
Its proof, omitted for brevity, involves exactly the same ideas as for MF.
\begin{proposition} \label{prop:MF:suffcond_MDA_BregmanPP}
    Under the assumptions of \propositionref{prop:MF:suffcond_MF}, 
    if additionally $\phi$ is strongly convex,
    then MDA, MP and Bregman PP
    converge locally exponentially to $z^*$ for any small enough step-size
    at a rate $\eta \tmu_{\tM} + O(\eta^2)$.
\end{proposition}



\section{Details for \texorpdfstring{\sectionref{sec:particleMNE}}{Section~\ref{sec:particleMNE}}} \label{apx:details_particleMNE}

In this appendix we 
provide details for \figureref{fig:particleMNE:manyalphas}, 
derive the first-order term (in $\eta$) in the local convergence rate of MP and EG, 
and discuss a counter-example.
For our purpose it is sufficient to consider only the continuous-time flows of those algorithms.

\subsection{Details for \texorpdfstring{\figureref{fig:particleMNE:manyalphas}}{Figure~\ref{fig:particleMNE:manyalphas}}} \label{subsec:details_particleMNE:for_figure}

In the numerical experiment reported in \figureref{fig:particleMNE:manyalphas},
we used the same random payoff functions as in \cite[Section~4]{wang_exponentially_2022}:
\begin{equation} \label{eq:details_particleMNE:def_rfpayoffs}
    f(x, y) = \Re \sum_{-K \leq k \leq K} \sum_{-L \leq l \leq L} c_{kl} e^{2\pi i (kx + ly)}
    ~~~~ \text{for $(c_{kl}) \in \CC^{(2K+1) \times (2L+1)}$,}
\end{equation}
with $K=L=2$ and $\Re(c_{kl})$, $\Im(c_{kl})$ drawn independently from the standard Gaussian distribution.
We measured the distance between the solution $z^*=(a^*, x^*, b^*, y^*)$ and the iterates $z^k$ by
\begin{equation}
    \norm{z-z^*}^2
    \coloneqq \norm{a-a^*}^2 + \norm{b-b^*}^2 + \norm{x-x^*}^2 + \norm{y-y^*}^2.
\end{equation}
For $z$ in a neighborhood of $z^*$,
$\norm{z-z^*}$
is an upper bound on the duality gap of the original MNE problem up to a constant dependent only on $f$
\citep[Proposition~3.1 and Claim~C.3]{wang_exponentially_2022}.
%

\subsection{First-order term in the local convergence rate of MP} \label{subsec:details_particleMNE:MP}

Denote the joint variable as $z = (a, b) \in \Delta_N \times \Delta_M$ and let $z^*$ such that $(\mu^*, \nu^*) = \left( \sum_I a^*_I \delta_{x^*_I}, \sum_J b^*_J \delta_{y^*_J} \right)$ is the MNE.
The continuous-time flow of MP is MF (\definitionref{def:MF:MF})
with
$\phi(z) = \sum_I z_I \log z_I$, 
$A = \begin{bmatrix}
    \bmone_N^\top & 0 \\
    0 & \bmone_M^\top
\end{bmatrix} \in \RR^{2 \times (N+M)}$
and $b= \begin{pmatrix} 1 \\ 1 \end{pmatrix}$.
So by \propositionref{prop:MF:suffcond_MDA_BregmanPP}, 
the first-order term in the local convergence rate of MP is
$\eta \tmu_{M_{\MP}}$ for some $M_{\MP}$ which we now compute.

Following the notations of \appendixref{apx:MF},
we have
$\Phi_{z^*} =
\Diag\left(
    \left( \frac{1}{a_I^*} \right)_I,
    \left( \frac{1}{b_J^*} \right)_J
\right)$,
and by straightforward calculations
$P_{z^*} = \Diag\left( I_N - \bmone_N (a^*)^\top, I_M - \bmone_M (b^*)^\top \right)$,
so that
$R_{z^*} = \Phi_{z^*}^{-\half} P_{z^*} \Phi_{z^*}^{\half} 
= \Diag\left( I_N-\sqrt{a^*} \sqrt{a^*}^\top, I_M-\sqrt{b^*}\sqrt{b^*}^\top \right)$.
In order to project out the constraints explicitly, 
let $\Pi_a$ any matrix in $\RR^{(N-1) \times N}$ with orthonormal rows, i.e., $\Pi_a \Pi_a^\top = I_{N-1}$, and such that $\Pi_a^\top \Pi_a  = I_N - \sqrt{a^*} \sqrt{a^*}^\top$, and likewise for $\Pi_b \in \RR^{(M-1) \times M}$.
For example the rows of $\Pi_a$ may be obtained by completing the unit vector $\sqrt{a^*}$ into an orthonormal basis and removing $\sqrt{a^*}$ from the basis.
Then $R_{z^*}$ as a linear projector from $\RR^{N+M}$ to $\{ \sqrt{a^*} \}^\perp \times \{ \sqrt{b^*} \}^\perp \simeq \RR^{N+M-2}$ can be written as the matrix $R_{z^*} = \Diag(\Pi_a, \Pi_b)$.%
\footnote{
    The expression ``$\Diag(\Pi_a, \Pi_b)$'', as well as ``$\Diag\left( \Pi_a, I_N, \Pi_b, I_M \right)$'' in \autoref{subsec:details_particleMNE:CPMP}, constitutes a slight abuse of our notation ``$\Diag$'', since $\Pi_a$ and $\Pi_b$ are not square matrices. To remove any ambiguity: by $\Diag(\Pi_a, \Pi_b)$ we mean the matrix 
    $\begin{bmatrix}
        \Pi_a & \bmzero_{(N-1) \times M} \\
        \bmzero_{(M-1) \times N} & \Pi_b
    \end{bmatrix}$,
    and similarly for $\Diag\left( \Pi_a, I_N, \Pi_b, I_M \right)$.
}
So $\tM$ as a linear operator over $\{ \sqrt{a^*} \}^\perp \times \{ \sqrt{b^*} \}^\perp$ can be written as the matrix, denoting $D_a = \Diag(\sqrt{a^*})$ and $D_b = \Diag(\sqrt{b^*})$ for concision,
\begin{align*}
    M_{\MP} &\coloneqq R_{z^*} \Phi_{z^*}^{-\half} M(z^*) \Phi_{z^*}^{-\half} R_{z^*}^\top \\
    &= 
    \Diag\left( \Pi_a, \Pi_b \right) 
    \Diag\left( \sqrt{a^*}, \sqrt{b^*} \right)
    \begin{bmatrix}
        \bmzero & P \\
        -P^\top & \bmzero
    \end{bmatrix}
    \Diag\left( \sqrt{a^*}, \sqrt{b^*} \right)
    \Diag\left( \Pi_a^\top, \Pi_b^\top \right) \\
    &= \begin{bmatrix}
        \bmzero & \Pi_a D_a P D_b \Pi_b^\top \\
        -\Pi_b D_b P^\top D_a \Pi_a^\top & \bmzero
    \end{bmatrix}.
\end{align*}

\subsection{First-order term in the local convergence rate of EG (conic particle methods)}
\label{subsec:details_particleMNE:CPMP}

As in the main text, we write
``EG'' to refer to the Conic Particle Mirror Prox algorithm of \cite{wang_exponentially_2022}.

Denote the joint variable as $z = (a, x, b, y) \in \Delta_N \times (\TT^1)^N \times \Delta_M \times (\TT^1)^M$ and let any $z^*$ such that $(\mu^*, \nu^*) = \left( \sum_I a^*_I \delta_{x^*_I}, \sum_J b^*_J \delta_{y^*_J} \right)$ is the MNE.
The flow is given by the system of ODEs
\begin{align}
\label{eq:exp:CP_MF}
        \frac{da}{dt} &= -\Phi_a^{-1} P_a \nabla_a F(z) 
        & &\text{and} &
        \frac{db}{dt} &= \Phi_b^{-1} P_b \nabla_b F(z) \\
        \frac{dx}{dt} &= -\gamma \Diag\left( \frac{1}{a_I} \right) \nabla_x F(z)
        & & &
        \frac{dy}{dt} &= \gamma \Diag\left( \frac{1}{b_J} \right) \nabla_y F(z).
\end{align}
Here
$F(z) = \sum_I \sum_J a_I b_J f(x_I, y_J)$,
$\gamma$ is a constant parameter,
$\Phi_a = \Diag\left( \frac{1}{a_I} \right)$
and
$P_a = I_N - \bmone_N a^\top$,
and $\Phi_b$ and $P_b$ are defined similarly.

In general, the flow~\eqref{eq:exp:CP_MF} does not match the structure of Mirror Flow because of the factor ``$\Diag(\frac{1}{a_I})$'' in the equation for $\frac{dx}{dt}$ 
\citep[Section~2.4]{gunasekar_mirrorless_2021}.
However, by adapting 
the reasoning of \appendixref{apx:MF} -- namely only \lemmaref{lm:MF:jac_g_eff} needs to be adapted --
it is easy to show that the statement of \propositionref{prop:MF:suffcond_MF} holds also for this dynamics.
By the same adaptation 
one can show
that the statement of \propositionref{prop:MF:suffcond_MDA_BregmanPP} holds also for EG.
Hence the first-order term in the local convergence rate of EG is $\eta \tmu_{M_{\gamma}}$ 
for some $M_{\gamma}$ which we now compute.

The statement of \propositionref{prop:MF:suffcond_MF} applies to
$\tM = R_{z^*} \cdot \Phi_{z^*}^{-\half} M(z^*) \Phi_{z^*}^{-\half} \cdot R_{z^*}^\top$
where
$\Phi_{z^*} =
\Diag\left(
    \left( \frac{1}{a_I^*} \right)_I,
    \frac{1}{\gamma} a^*,
    \left( \frac{1}{b_J^*} \right)_J,
    \frac{1}{\gamma} b^*
\right)$
and
$P_{z^*} = \Diag\left( P_{a^*}, I_N, P_{b^*}, I_M \right)$
so that
$R_{z^*} = \Phi_{z^*}^{-\half} P_{z^*} \Phi_{z^*}^{\half} 
= \Diag\left( I_N-\sqrt{a^*} \sqrt{a^*}^\top, I_N, I_M-\sqrt{b^*}\sqrt{b^*}^\top, I_M \right)$.
In order to project out the constraints explicitly, 
let $\Pi_a \in \RR^{(N-1) \times N}$, $\Pi_b \in \RR^{(M-1) \times M}$ the same matrices as in the previous subsection.
Then $R_{z^*}$ as a linear projector from $\RR^{2N+2M}$ to $\{ \sqrt{a^*} \}^\perp \times \RR^N \times \{ \sqrt{b^*} \}^\perp \times \RR^M \simeq \RR^{2N+2M-2}$ can be written as the matrix
$R_{z^*} = \Diag\left( \Pi_a, I_N, \Pi_b, I_M \right)$.
Moreover by \cite[Claim~C.2]{wang_exponentially_2022}, dropping superscript *'s for concision only in this equation,
\begin{align*}
    M(z^*)
    &= \begin{bmatrix}
        \bmzero & \bmzero & P & \partial_y P \Diag(b) \\
        \bmzero & \Diag(a) \Diag(\partial_{xx}^2 P b) & \Diag(a) \partial_x P & \Diag(a) \partial_{xy}^2 P \Diag(b) \\
        -P^\top & -\left( \Diag(a) \partial_x P \right)^\top & \bmzero & \bmzero \\
        -\left( \partial_y P \Diag(b) \right)^\top & -\left( \Diag(a) \partial_{xy}^2 P \Diag(b) \right)^\top & \bmzero & -\Diag(b) \Diag(\partial_{yy}^2 P^\top a)
    \end{bmatrix}
\end{align*}
where $\left[ \partial_x P \right]_{IJ} = \partial_x f(x^*_I, y^*_J)$, and likewise for $\partial_y P$, $\partial_{xx}^2 P$, $\partial_{yy}^2 P$, $\partial_{xy}^2 P$.
So, finally, $\tM = R_{z^*} \Phi_{z^*}^{-\half} M(z^*) \Phi_{z^*}^{-\half} R_{z^*}^\top$ as a linear operator over $\{ \sqrt{a^*} \}^\perp \times \RR^N \times \{ \sqrt{b^*} \}^\perp \times \RR^M$ can be written as the matrix,
denoting $D_a = \Diag(\sqrt{a^*})$ and $D_b = \Diag(\sqrt{b^*})$ for concision,
\begin{align*}
    M_{\gamma}
    &\coloneqq \begin{bmatrix}
        \bmzero & \bmzero & \Pi_a D_a P D_b \Pi_b^\top & \sqrt{\gamma} \ \Pi_a D_a [\partial_y P] D_b \\
        \bmzero & \gamma \Diag(\partial_{xx}^2 P b^*) & \sqrt{\gamma} D_a [\partial_x P] D_b \Pi_b^\top & \gamma D_a [\partial_{xy}^2 P] D_b \\
        -(*)^\top & -\left( * \right)^\top & \bmzero & \bmzero \\
        -\left( * \right)^\top & -\left( * \right)^\top & \bmzero & -\gamma \Diag(\partial_{yy}^2 P^\top a^*)
    \end{bmatrix}.
\end{align*}

\subsection{An example where the first-order term in the local convergence rate of EG is zero}
\label{subsec:details_particleMNE:counter_example}

Consider the payoff function defined by 
\eqref{eq:details_particleMNE:def_rfpayoffs}
with
$c_{20} = c_{02} = -i$,
$c_{11} = 2$ and
$c_{kl} = 0$ otherwise,
i.e.,
$
    f(x, y) = \sin(4\pi x) + \sin(4\pi y) + 2 \cos(2\pi x + 2\pi y)
$.
As shown in \cite[Example~4.1]{wang_exponentially_2022},
the MNE is unique
and given by
$a^* = b^* = \left( \frac{1}{2}, \frac{1}{2} \right)$,
$x^* = \left( \frac{3}{8}, \frac{7}{8} \right)$,
and $y^* = \left( \frac{1}{8}, \frac{5}{8} \right)$.
So we can compute $M_{\gamma}$ explicitly in this case: 
we find
\begin{equation}
    P = \begin{pmatrix}
        -2 & 2 \\
        2 & -2
    \end{pmatrix},~
    \partial_x P = \partial_y P = 0,~
    \partial_{xx}^2 P \begin{pmatrix} \half \\ \half \end{pmatrix}
    = \partial_{yy}^2 P \begin{pmatrix} \half \\ \half \end{pmatrix} = \begin{pmatrix} 16 \pi^2 \\ 16 \pi^2 \end{pmatrix},~
    \partial_{xy} P = \begin{pmatrix}
        8 \pi^2 & -8 \pi^2 \\
        -8 \pi^2 & 8 \pi^2
    \end{pmatrix}
\end{equation}
and
$D_a = D_b = \frac{1}{\sqrt{2}} I$,
and so
(for a certain choice of $\Pi_a$ and $\Pi_b$, each of which is anyway determined up to a sign)
\begin{align*}
    M_{\gamma} =
    \left[\begin{array}{c c c | c c c}
        0 & 0 & 0 & -2 & 0 & 0 \\
        0 & \gamma (4\pi)^2 & 0 & 0 & \gamma (2\pi)^2 & -\gamma (2\pi)^2 \\
        0 & 0 & \gamma (4\pi)^2 & 0 & -\gamma (2\pi)^2 & \gamma (2\pi)^2 \\
        \hline
        2 & 0 & 0 & 0 & 0 & 0 \\
        0 & -\gamma (2\pi)^2 & \gamma (2\pi)^2 & 0 & \gamma (4\pi)^2 & 0\\
        0 & \gamma (2\pi)^2 & -\gamma (2\pi)^2 & 0 & 0 & \gamma (4\pi)^2
    \end{array}\right].
\end{align*}
This matrix clearly does not satisfy condition $(iii)$ of \theoremref{thm:quadr:eigvals}, so 
$\tmu_{M_{\gamma}} = 0$.

\subsection{Proof of \autoref{prop:particleMNE:scale_gamma2}}
\label{subsec:details_particleMNE:proof_scale_gamma2}

For ease of reference, we restate the proposition below.
\particleMNEscaleGammaSq*

Pose $\alpha = \sqrt{\gamma}$, $M_0 = A_0 + A_\eps$, $M_1 = A_1$ and $M_2 = 2 (A_2 + S_2)$,
where $A_\eps$ is any antisymmetric matrix such that $\spnorm{A_\eps} \leq \eps$ and $M_0$ has distinct eigenvalues.
We will prove the proposition by applying the spectral expansions of \autoref{subsec:DT_expan:proof_gen_expan_formula} to the matrix curve $M_\alpha = M_0 + \alpha M_1 + \frac{\alpha^2}{2} M_2$.

Adopting the notations of that section, we have the expansion for the eigenvalues $\lambda_k(\alpha)$ of $M_\alpha$:
\begin{align}
    \lambda_k(\alpha) = \lambda_k(0) + \alpha \dot{\lambda}_k(0) + \frac{\alpha^2}{2} \ddot{\lambda}_k(0) + \frac{\alpha^3}{3!} \dddot{\lambda}_k(0) + \frac{\alpha^4}{4!} \ddddot{\lambda}_k(0) + O(\alpha^5),
\end{align}
with $\{ \lambda_k(0) \}_k = \spectrum(M_0) \subset i \RR$ by antisymmetry, 
$\dot{\lambda}_k(0) = v^*_{0k} A_1 v_{0k} \in i \RR$ -- where $v_{0k}$ are the eigenvectors of $M_0$ -- by normality of $M_0$ and antisymmetry of $A_1$, and
\begin{align}
    \ddot{\lambda}_k(0) = 2 v^*_{0k} (A_2 + S_2) v_{0k} 
    + 2 \sum_{j \neq k}~ \underbrace{
        \frac{(v^*_{0k} A_1 v_{0j}) (v^*_{0j} A_1 v_{0k})}{\lambda_k(0) - \lambda_j(0)}
    }_{\in i \RR}
\end{align}
since $v^*_{0k} A_1 v_{0j} = -(v^*_{0j} A_1 v_{0k})^*$,
and
\begin{multline*}
    \dddot{\lambda}_k(0) =
    3 \sum_{j \neq k} \frac{1}{\lambda_k(0)-\lambda_j(0)}
    \left[
        (v^*_{0k} A_1 v_{0j}) (v^*_{0j} M_2 v_{0k}) 
        + (v^*_{0k} M_2 v_{0j}) (v^*_{0j} A_1 v_{0k})
    \right] \\
    + 6~
    \underbrace{
        \sum_{j, l \neq k} \frac{(v^*_{0k} A_1 v_{0j}) (v^*_{0j} A_1 v_{0l}) (v^*_{0l} A_1 v_{0k})}{(\lambda_k(0)-\lambda_j(0)) (\lambda_k(0)-\lambda_l(0))}
    }_{\in i \RR}
    ~-~ 6 \sum_{j \neq k}~ 
    \underbrace{
         \frac{(v^*_{0k} A_1 v_{0j}) (v^*_{0j} A_1 v_{0k})}{(\lambda_k(0)-\lambda_j(0))^2} v^*_{0k} A_1 v_{0k}
    }_{\in i \RR}
\end{multline*}
as one can check by computing the convex conjugate of each of the underbraced expressions.

Now, let $\KKK$ the set of indices corresponding to the non-zero eigenvalues of $A_0$, and ${ \Kc = \{1,...,d\} \setminus \KKK }$. Note that the eigenvectors of $A_0$ are of the form
\begin{align}
    \forall k \in \KKK,~ 
    \tv_k = \begin{pmatrix}
        * \\ 0 \\  * \\ 0
    \end{pmatrix}
    ~~~~\text{and}~~~~
    \forall h \in \Kc,~
    \tv_h = \begin{pmatrix}
        0 \\ * \\ 0 \\  *
    \end{pmatrix}.
\end{align}
Furthermore, note that
\begin{align}
    A_1 \begin{pmatrix}
        * \\ 0 \\  * \\ 0
    \end{pmatrix}
    &= \begin{pmatrix}
        0 \\ * \\ 0 \\  *
    \end{pmatrix},
&
    S_2 \begin{pmatrix}
        * \\ 0 \\  * \\ 0
    \end{pmatrix}
    &= \bmzero,
&
    A_2 \begin{pmatrix}
        * \\ 0 \\  * \\ 0
    \end{pmatrix}
    &= \bmzero,
\\
    A_1 \begin{pmatrix}
        0 \\ * \\ 0 \\  *
    \end{pmatrix}
    &= \begin{pmatrix}
        * \\ 0 \\  * \\ 0
    \end{pmatrix},
&
    A_2 \begin{pmatrix}
        0 \\ * \\ 0 \\  *
    \end{pmatrix}
    &= \begin{pmatrix}
        0 \\ * \\ 0 \\  *
    \end{pmatrix},
&
    A_2 \begin{pmatrix}
        0 \\ * \\ 0 \\  *
    \end{pmatrix}
    &= \begin{pmatrix}
        0 \\ * \\ 0 \\  *
    \end{pmatrix}.
\end{align}
Using this structure, one can check that for all $k \in \KKK$,
$\Re \ddot{\lambda}_k(0) = o_\eps(1)$
and
$\Re \dddot{\lambda}_k(0) = o_\eps(1)$.

Thus, the spectrum of $M_\alpha$ consists of eigenvalues (corresponding to indices $h \in \Kc$) with $\Re \ddot{\lambda}_h(0)$ non-zero a priori, in which case $\Re \lambda_h(\alpha) = \Theta(\alpha^2) = \Theta(\gamma)$,
and of eigenvalues (corresponding to indices $k \in \KKK$) with $\Re \lambda_k(\alpha) = O(\alpha^4) + o_\eps(1) = O(\gamma^2) + o_\eps(1)$.
By letting $\eps \to 0$ and using that eigenvalues are continuous, this shows that $\tmu_{M_\gamma} = O(\gamma^2)$ as $\gamma \to 0$.


\fi

\end{document}